\documentclass[11pt,a4paper]{article}

\usepackage{amsmath,amssymb,latexsym,pstricks,mathrsfs,comment,amsthm,graphicx,tikz,tikz-cd,enumerate,accents,pgffor,cite,wrapfig,multicol,float,slashed}
\usepackage{bibspacing}
\usepackage[T1]{fontenc}
\usepackage{ wasysym }
\usepackage{stmaryrd}
\usepackage[colorlinks]{hyperref}

\newcommand{\nc}{\newcommand}
\nc{\rnc}{\renewcommand}

\usepackage{geometry}  \rnc{\ss}{\smallskip} \nc{\ms}{\medskip} \nc{\bs}{\bigskip} \nc{\nss}{\vspace{-3mm}}

\let\oldproofname=\proofname
\renewcommand{\proofname}{\rm\bf{\oldproofname}}

\makeatletter

\DeclareMathSymbol{\widehatsym}{\mathord}{largesymbols}{"62}
\newcommand\lowerwidehatsym{%
  \text{\smash{\raisebox{-1.3ex}{%
    $\widehatsym$}}}}
\newcommand\fixwidehat[1]{%
  \mathchoice
    {\accentset{\displaystyle\lowerwidehatsym}{#1}}
    {\accentset{\textstyle\lowerwidehatsym}{#1}}
    {\accentset{\scriptstyle\lowerwidehatsym}{#1}}
    {\accentset{\scriptscriptstyle\lowerwidehatsym}{#1}}
}
\rnc{\widehat}{\fixwidehat}

\begin{document}

%\rnc{\phi}{\includegraphics[width=3.5mm]{Om.png}}
\nc{\pfitem}[1]{\medskip \noindent (#1).}
\nc{\Set}{\operatorname{{\bf Set}}}
\nc{\Setp}{\operatorname{{\bf Set}}^+}
\nc{\TXYa}{\T_{XY}^a}
\nc{\IXYa}{\I_{XY}^a}
\nc{\PT}{\P\T}
\rnc{\emptyset}{\varnothing}
\nc{\leqJa}{\leq_{\J^a}}
\nc{\leqRa}{\leq_{\R^a}}
\nc{\leqLa}{\leq_{\L^a}}
\rnc{\th}{\theta}
\nc{\starb}{\oast}
%\nc{\starb}{\mathrel{\includegraphics[width=3mm]{star1.png}}}
%\nc{\starb}{\mathrel{\includegraphics[width=3mm]{star2.png}}}
%\nc{\starb}{\mathrel{\odot}}
%\nc{\starb}{\mathrel{\bigstar}}
%\nc{\starb}{\mathrel{\ast}}
%\nc{\starb}{\mathrel{\includegraphics[width=2.5mm]{starofdavid.jpg}}}
%\nc{\starb}{\mathrel{\includegraphics[width=2.5mm]{starofdavid1.png}}}
%\nc{\starb}{\mathrel{{\lower0.19 ex\hbox{$\includegraphics[width=2.5mm]{starofdavid2.jpg}$}}}}
%\nc{\starb}{\mathrel{{\lower0.15 ex\hbox{$\includegraphics[width=2.5mm]{starofdavid2.jpg}$}}}} % good SOD
\nc{\bPTfr}{\bPT^{\operatorname{fr}}}
\nc{\Sym}{\operatorname{Sym}}
\nc{\PTXYa}{\PT_{XY}^a}
\nc{\PTXY}{\PT_{XY}}
\nc{\TXY}{\T_{XY}}
\nc{\IXY}{\I_{XY}}
\nc{\PTYX}{\PT_{YX}}
\nc{\tran}[2]{\left(\begin{smallmatrix} #1\\#2 \end{smallmatrix}\right)}
\nc{\gt}{\widetilde{g}}
\nc{\cE}{\mathcal E}
\nc{\EXYa}{\cE_{XY}^a}
\nc{\stirlingii}{\genfrac{[}{]}{0pt}{}}
\nc{\sh}{\operatorname{sh}}
\nc{\col}{\operatorname{col}}
\nc{\defect}{\operatorname{def}}
\nc{\codefect}{\operatorname{codef}}
\nc{\sP}{\mathscr P}
\nc{\sE}{\mathscr E}
\nc{\sF}{\mathscr F}
\nc{\JPa}{\J^{P^a}}
\nc{\leqJPa}{\leq_{\JPa}}
\nc{\A}{B} 
\nc{\B}{A} 
\nc{\sectiontitle}[1]{\section{\boldmath #1}}
\nc{\subsectiontitle}[1]{\subsection[#1]{\boldmath #1}}

%\nc{\starb}{\mathrel{\includegraphics[width=3mm]{pentagram.png}}}
%\nc{\starb}{\mbox{\FiveStarOpen}}
\nc{\MaxE}{\operatorname{Max}_\preceq}

\nc{\DClass}[5]{
\begin{tikzpicture}[scale=#5]
\foreach \x/\y in {#3} {\fill[lightgray](\y-1,#2-\x)--(\y,#2-\x)--(\y,#2-\x+1)--(\y-1,#2-\x+1)--(\y-1,#2-\x); \draw(\y-.5,#2-\x+.5)node{$\scriptscriptstyle{#4}$};}
\foreach \x in {0,1,...,#1} {\draw(\x,0)--(\x,#2);}
\foreach \x in {0,1,...,#2} {\draw(0,\x)--(#1,\x);}
\draw[ultra thick] (0,0)--(#1,0)--(#1,#2)--(0,#2)--(0,0)--(#1,0);
\end{tikzpicture}
}

\nc{\DaClass}[7]{
\begin{tikzpicture}[scale=#5]
\foreach \x/\y in {#3} {\fill[lightgray](\y-1,#2-\x)--(\y,#2-\x)--(\y,#2-\x+1)--(\y-1,#2-\x+1)--(\y-1,#2-\x); \draw(\y-.5,#2-\x+.5)node{$\scriptscriptstyle{#4}$};}
\foreach \x in {0,1,...,#1} {\draw(\x,0)--(\x,#2);}
\foreach \x in {0,1,...,#2} {\draw(0,\x)--(#1,\x);}
\foreach \x in {#6} {\draw[ultra thick](\x,0)--(\x,#2);}
\foreach \x in {#7} {\draw[ultra thick](0,\x)--(#1,\x);}
\draw[ultra thick] (0,0)--(#1,0)--(#1,#2)--(0,#2)--(0,0)--(#1,0);
\end{tikzpicture}
}

\nc{\bT}{\operatorname{\bf T}}
\nc{\bPT}{\operatorname{\bf PT}}
\nc{\bI}{\operatorname{\bf I}}

\nc{\M}{\mathcal M}
\nc{\G}{\mathcal G}
\nc{\F}{\mathbb F}
\nc{\MnJ}{\mathcal M_n^J}
\nc{\EnJ}{\mathcal E_n^J}
\nc{\Mat}{\operatorname{Mat}}
\nc{\RegMnJ}{\Reg(\MnJ)}
\nc{\row}{\mathfrak r}
%\nc{\col}{\mathfrak c}
\nc{\Row}{\operatorname{Row}}
\nc{\Col}{\operatorname{Col}}
\nc{\Span}{\operatorname{span}}
\nc{\mat}[4]{\left[\begin{matrix}#1&#2\\#3&#4\end{matrix}\right]}
\nc{\tmattwo}[2]{\left[\begin{smallmatrix}#1\\#2\end{smallmatrix}\right]}
\nc{\tmat}[4]{\left[\begin{smallmatrix}#1&#2\\#3&#4\end{smallmatrix}\right]}
\nc{\ttmat}[4]{{\tiny \left[\begin{smallmatrix}#1&#2\\#3&#4\end{smallmatrix}\right]}}
\nc{\tmatt}[9]{\left[\begin{smallmatrix}#1&#2&#3\\#4&#5&#6\\#7&#8&#9\end{smallmatrix}\right]}
\nc{\ttmatt}[9]{{\tiny \left[\begin{smallmatrix}#1&#2&#3\\#4&#5&#6\\#7&#8&#9\end{smallmatrix}\right]}}
\nc{\MnGn}{\M_n\sm\G_n}
\nc{\MrGr}{\M_r\sm\G_r}
\nc{\qbin}[2]{\left[\begin{matrix}#1\\#2\end{matrix}\right]_q}
\nc{\tqbin}[2]{\left[\begin{smallmatrix}#1\\#2\end{smallmatrix}\right]_q}
\nc{\qbinx}[3]{\left[\begin{matrix}#1\\#2\end{matrix}\right]_{#3}}
\nc{\tqbinx}[3]{\left[\begin{smallmatrix}#1\\#2\end{smallmatrix}\right]_{#3}}
\nc{\MNJ}{\M_nJ}
\nc{\JMN}{J\M_n}
\nc{\RegMNJ}{\Reg(\MNJ)}
\nc{\RegJMN}{\Reg(\JMN)}
\nc{\RegMMNJ}{\Reg(\MMNJ)}
\nc{\RegJMMN}{\Reg(\JMMN)}
\nc{\Wb}{\overline{W}}
\nc{\Xb}{\overline{X}}
\nc{\Yb}{\overline{Y}}
\nc{\Zb}{\overline{Z}}
\nc{\Sib}{\overline{\Si}}
\nc{\Om}{\Omega}
\nc{\Omb}{\overline{\Om}}
\nc{\Gab}{\overline{\Ga}}
\nc{\qfact}[1]{[#1]_q!}
\nc{\smat}[2]{\left[\begin{matrix}#1&#2\end{matrix}\right]}
\nc{\tsmat}[2]{\left[\begin{smallmatrix}#1&#2\end{smallmatrix}\right]}
\nc{\hmat}[2]{\left[\begin{matrix}#1\\#2\end{matrix}\right]}
\nc{\thmat}[2]{\left[\begin{smallmatrix}#1\\#2\end{smallmatrix}\right]}
\nc{\LVW}{\mathcal L(V,W)}
\nc{\KVW}{\mathcal K(V,W)}
\nc{\LV}{\mathcal L(V)}
\nc{\RegLVW}{\Reg(\LVW)}
\nc{\sM}{\mathscr M}
\nc{\sN}{\mathscr N}
\rnc{\iff}{\ \Leftrightarrow\ }
\nc{\Hom}{\operatorname{Hom}}
\nc{\End}{\operatorname{End}}
\nc{\Aut}{\operatorname{Aut}}
\nc{\Lin}{\mathcal L}
\nc{\Hommn}{\Hom(V_m,V_n)}
\nc{\Homnm}{\Hom(V_n,V_m)}
\nc{\Homnl}{\Hom(V_n,V_l)}
\nc{\Homkm}{\Hom(V_k,V_m)}
\nc{\Endm}{\End(V_m)}
\nc{\Endn}{\End(V_n)}
\nc{\Endr}{\End(V_r)}
\nc{\Autm}{\Aut(V_m)}
\nc{\Autn}{\Aut(V_n)}
\nc{\MmnJ}{\M_{mn}^J}
\nc{\MmnA}{\M_{mn}^A}
\nc{\MmnB}{\M_{mn}^B}
\nc{\Mmn}{\M_{mn}}
\nc{\Mkl}{\M_{kl}}
\nc{\Mnm}{\M_{nm}}
\nc{\EmnJ}{\mathcal E_{mn}^J}
\nc{\MmGm}{\M_m\sm\G_m}
\nc{\RegMmnJ}{\Reg(\MmnJ)}
\rnc{\implies}{\ \Rightarrow\ }
\nc{\DMmn}[1]{D_{#1}(\Mmn)}
\nc{\DMmnJ}[1]{D_{#1}(\MmnJ)}
\nc{\MMNJ}{\Mmn J}
\nc{\JMMN}{J\Mmn}
\nc{\JMMNJ}{J\Mmn J}
\nc{\Inr}{\mathcal I(V_n,W_r)}
\nc{\Lnr}{\mathcal L(V_n,W_r)}
\nc{\Knr}{\mathcal K(V_n,W_r)}
\nc{\Imr}{\mathcal I(V_m,W_r)}
\nc{\Kmr}{\mathcal K(V_m,W_r)}
\nc{\Lmr}{\mathcal L(V_m,W_r)}
\nc{\Kmmr}{\mathcal K(V_m,W_{m-r})}
\nc{\tr}{{\operatorname{T}}}
\nc{\MMN}{\MmnA(\F_1)}
\nc{\MKL}{\Mkl^B(\F_2)}
\nc{\RegMMN}{\Reg(\MmnA(\F_1))}
\nc{\RegMKL}{\Reg(\Mkl^B(\F_2))}
\nc{\gRhA}{\widehat{\mathscr R}^A}
\nc{\gRhB}{\widehat{\mathscr R}^B}
\nc{\gLhA}{\widehat{\mathscr L}^A}
\nc{\gLhB}{\widehat{\mathscr L}^B}
\nc{\timplies}{\Rightarrow}
\nc{\tiff}{\Leftrightarrow}
\nc{\Sija}{S_{ij}^a}
\nc{\dmat}[8]{\draw(#1*1.5,#2)node{$\left[\begin{smallmatrix}#3&#4&#5\\#6&#7&#8\end{smallmatrix}\right]$};}
\nc{\bdmat}[8]{\draw(#1*1.5,#2)node{${\mathbf{\left[\begin{smallmatrix}#3&#4&#5\\#6&#7&#8\end{smallmatrix}\right]}}$};}
\nc{\rdmat}[8]{\draw(#1*1.5,#2)node{\rotatebox{90}{$\left[\begin{smallmatrix}#3&#4&#5\\#6&#7&#8\end{smallmatrix}\right]$}};}
\nc{\rldmat}[8]{\draw(#1*1.5-0.375,#2)node{\rotatebox{90}{$\left[\begin{smallmatrix}#3&#4&#5\\#6&#7&#8\end{smallmatrix}\right]$}};}
\nc{\rrdmat}[8]{\draw(#1*1.5+.375,#2)node{\rotatebox{90}{$\left[\begin{smallmatrix}#3&#4&#5\\#6&#7&#8\end{smallmatrix}\right]$}};}
\nc{\rfldmat}[8]{\draw(#1*1.5-0.375+.15,#2)node{\rotatebox{90}{$\left[\begin{smallmatrix}#3&#4&#5\\#6&#7&#8\end{smallmatrix}\right]$}};}
\nc{\rfrdmat}[8]{\draw(#1*1.5+.375-.15,#2)node{\rotatebox{90}{$\left[\begin{smallmatrix}#3&#4&#5\\#6&#7&#8\end{smallmatrix}\right]$}};}
\nc{\xL}{[x]_{\! _\gL}}\nc{\yL}{[y]_{\! _\gL}}\nc{\xR}{[x]_{\! _\gR}}\nc{\yR}{[y]_{\! _\gR}}\nc{\xH}{[x]_{\! _\gH}}\nc{\yH}{[y]_{\! _\gH}}\nc{\XK}{[X]_{\! _\gK}}\nc{\xK}{[x]_{\! _\gK}}
%\nc{\xL}{[x]_\L}\nc{\yL}{[y]_\L}\nc{\xR}{[x]_\R}\nc{\yR}{[y]_\R}\nc{\xH}{[x]_\H}\nc{\yH}{[y]_\H}
\nc{\RegSija}{\Reg(\Sija)}
\nc{\MnmK}{\M_{nm}^K}
\nc{\cC}{\mathcal C}
\nc{\cR}{\mathcal R}
\nc{\Ckl}{\cC_k(l)}
\nc{\Rkl}{\cR_k(l)}
\nc{\Cmr}{\cC_m(r)}
\nc{\Rmr}{\cR_m(r)}
\nc{\Cnr}{\cC_n(r)}
\nc{\Rnr}{\cR_n(r)}
\nc{\Z}{\mathbb Z}

\nc{\Reg}{\operatorname{Reg}}
\nc{\RP}{\operatorname{RP}}
\nc{\MI}{\operatorname{MI}}
\nc{\TXa}{\T_X^a}
\nc{\TXA}{\T(X,A)}
\nc{\TXal}{\T(X,\al)}
\nc{\RegTXa}{\Reg(\TXa)}
\nc{\RegTXA}{\Reg(\TXA)}
\nc{\RegTXal}{\Reg(\TXal)}
\nc{\PalX}{\P_\al(X)}
\nc{\EAX}{\E_A(X)}
\nc{\Bb}{\overline{B}}
\nc{\bb}{\overline{b}}
\nc{\bw}{{\bf w}}
\nc{\bz}{{\bf z}}
\nc{\TASA}{\T_A\sm\S_A}
\nc{\Ub}{\overline{U}}
\nc{\Vb}{\overline{V}}
\nc{\eb}{\overline{e}}
\nc{\EXa}{\E_X^a}
\nc{\oijr}{1\leq i<j\leq r}
\nc{\veb}{\overline{\ve}}
\nc{\bbT}{\mathbb T}
\nc{\Surj}{\operatorname{Surj}}
\nc{\Sone}{S^{(1)}}
\nc{\fillbox}[2]{\draw[fill=gray!30](#1,#2)--(#1+1,#2)--(#1+1,#2+1)--(#1,#2+1)--(#1,#2);}
\nc{\raa}{\rangle_J}
\nc{\raJ}{\rangle_J}
%\rnc{\star}{\star}
\nc{\Ea}{E_J}
\nc{\EJ}{E_J}
\nc{\ep}{\epsilon} \nc{\ve}{\varepsilon}
%\nc{\ep}{\varepsilon} \nc{\ve}{\eta}
\nc{\IXa}{\I_X^a}
\nc{\RegIXa}{\Reg(\IXa)}
\nc{\JXa}{\J_X^a}
\nc{\RegJXa}{\Reg(\JXa)}
\nc{\IXA}{\I(X,A)}
\nc{\IAX}{\I(A,X)}
\nc{\RegIXA}{\Reg(\IXA)}
\nc{\RegIAX}{\Reg(\IAX)}
\nc{\trans}[2]{\left(\begin{smallmatrix} #1 \\ #2 \end{smallmatrix}\right)}
%\nc{\trans}[2]{\big(\begin{smallmatrix} #1 \\ #2 \end{smallmatrix}\big)}
%\nc{\trans}[2]{\Big(\begin{smallmatrix} #1 \\ #2 \end{smallmatrix}\Big)}
\nc{\bigtrans}[2]{\left(\begin{matrix} #1 \\ #2 \end{matrix}\right)}
\nc{\lmap}[1]{\mapstochar \xrightarrow {\ #1\ }}
\nc{\EaTXa}{E}

\nc{\gL}{\mathrel{\mathscr L}}
\nc{\gR}{\mathrel{\mathscr R}}
\nc{\gH}{\mathrel{\mathscr H}}
\nc{\gJ}{\mathrel{\mathscr J}}
\nc{\gD}{\mathrel{\mathscr D}}
\nc{\gK}{\mathrel{\mathscr K}}
\nc{\gLa}{\mathrel{\mathscr L^a}}
\nc{\gRa}{\mathrel{\mathscr R^a}}
\nc{\gHa}{\mathrel{\mathscr H^a}}
\nc{\gJa}{\mathrel{\mathscr J^a}}
\nc{\gDa}{\mathrel{\mathscr D^a}}
\nc{\gKa}{\mathrel{\mathscr K^a}}
\nc{\gLJ}{\mathrel{\mathscr L^J}}
\nc{\gRJ}{\mathrel{\mathscr R^J}}
\nc{\gHJ}{\mathrel{\mathscr H^J}}
\nc{\gJJ}{\mathrel{\mathscr J^J}}
\nc{\gDJ}{\mathrel{\mathscr D^J}}
\nc{\gKJ}{\mathrel{\mathscr K^J}}
\nc{\gLh}{\mathrel{\widehat{\mathscr L}}}
\nc{\gRh}{\mathrel{\widehat{\mathscr R}}}
\nc{\gHh}{\mathrel{\widehat{\mathscr H}}}
\nc{\gJh}{\mathrel{\widehat{\mathscr J}}}
\nc{\gDh}{\mathrel{\widehat{\mathscr D}}}
\nc{\gKh}{\mathrel{\widehat{\mathscr K}}}
\nc{\Lh}{\widehat{L}}
\nc{\Rh}{\widehat{R}}
\nc{\Hh}{\widehat{H}}
\nc{\Jh}{\widehat{J}}
\nc{\Dh}{\widehat{D}}
\nc{\Kh}{\widehat{K}}
\nc{\gLb}{\mathrel{\widehat{\mathscr L}}}
\nc{\gRb}{\mathrel{\widehat{\mathscr R}}}
\nc{\gHb}{\mathrel{\widehat{\mathscr H}}}
\nc{\gJb}{\mathrel{\widehat{\mathscr J}}}
\nc{\gDb}{\mathrel{\widehat{\mathscr D}}}
\nc{\gKb}{\mathrel{\widehat{\mathscr K}}}
\nc{\Lb}{\mathrel{\widehat{L}^J}}
\nc{\Rb}{\mathrel{\widehat{R}^J}}
\nc{\Hb}{\mathrel{\widehat{H}^J}}
\nc{\Jb}{\mathrel{\widehat{J}^J}}
\nc{\Db}{\mathrel{\overline{D}}}
\nc{\Kb}{\mathrel{\widehat{K}}}

%\nc{\ub}{\overline{u}}
%\nc{\vb}{\overline{v}}
%\nc{\wb}{\overline{w}}
\nc{\xb}{\overline{x}}
\nc{\yb}{\overline{y}}
\nc{\zb}{\overline{z}}
%\nc{\pb}{\overline{p}}
\nc{\qb}{\overline{q}}

\hyphenation{mon-oid mon-oids}

\nc{\itemit}[1]{\item[\emph{(#1)}]}
\nc{\itemnit}[1]{\item[(#1)]}
\nc{\E}{\mathbb E}
\nc{\TX}{\T(X)}
\nc{\TXP}{\T(X,\P)}
\nc{\EX}{\E(X)}
\nc{\EXP}{\E(X,\P)}
\nc{\SX}{\S(X)}
\nc{\SXP}{\S(X,\P)}
\nc{\Sing}{\operatorname{Sing}}
%\nc{\Sing}{\E}
\nc{\idrank}{\operatorname{idrank}}
\nc{\SingXP}{\Sing(X,\P)}
\nc{\De}{\Delta}
\nc{\sgp}{\operatorname{sgp}}
\nc{\mon}{\operatorname{mon}}
\nc{\Dn}{\mathcal D_n}
\nc{\Dm}{\mathcal D_m}

\nc{\lline}[1]{\draw(3*#1,0)--(3*#1+2,0);}
\nc{\uline}[1]{\draw(3*#1,5)--(3*#1+2,5);}
\nc{\thickline}[2]{\draw(3*#1,5)--(3*#2,0); \draw(3*#1+2,5)--(3*#2+2,0) ;}
\nc{\thicklabel}[3]{\draw(3*#1+1+3*#2*0.15-3*#1*0.15,4.25)node{{\tiny $#3$}};}

\nc{\slline}[3]{\draw(3*#1+#3,0+#2)--(3*#1+2+#3,0+#2);}
\nc{\suline}[3]{\draw(3*#1+#3,5+#2)--(3*#1+2+#3,5+#2);}
\nc{\sthickline}[4]{\draw(3*#1+#4,5+#3)--(3*#2+#4,0+#3); \draw(3*#1+2+#4,5+#3)--(3*#2+2+#4,0+#3) ;}
\nc{\sthicklabel}[5]{\draw(3*#1+1+3*#2*0.15-3*#1*0.15+#5,4.25+#4)node{{\tiny $#3$}};}

\nc{\stll}[5]{\sthickline{#1}{#2}{#4}{#5} \sthicklabel{#1}{#2}{#3}{#4}{#5}}
\nc{\tll}[3]{\stll{#1}{#2}{#3}00}

\nc{\mfourpic}[9]{
\slline1{#9}0
\slline3{#9}0
\slline4{#9}0
\slline5{#9}0
\suline1{#9}0
\suline3{#9}0
\suline4{#9}0
\suline5{#9}0
\stll1{#1}{#5}{#9}{0}
\stll3{#2}{#6}{#9}{0}
\stll4{#3}{#7}{#9}{0}
\stll5{#4}{#8}{#9}{0}
\draw[dotted](6,0+#9)--(8,0+#9);
\draw[dotted](6,5+#9)--(8,5+#9);
}
\nc{\vdotted}[1]{
\draw[dotted](3*#1,10)--(3*#1,15);
\draw[dotted](3*#1+2,10)--(3*#1+2,15);
}

\nc{\Clab}[2]{
\sthicklabel{#1}{#1}{{}_{\phantom{#1}}C_{#1}}{1.25+5*#2}0
}
\nc{\sClab}[3]{
\sthicklabel{#1}{#1}{{}_{\phantom{#1}}C_{#1}}{1.25+5*#2}{#3}
}
\nc{\Clabl}[3]{
\sthicklabel{#1}{#1}{{}_{\phantom{#3}}C_{#3}}{1.25+5*#2}0
}
\nc{\sClabl}[4]{
\sthicklabel{#1}{#1}{{}_{\phantom{#4}}C_{#4}}{1.25+5*#2}{#3}
}
\nc{\Clabll}[3]{
\sthicklabel{#1}{#1}{C_{#3}}{1.25+5*#2}0
}
\nc{\sClabll}[4]{
\sthicklabel{#1}{#1}{C_{#3}}{1.25+5*#2}{#3}
}

\nc{\mtwopic}[6]{
\slline1{#6*5}{#5}
\slline2{#6*5}{#5}
\suline1{#6*5}{#5}
\suline2{#6*5}{#5}
\stll1{#1}{#3}{#6*5}{#5}
\stll2{#2}{#4}{#6*5}{#5}
}
\nc{\mtwopicl}[6]{
\slline1{#6*5}{#5}
\slline2{#6*5}{#5}
\suline1{#6*5}{#5}
\suline2{#6*5}{#5}
\stll1{#1}{#3}{#6*5}{#5}
\stll2{#2}{#4}{#6*5}{#5}
\sClabl1{#6}{#5}{i}
\sClabl2{#6}{#5}{j}
}

%\nc{\keru}{\ker_{\operatorname{u}}} \nc{\kerl}{\ker_{\operatorname{l}}}
%\nc{\dom}{\operatorname{dom}_u} \nc{\codom}{\operatorname{dom}_l}
\nc{\keru}{\operatorname{ker}^\wedge} \nc{\kerl}{\operatorname{ker}_\vee}%\nc{\codom}{\operatorname{im}}

\nc{\coker}{\operatorname{coker}}
%\nc{\KER}{\operatorname{KER}}
\nc{\KER}{\ker}
\nc{\N}{\mathbb N}
\nc{\LaBn}{L_\al(\B_n)}
\nc{\RaBn}{R_\al(\B_n)}
\nc{\LaPBn}{L_\al(\PB_n)}
\nc{\RaPBn}{R_\al(\PB_n)}
\nc{\rhorBn}{\rho_r(\B_n)}
\nc{\DrBn}{D_r(\B_n)}
\nc{\DrPn}{D_r(\P_n)}
\nc{\DrPBn}{D_r(\PB_n)}
\nc{\DrKn}{D_r(\K_n)}
\nc{\alb}{\al_{\vee}}
\nc{\beb}{\be^{\wedge}}
\nc{\bnf}{\bn^\flat}
\nc{\Bal}{\operatorname{Bal}}
\nc{\Red}{\operatorname{Red}}
\nc{\Pnxi}{\P_n^\xi}
\nc{\Bnxi}{\B_n^\xi}
\nc{\PBnxi}{\PB_n^\xi}
\nc{\Knxi}{\K_n^\xi}
\nc{\C}{\mathscr C}
\nc{\exi}{e^\xi}
\nc{\Exi}{E^\xi}
\nc{\eximu}{e^\xi_\mu}
\nc{\Eximu}{E^\xi_\mu}
\nc{\REF}{ {\red [Ref?]} }
\nc{\GL}{\operatorname{GL}}
\rnc{\O}{\operatorname{O}}

\nc{\vtx}[2]{\fill (#1,#2)circle(.2);}
\nc{\lvtx}[2]{\fill (#1,0)circle(.2);}
\nc{\uvtx}[2]{\fill (#1,1.5)circle(.2);}

\nc{\Eq}{\mathfrak{Eq}}
%\nc{\Gau}{\Ga_{\operatorname{u}}} \nc{\Gal}{\Ga_{\operatorname{l}}}
\nc{\Gau}{\Ga^\wedge} \nc{\Gal}{\Ga_\vee}
%\nc{\Lamu}{\Lam_{\operatorname{u}}} \nc{\Laml}{\Lam_{\operatorname{l}}}
\nc{\Lamu}{\Lam^\wedge} \nc{\Laml}{\Lam_\vee}
\nc{\bX}{{\bf X}}
\nc{\bY}{{\bf Y}}
\nc{\ds}{\displaystyle}

\nc{\uvert}[1]{\fill (#1,1.5)circle(.2);}
\nc{\uuvert}[1]{\fill (#1,3)circle(.2);}
\nc{\uuuvert}[1]{\fill (#1,4.5)circle(.2);}
\rnc{\lvert}[1]{\fill (#1,0)circle(.2);}
%\nc{\overt}[1]{\fill (#1,0)circle(.1);}
\nc{\overtl}[3]{\node[vertex] (#3) at (#1,0) {  {\tiny $#2$} };}
\nc{\cv}[2]{\draw(#1,1.5) to [out=270,in=90] (#2,0);}
\nc{\cvs}[2]{\draw(#1,1.5) to [out=270+30,in=90+30] (#2,0);}
\nc{\ucv}[2]{\draw(#1,3) to [out=270,in=90] (#2,1.5);}
\nc{\uucv}[2]{\draw(#1,4.5) to [out=270,in=90] (#2,3);}
\nc{\textpartn}[1]{{\lower0.45 ex\hbox{\begin{tikzpicture}[xscale=.2,yscale=0.2] #1 \end{tikzpicture}}}}
\nc{\textpartnx}[2]{{\lower1.0 ex\hbox{\begin{tikzpicture}[xscale=.3,yscale=0.3] 
\foreach \x in {1,...,#1}
{ \uvert{\x} \lvert{\x} }
#2 \end{tikzpicture}}}}
\nc{\disppartnx}[2]{{\lower1.0 ex\hbox{\begin{tikzpicture}[scale=0.3] 
\foreach \x in {1,...,#1}
{ \uvert{\x} \lvert{\x} }
#2 \end{tikzpicture}}}}
\nc{\disppartnxd}[2]{{\lower2.1 ex\hbox{\begin{tikzpicture}[scale=0.3] 
\foreach \x in {1,...,#1}
{ \uuvert{\x} \uvert{\x} \lvert{\x} }
#2 \end{tikzpicture}}}}
\nc{\disppartnxdn}[2]{{\lower2.1 ex\hbox{\begin{tikzpicture}[scale=0.3] 
\foreach \x in {1,...,#1}
{ \uuvert{\x} \lvert{\x} }
#2 \end{tikzpicture}}}}
%\nc{\disppartnxdd}[2]{{\lower2.1 ex\hbox{\begin{tikzpicture}[scale=0.3] 
%\foreach \x in {1,...,#1}
%{ \uuuvert{\x} \uuvert{\x} \uvert{\x} \lvert{\x} }
%#2 \end{tikzpicture}}}}
\nc{\disppartnxdd}[2]{{\lower3.6 ex\hbox{\begin{tikzpicture}[scale=0.3] 
\foreach \x in {1,...,#1}
{ \uuuvert{\x} \uuvert{\x} \uvert{\x} \lvert{\x} }
#2 \end{tikzpicture}}}}

\nc{\dispgax}[2]{{\lower0.0 ex\hbox{\begin{tikzpicture}[scale=0.3] 
#2
\foreach \x in {1,...,#1}
{\lvert{\x} }
 \end{tikzpicture}}}}
\nc{\textgax}[2]{{\lower0.4 ex\hbox{\begin{tikzpicture}[scale=0.3] 
#2
\foreach \x in {1,...,#1}
{\lvert{\x} }
 \end{tikzpicture}}}}
\nc{\textlinegraph}[2]{{\raise#1 ex\hbox{\begin{tikzpicture}[scale=0.8] 
#2
 \end{tikzpicture}}}}
\nc{\textlinegraphl}[2]{{\raise#1 ex\hbox{\begin{tikzpicture}[scale=0.8] 
\tikzstyle{vertex}=[circle,draw=black, fill=white, inner sep = 0.07cm]
#2
 \end{tikzpicture}}}}
\nc{\displinegraph}[1]{{\lower0.0 ex\hbox{\begin{tikzpicture}[scale=0.6] 
#1
 \end{tikzpicture}}}}
 
\nc{\disppartnthreeone}[1]{{\lower1.0 ex\hbox{\begin{tikzpicture}[scale=0.3] 
\foreach \x in {1,2,3,5,6}
{ \uvert{\x} }
\foreach \x in {1,2,4,5,6}
{ \lvert{\x} }
\draw[dotted] (3.5,1.5)--(4.5,1.5);
\draw[dotted] (2.5,0)--(3.5,0);
#1 \end{tikzpicture}}}}

\nc{\partn}[4]{\left( \begin{array}{c|c} %fine
#1 \ & \ #3 \ \ \\ \cline{2-2}
#2 \ & \ #4 \ \
\end{array} \!\!\! \right)}
\nc{\partnlong}[6]{\partn{#1}{#2}{#3,\ #4}{#5,\ #6}} %fine
\nc{\partnsh}[2]{\left( \begin{array}{c} %fine
#1 \\
#2 
\end{array} \right)}
\nc{\partncodefz}[3]{\partn{#1}{#2}{#3}{\emptyset}}
\nc{\partndefz}[3]{{\partn{#1}{#2}{\emptyset}{#3}}}
\nc{\partnlast}[2]{\left( \begin{array}{c|c}
#1 \ &  \ #2 \\
#1 \ &  \ #2
\end{array} \right)}

\nc{\uuarcx}[3]{\draw(#1,3)arc(180:270:#3) (#1+#3,3-#3)--(#2-#3,3-#3) (#2-#3,3-#3) arc(270:360:#3);}
\nc{\uuarc}[2]{\uuarcx{#1}{#2}{.4}}
\nc{\uuuarcx}[3]{\draw(#1,4.5)arc(180:270:#3) (#1+#3,4.5-#3)--(#2-#3,4.5-#3) (#2-#3,4.5-#3) arc(270:360:#3);}
\nc{\uuuarc}[2]{\uuuarcx{#1}{#2}{.4}}
\nc{\darcx}[3]{\draw(#1,0)arc(180:90:#3) (#1+#3,#3)--(#2-#3,#3) (#2-#3,#3) arc(90:0:#3);}
\nc{\darc}[2]{\darcx{#1}{#2}{.4}}
\nc{\udarcx}[3]{\draw(#1,1.5)arc(180:90:#3) (#1+#3,1.5+#3)--(#2-#3,1.5+#3) (#2-#3,1.5+#3) arc(90:0:#3);}
\nc{\udarc}[2]{\udarcx{#1}{#2}{.4}}
\nc{\uudarcx}[3]{\draw(#1,3)arc(180:90:#3) (#1+#3,3+#3)--(#2-#3,3+#3) (#2-#3,3+#3) arc(90:0:#3);}
\nc{\uudarc}[2]{\uudarcx{#1}{#2}{.4}}
\nc{\uarcx}[3]{\draw(#1,1.5)arc(180:270:#3) (#1+#3,1.5-#3)--(#2-#3,1.5-#3) (#2-#3,1.5-#3) arc(270:360:#3);}
\nc{\uarc}[2]{\uarcx{#1}{#2}{.4}}
\nc{\darcxhalf}[3]{\draw(#1,0)arc(180:90:#3) (#1+#3,#3)--(#2,#3) ;}
\nc{\darchalf}[2]{\darcxhalf{#1}{#2}{.4}}
\nc{\uarcxhalf}[3]{\draw(#1,1.5)arc(180:270:#3) (#1+#3,1.5-#3)--(#2,1.5-#3) ;}
\nc{\uarchalf}[2]{\uarcxhalf{#1}{#2}{.4}}
\nc{\uarcxhalfr}[3]{\draw (#1+#3,1.5-#3)--(#2-#3,1.5-#3) (#2-#3,1.5-#3) arc(270:360:#3);}
\nc{\uarchalfr}[2]{\uarcxhalfr{#1}{#2}{.4}}

\nc{\bdarcx}[3]{\draw[blue](#1,0)arc(180:90:#3) (#1+#3,#3)--(#2-#3,#3) (#2-#3,#3) arc(90:0:#3);}
\nc{\bdarc}[2]{\darcx{#1}{#2}{.4}}
\nc{\rduarcx}[3]{\draw[red](#1,0)arc(180:270:#3) (#1+#3,0-#3)--(#2-#3,0-#3) (#2-#3,0-#3) arc(270:360:#3);}
\nc{\rduarc}[2]{\uarcx{#1}{#2}{.4}}
\nc{\duarcx}[3]{\draw(#1,0)arc(180:270:#3) (#1+#3,0-#3)--(#2-#3,0-#3) (#2-#3,0-#3) arc(270:360:#3);}
\nc{\duarc}[2]{\uarcx{#1}{#2}{.4}}

\nc{\uv}[1]{\fill (#1,2)circle(.1);}
\nc{\lv}[1]{\fill (#1,0)circle(.1);}
\nc{\stline}[2]{\draw(#1,2)--(#2,0);}
\nc{\tlab}[2]{\draw(#1,2)node[above]{\tiny $#2$};}
\nc{\tudots}[1]{\draw(#1,2)node{$\cdots$};}
\nc{\tldots}[1]{\draw(#1,0)node{$\cdots$};}

\nc{\uvw}[1]{\fill[white] (#1,2)circle(.1);}
\nc{\huv}[1]{\fill (#1,1)circle(.1);}
\nc{\llv}[1]{\fill (#1,-2)circle(.1);}
\nc{\arcup}[2]{
\draw(#1,2)arc(180:270:.4) (#1+.4,1.6)--(#2-.4,1.6) (#2-.4,1.6) arc(270:360:.4);
}
\nc{\harcup}[2]{
\draw(#1,1)arc(180:270:.4) (#1+.4,.6)--(#2-.4,.6) (#2-.4,.6) arc(270:360:.4);
}
\nc{\arcdn}[2]{
\draw(#1,0)arc(180:90:.4) (#1+.4,.4)--(#2-.4,.4) (#2-.4,.4) arc(90:0:.4);
}
\nc{\cve}[2]{
\draw(#1,2) to [out=270,in=90] (#2,0);
}
\nc{\hcve}[2]{
\draw(#1,1) to [out=270,in=90] (#2,0);
}
\nc{\catarc}[3]{
\draw(#1,2)arc(180:270:#3) (#1+#3,2-#3)--(#2-#3,2-#3) (#2-#3,2-#3) arc(270:360:#3);
}

\nc{\arcr}[2]{
\draw[red](#1,0)arc(180:90:.4) (#1+.4,.4)--(#2-.4,.4) (#2-.4,.4) arc(90:0:.4);
}
\nc{\arcb}[2]{
\draw[blue](#1,2-2)arc(180:270:.4) (#1+.4,1.6-2)--(#2-.4,1.6-2) (#2-.4,1.6-2) arc(270:360:.4);
}
\nc{\loopr}[1]{
\draw[blue](#1,-2) edge [out=130,in=50,loop] ();
}
\nc{\loopb}[1]{
\draw[red](#1,-2) edge [out=180+130,in=180+50,loop] ();
}
%\nc{\arcr}[2]{
%\draw[red](#1,0-2)arc(180:90:.4) (#1+.4,.4-2)--(#2-.4,.4-2) (#2-.4,.4-2) arc(90:0:.4);
%}
%\nc{\arcb}[2]{
%\draw[blue](#1,2-2-2)arc(180:270:.4) (#1+.4,1.6-2-2)--(#2-.4,1.6-2-2) (#2-.4,1.6-2-2) arc(270:360:.4);
%}
%\nc{\loopr}[1]{
%\draw[red](#1,0-2) edge [out=130,in=50,loop] ();
%}
%\nc{\loopb}[1]{
%\draw[blue](#1,0-2) edge [out=180+130,in=180+50,loop] ();
%}
\nc{\redto}[2]{\draw[red](#1,0)--(#2,0);}
\nc{\bluto}[2]{\draw[blue](#1,0)--(#2,0);}
\nc{\dotto}[2]{\draw[dotted](#1,0)--(#2,0);}
\nc{\lloopr}[2]{\draw[red](#1,0)arc(0:360:#2);}
\nc{\lloopb}[2]{\draw[blue](#1,0)arc(0:360:#2);}
\nc{\rloopr}[2]{\draw[red](#1,0)arc(-180:180:#2);}
\nc{\rloopb}[2]{\draw[blue](#1,0)arc(-180:180:#2);}
\nc{\uloopr}[2]{\draw[red](#1,0)arc(-270:270:#2);}
\nc{\uloopb}[2]{\draw[blue](#1,0)arc(-270:270:#2);}
\nc{\dloopr}[2]{\draw[red](#1,0)arc(-90:270:#2);}
\nc{\dloopb}[2]{\draw[blue](#1,0)arc(-90:270:#2);}
\nc{\llloopr}[2]{\draw[red](#1,0-2)arc(0:360:#2);}
\nc{\llloopb}[2]{\draw[blue](#1,0-2)arc(0:360:#2);}
\nc{\lrloopr}[2]{\draw[red](#1,0-2)arc(-180:180:#2);}
\nc{\lrloopb}[2]{\draw[blue](#1,0-2)arc(-180:180:#2);}
\nc{\ldloopr}[2]{\draw[red](#1,0-2)arc(-270:270:#2);}
\nc{\ldloopb}[2]{\draw[blue](#1,0-2)arc(-270:270:#2);}
\nc{\luloopr}[2]{\draw[red](#1,0-2)arc(-90:270:#2);}
\nc{\luloopb}[2]{\draw[blue](#1,0-2)arc(-90:270:#2);}

\nc{\larcb}[2]{
\draw[blue](#1,0-2)arc(180:90:.4) (#1+.4,.4-2)--(#2-.4,.4-2) (#2-.4,.4-2) arc(90:0:.4);
}
\nc{\larcr}[2]{
\draw[red](#1,2-2-2)arc(180:270:.4) (#1+.4,1.6-2-2)--(#2-.4,1.6-2-2) (#2-.4,1.6-2-2) arc(270:360:.4);
}

\rnc{\H}{\mathrel{\mathscr H}}
\rnc{\L}{\mathrel{\mathscr L}}
\nc{\R}{\mathrel{\mathscr R}}
\nc{\D}{\mathrel{\mathscr D}}
\nc{\J}{\mathrel{\mathscr J}}
\nc{\leqR}{\mathrel{\leq_{\R}}}
\nc{\leqL}{\mathrel{\leq_{\L}}}
\nc{\leqJ}{\mathrel{\leq_{\J}}}

\nc{\ssim}{\mathrel{\raise0.25 ex\hbox{\oalign{$\approx$\crcr\noalign{\kern-0.84 ex}$\approx$}}}}
\nc{\POI}{\mathcal{POI}}
\nc{\wb}{\overline{w}}
\nc{\ub}{\overline{u}}
\nc{\vb}{\overline{v}}
\nc{\fb}{\overline{f}}
\nc{\gb}{\overline{g}}
\nc{\hb}{\overline{h}}
\nc{\pb}{\overline{p}}
\rnc{\sb}{\overline{s}}
\nc{\XR}{\pres{X}{R\,}}
\nc{\YQ}{\pres{Y}{Q}}
\nc{\ZP}{\pres{Z}{P\,}}
\nc{\XRone}{\pres{X_1}{R_1}}
\nc{\XRtwo}{\pres{X_2}{R_2}}
\nc{\XRthree}{\pres{X_1\cup X_2}{R_1\cup R_2\cup R_3}}
\nc{\er}{\eqref}
\nc{\larr}{\mathrel{\hspace{-0.35 ex}>\hspace{-1.1ex}-}\hspace{-0.35 ex}}
\nc{\rarr}{\mathrel{\hspace{-0.35 ex}-\hspace{-0.5ex}-\hspace{-2.3ex}>\hspace{-0.35 ex}}}
\nc{\lrarr}{\mathrel{\hspace{-0.35 ex}>\hspace{-1.1ex}-\hspace{-0.5ex}-\hspace{-2.3ex}>\hspace{-0.35 ex}}}
\nc{\nn}{\tag*{}}
\nc{\epfal}{\tag*{$\Box$}}
\nc{\tagd}[1]{\tag*{(#1)$'$}}
\nc{\ldb}{[\![}
\nc{\rdb}{]\!]}
\nc{\sm}{\setminus}
\nc{\I}{\mathcal I}
\nc{\InSn}{\I_n\setminus\S_n}
%\nc{\dom}{\operatorname{dom}_{\operatorname{u}}} \nc{\codom}{\operatorname{dom}_{\operatorname{l}}}
%\nc{\dom}{\operatorname{dom}_u} \nc{\codom}{\operatorname{dom}_l}
\nc{\dom}{\operatorname{dom}} \nc{\codom}{\operatorname{dom}}%\nc{\codom}{\operatorname{im}}
\nc{\ojin}{1\leq j<i\leq n}
%\nc{\R}{\mathcal R}
%\rnc{\L}{\mathcal L}
\nc{\eh}{\widehat{e}}
\nc{\wh}{\widehat{w}}
\nc{\uh}{\widehat{u}}
\nc{\vh}{\widehat{v}}
%\nc{\sh}{\widehat{s}}
\nc{\fh}{\widehat{f}}
\nc{\textres}[1]{\text{\emph{#1}}}
\nc{\aand}{\emph{\ and \ }}
\nc{\iif}{\emph{\ if \ }}
\nc{\textlarr}{\ \larr\ }
\nc{\textrarr}{\ \rarr\ }
\nc{\textlrarr}{\ \lrarr\ }

\nc{\comma}{,\ }

\nc{\COMMA}{,\qquad}
\nc{\COMMa}{,\quad }
\nc{\COMma}{,\ \ \ }
\nc{\COmma}{,\ \ }
\nc{\TnSn}{\T_n\setminus\S_n} 
\nc{\TmSm}{\T_m\setminus\S_m} 
\nc{\TXSX}{\T_X\setminus\S_X} 
\rnc{\S}{\mathcal S}

\nc{\T}{\mathcal T} 
\rnc{\P}{\mathcal P} 
\nc{\K}{\mathrel\mathscr K}
\nc{\PB}{\mathcal{PB}} 
\nc{\rank}{\operatorname{rank}}

\nc{\mtt}{\!\!\!\mt\!\!\!}

\nc{\sub}{\subseteq}
\nc{\la}{\langle}
\nc{\ra}{\rangle}
\nc{\mt}{\mapsto}
\nc{\im}{\mathrm{im}}
\nc{\id}{\mathrm{id}}
\nc{\bn}{\mathbf{n}}
\nc{\ba}{\mathbf{a}}
\nc{\bl}{\mathbf{l}}
\nc{\bm}{\mathbf{m}}
\nc{\bk}{\mathbf{k}}
\nc{\br}{\mathbf{r}}
\nc{\al}{\alpha}
\nc{\be}{\beta}
\nc{\ga}{\gamma}
\nc{\Ga}{\Gamma}
\nc{\de}{\delta}
\nc{\ka}{\kappa}
\nc{\lam}{\lambda}
\nc{\Lam}{\Lambda}
\nc{\si}{\sigma}
\nc{\Si}{\Sigma}
\nc{\oijn}{1\leq i<j\leq n}
\nc{\oijm}{1\leq i<j\leq m}

\nc{\comm}{\rightleftharpoons}
\nc{\AND}{\qquad\text{and}\qquad}
\nc{\ANd}{\quad\text{and}\quad}

\nc{\bit}{\begin{itemize}}
\nc{\bitbmc}{\begin{itemize}\begin{multicols}}
\nc{\bmc}{\begin{itemize}\begin{multicols}}
\nc{\emc}{\end{multicols}\end{itemize}}
\nc{\eit}{\end{itemize}}
\nc{\ben}{\begin{enumerate}}
\nc{\een}{\end{enumerate}}
\nc{\eitres}{\end{itemize}}
\nc{\eitmc}{\end{itemize}}

\nc{\set}[2]{\{ {#1} : {#2} \}} 
\nc{\bigset}[2]{\big\{ {#1}: {#2} \big\}} 
\nc{\Bigset}[2]{\left\{ \,{#1} :{#2}\, \right\}}

\nc{\pres}[2]{\la {#1} \,|\, {#2} \ra}
\nc{\bigpres}[2]{\big\la {#1} \,\big|\, {#2} \big\ra}
\nc{\Bigpres}[2]{\Big\la \,{#1}\, \,\Big|\, \,{#2}\, \Big\ra}
\nc{\Biggpres}[2]{\Bigg\la {#1} \,\Bigg|\, {#2} \Bigg\ra}

\newcommand{\pf}{\begin{proof}}
\newcommand{\epf}{\end{proof}}
%\nc{\pf}{\noindent{\bf Proof.}  }
%\nc{\epf}{\hfill$\Box$\bigskip}
\nc{\epfres}{\qed}
\nc{\pfnb}{\pf}
\nc{\epfnb}{\bigskip}
\nc{\pfthm}[1]{\bigskip \noindent{\bf Proof of Theorem \ref{#1}}\,\,  } 
\nc{\pfprop}[1]{\bigskip \noindent{\bf Proof of Proposition \ref{#1}}\,\,  } 
%\nc{\pfthm}{\noindent{\bf Proof of Theorem \ref{mainthm} modulo Propositions \ref{prop1} and \ref{prop2}}\,\,  } 
\nc{\epfreseq}{\tag*{$\Box$}}

\makeatletter
\newcommand\footnoteref[1]{\protected@xdef\@thefnmark{\ref{#1}}\@footnotemark}
\makeatother

\numberwithin{equation}{section}

\newtheorem{thm}[equation]{Theorem}
\newtheorem{lemma}[equation]{Lemma}
\newtheorem{cor}[equation]{Corollary}
\newtheorem{prop}[equation]{Proposition}

\theoremstyle{definition}

\newtheorem{rem}[equation]{Remark}
\newtheorem{defn}[equation]{Definition}
\newtheorem{eg}[equation]{Example}
\newtheorem{ass}[equation]{Assumption}

\title{Sandwich semigroups in locally small categories II: Transformations}
%\title{Mmmmm... Sandwiches!}
\author{
Igor Dolinka\footnote{Department of Mathematics and Informatics, University of Novi Sad, Novi Sad, Serbia. {\it Emails:} {\tt dockie\,@\,dmi.uns.ac.rs}, {\tt ivana.djurdjev\,@\,dmi.uns.ac.rs}}, 
Ivana \DJ ur\dj ev,\hspace{-.2em}${}^*$\footnote{Mathematical Institute of the Serbian Academy of Sciences and Arts, Beograd, Serbia. {\it Email:} {\tt djurdjev.ivana\,@\,gmail.com}}\;
James East\footnote{Centre for Research in Mathematics, School of Computing, Engineering and Mathematics, Western Sydney University, Sydney, Australia. {\it Email:} {\tt j.east\,@\,westernsydney.edu.au}}, 
Preeyanuch Honyam\footnote{Department of Mathematics, Chiang Mai University, Chiang Mai, Thailand. {\it Emails:} {\tt preeyanuch.h\,@\,cmu.ac.th}, {\tt kritsada.s\,@\,cmu.ac.th}, {\tt jintana.s\,@\,cmu.ac.th}}, \\ 
Kritsada Sangkhanan,\hspace{-.2em}${}^\mathsection$
Jintana Sanwong,\hspace{-.2em}${}^\mathsection$
Worachead Sommanee\footnote{Department of Mathematics and Statistics, Chiang Mai Rajabhat University, Chiang Mai, Thailand. {\it Email:} {\tt worachead\_som\,@\,cmru.ac.th}}
}

\date{}

\maketitle

\vspace{-0.5cm}

\begin{abstract}
%In the prequel, we developed a general theory of sandwich semigroups in arbitrary locally small categories.  Here we apply this to several concrete categories of transformations and partial transformations.  

Fix sets $X$ and $Y$, and write $\PTXY$ for the set of all partial functions $X\to Y$.  Fix a partial function~${a:Y\to X}$, and define the operation $\star_a$ on $\PTXY$ by $f\star_ag=fag$ for $f,g\in\PTXY$.  The \emph{sandwich semigroup} $(\PTXY,\star_a)$ is denoted $\PTXYa$.  We apply general results from Part I to thoroughly describe the structural and combinatorial properties of $\PTXYa$, as well as its regular and idempotent-generated subsemigroups, $\Reg(\PTXYa)$ and $\mathbb E(\PTXYa)$.  After describing regularity, stability and Green's relations and preorders, we exhibit $\Reg(\PTXYa)$ as a pullback product of certain regular subsemigroups of the (non-sandwich) partial transformation semigroups $\PT_X$ and $\PT_Y$, and as a kind of ``inflation'' of $\PT_A$, where~$A$ is the image of the sandwich element $a$.  We also calculate the rank (minimal size of a generating set) and, where appropriate, the idempotent rank (minimal size of an idempotent generating set) of~$\PTXYa$,~$\Reg(\PTXYa)$ and $\mathbb E(\PTXYa)$.  The same program is also carried for sandwich semigroups of totally defined functions and for injective partial functions.  Several corollaries are obtained for various (non-sandwich) semigroups of (partial) transformations with restricted image, domain and/or kernel.

%We develop a theory of sandwich semigroups in arbitrary locally small categories (Part I), and apply this to several concrete categories of transformations and partial transformations (Part II).

%We develop a theory of sandwich semigroups in arbitrary (locally small) categories.
%%, and apply this to several concrete categories of transformations and partial transformations in the sequel.
%
%We apply the general framework of Part I to several concrete categories of transformations and partial transformations.

{\it Keywords}: Categories, sandwich semigroups, transformation semigroups, rank, idempotent rank.

MSC: 20M50; 18B40; 20M10; 20M17; 20M20; 05E15.
\end{abstract}

\section*{Introduction}\label{sect:intro}

%This article is the second in a series on sandwich semigroups in arbitrary (locally small) categories.  The previous article \cite{Sandwiches1} 

This is a continuation of the article \cite{Sandwiches1}, in which we developed a theory of sandwich semigroups in arbitrary (locally small) categories.  Here we apply the general results of \cite{Sandwiches1} to three concrete categories of (partial) transformations.  These are the categories $\PT$, $\T$ and $\I$, consisting of all partial transformations, all full transformations, and all injective partial transformations, respectively.  Endomorphism monoids in these categories are the partial transformation semigroups $\PT_X$, the full transformation semigroups $\T_X$, and the symmetric inverse monoids $\I_X$, respectively.  These monoids have had a profound impact on semigroup theory from its very beginning, and the literature on them is enormous; see \cite{GMbook} for a recent monograph concerning the finite case.  Sandwich semigroups in the categories $\PT$, $\T$ and $\I$ have been studied in various contexts, by a number of authors.  The original studies were carried out in a series of articles by Magill in the 1960s and 1970s; see for example \cite{MS1975,Magill1967}.  A recent study may be found in \cite{MGS2013}, and the introductions of \cite{DElinear,DEvariants,Sandwiches1} may be consulted for more references, historical details and discussion on the connections to other areas of mathematics.  As significant applications of our results, we note in the current article (see also \cite[Section~2.1]{Sandwiches1}) that several other (non-sandwich) semigroups of transformations arise as special cases of the sandwich semigroup construction: specifically, semigroups of (partial) transformations with restricted range, domain or kernel, as studied in \cite{SS2013,FS2014,MGS2010,Symons1975}, for example.  Thus, our general results on sandwich semigroups have natural corollaries in all of these semigroups, as we note throughout.

The article is organised as follows.  We begin with the basic definitions in Section \ref{sect:prelim_transformations}.  Sections \ref{sect:PT}, \ref{sect:T} and \ref{sect:I} then cover the categories $\PT$, $\T$ and $\I$, respectively, and their sandwich semigroups; each of these sections has its own introduction, in which a summary of the main results may be found.  Section \ref{sect:PT} gives a thorough treatment of the category $\PT$, while Sections \ref{sect:T} and \ref{sect:I} are comparatively brief; generally, the main results on the categories $\T$ and $\I$ are stated, with indications of how the arguments of Section \ref{sect:PT} may be adapted to prove them, and comments are provided on similarities and differences.  Throughout the article, we will also obtain a number of corollaries regarding the above-mentioned families of (non-sandwich) semigroups of (partial) transformations with restricted range, domain or kernel.  We work in standard ZFC set theory; see for example \cite[Chapters 1, 5 and 6]{Jech2003}.    The reader may refer to \cite{MacLane1998} for basics on category theory.  We keep the notation of \cite{Sandwiches1}.

%As the current article is a true continuation of \cite{Sandwiches1}, we will freely use the notation and results of the former article throughout.

%\newpage\tableofcontents

\sectiontitle{Preliminaries}\label{sect:prelim_transformations}

We denote by $\Set$ the class of all sets.  For sets $A,B\in\Set$, we write
\begin{align*}
\bT_{AB} &= \set{f}{f\text{ is a function $A\to B$}},\\
\bPT_{AB} &= \set{f}{f\text{ is a function $C\to B$ for some $C\sub A$}},\\
\bI_{AB} &= \set{f}{f\text{ is an injective function $C\to B$ for some $C\sub A$}}.
%\bI_{AB} &= \set{f\in\bPT_{AB}}{f\text{ is injective}}.
\end{align*}
Note that we regard an element of $\bPT_{AB}$ as a subset of $A\times B$ (satisfying appropriate conditions), so $\bPT_{AB}\cap\bPT_{CD}$ and $\bI_{AB}\cap\bI_{CD}$ are never empty (since the empty map $\emptyset$ belongs to $\bI_{AB}\sub \bPT_{AB}$ for all~$A,B$).  Similarly,~$\bT_{AB}\cap\bT_{CD}$ is non-empty if and only if $A=C=\emptyset$, or $A=C\not=\emptyset$ and $B\cap D\not=\emptyset$.
Note also that~$\bT_{A\emptyset}=\emptyset$ if $A\not=\emptyset$.  
For $A\in\Set$, we write $\bPT_A=\bPT_{AA}$, $\bT_A=\bT_{AA}$ and $\bI_A=\bI_{AA}$, noting that these are all semigroups under composition: they are the partial and full transformation semigroups, and the symmetric inverse semigroup on $A$, respectively.
At times, we will refer to sets such as $\bPT_{\al\be}$ or~$\bPT_\al$, where $\al$ and $\be$ are cardinals; in such instances, we simply regard $\al,\be$ as sets (so $\al$ is the set of all ordinals less than $\al$, etc.).

We define the class
\[
\PT = \bigset{(A,f,B)}{A,B\in\Set,\ f\in\bPT_{AB}},
\]
with partial product $\cdot$ defined by
\[
(A,f,B)\cdot(C,g,D) = \begin{cases}
(A,fg,D) &\text{if $B=C$}\\
\text{undefined} &\text{otherwise.}
\end{cases}
\]
As we have done elsewhere, we continue to write $xf$ for the image of $x$ under a mapping $f$; thus, $f$ is performed first in the composite $fg$.

We define mappings
\[
\lam:\PT\to\Set:(A,f,B)\mt A \AND \rho:\PT\to\Set:(A,f,B)\mt B.
\]
So $(A,f,B)\cdot(C,g,D)$ is defined if and only if $(A,f,B)\rho=(C,g,D)\lam$, and the (partial) associative law holds because of the associativity of binary relation composition:
\[
\big( (A,f,B)\cdot(B,g,C) \big) \cdot (C,h,D) = (A,f,B)\cdot \big( (B,g,C) \cdot (C,h,D) \big).
\]
For $A,B\in\Set$, 
\[
\PT_{AB} = \set{(A,f,B)}{f\in\bPT_{AB}}
\]
is a set.  So $(\PT,\cdot,\Set,\lam,\rho)$ is a partial semigroup.  

We write $\Setp=\Set\sm\{\emptyset\}$ for the class of all non-empty sets, and define the subclasses 
\[
\T = \bigset{(A,f,B)}{A,B\in\Setp,\ f\in\bT_{AB}} \AND
\I = \bigset{(A,f,B)}{A,B\in\Set,\ f\in\bI_{AB}}
\]
of $\PT$, noting that both $\T$ and $\I$ are closed under $\cdot$.  Writing $\lam,\rho$ also for the restrictions of these maps to~$\T$ and $\I$, we see that $(\T,\cdot,\Setp,\lam,\rho)$ and $(\I,\cdot,\Set,\lam,\rho)$ are both partial semigroups.  

As usual, we use the abbreviation $\PT\equiv(\PT,\cdot,\Set,\lam,\rho)$, and similarly for $\T$ and $\I$.
%abbreviations $\PT,\T,\I$ for the partial semigroups as defined above.  
Because of the identity maps $\id_A:A\to A:a\mt a$, each of $\PT,\T,\I$ is monoidal; equivalently, they are all (locally small) categories.

In what follows, we use the following standard notation.  If $f\in\bPT_{AB}$, we write $\dom(f)$, $\im(f)$, $\rank(f)$ and $\ker(f)$ for the \emph{domain}, \emph{image}, \emph{rank} and \emph{kernel} of $f$; the last two are defined by $\rank(f)=|\im(f)|$ and $\ker(f)=\set{(x,y)\in\dom(f)\times\dom(f)}{xf=yf}$.  So $\dom(f)\sub A$, $\im(f)\sub B$, $0\leq\rank(f)\leq\min(|A|,|B|)$, and $\ker(f)$ is an equivalence on $\dom(f)$.  We also write $f=\binom{F_i}{f_i}_{i\in I}$ to indicate that $\im(f)=\set{f_i}{i\in I}$ and $F_if=\set{xf}{x\in F_i}=\{f_i\}$ for all $i$; when we use this notation, we will always assume that $f_i\not=f_j$ if $i\not=j$ (that is, the map $I\to\im(f):i\mt f_i$ is injective).  Sometimes we just write $f=\binom{F_i}{f_i}$ with the indexing set $I$ being implied.  With the above notation, $\dom(f)=\bigcup_{i\in I}F_i$, and $\ker(f)=\bigcup_{i\in I}(F_i\times F_i)$.

\begin{prop}\label{prop:reg_PTTI}
The partial semigroups $\PT$, $\T$ and $\I$ are all regular categories.  
\end{prop}

%\newpage

\pf %It is clear that for any $A\in\Set$, the identity map $\id_A$ (which belongs to all three partial semigroups) satisfies condition (v).  
Let $(A,f,B)\in\PT$.  It suffices to show that:
\bit
\item[(i)] there exists $g\in\bI_{BA}$ such that $(A,f,B) = (A,f,B) \cdot (B,g,A) \cdot (A,f,B)$, and
\item[(ii)] if $A,B\not=\emptyset$, then there exists $h\in\bT_{BA}$ such that $(A,f,B) = (A,f,B) \cdot (B,h,A) \cdot (A,f,B)$.
\eit
We begin with (i).  If $f=\emptyset$ is the empty map, then we may clearly take $g=\emptyset$, so suppose $f\not=\emptyset$, and write $f=\binom{F_i}{f_i}_{i\in I}$.  For each $i\in I$, choose some $g_i\in F_i$, and put $g=\binom{f_i}{g_i}\in\bI_{BA}$; clearly, $(B,g,A)$ has the desired property.  If $A,B\not=\emptyset$, then any $h\in\bT_{BA}$ with $g\sub h$ satisfies (ii). \epf

%note that if $h\in\bPT_{BA}$ is arbitrary with the property that $g\sub h$, then clearly $(A,f,B) \cdot (B,g,A)=(A,f,B) \cdot (B,h,A)$; in particular, if $B\not=\emptyset$, then we may extend $h\in\bI_{BA}$ to $g\in\bT_{BA}$. \epf

Since $\PT$, $\T$ and $\I$ are all regular, every element of these categories is \emph{sandwich-regular}, so the general theory of \cite[Sections 2 and 3]{Sandwiches1} applies to any sandwich semigroup in these categories.  Also, since $\T$ and~$\I$ are regular subcategories of $\PT$, we may deduce several facts about $\T$ and $\I$ from corresponding facts about~$\PT$, using the results from \cite[Section~1.4]{Sandwiches1}: for example, information concerning Green's $\R$, $\L$ and~$\H$ relations \cite[Lemma 1.8]{Sandwiches1}; $\R$- and/or $\L$-stability of elements \cite[Lemma 1.9]{Sandwiches1}; the sets $P_1^a,P_2^a,P_3^a,P^a$ \cite[Lemma~1.10]{Sandwiches1}.
As such, we will first restrict our attention to the category $\PT$. %, and then explain how results on $\T$ and $\I$ may be deduced.

\sectiontitle{The category $\PT$}\label{sect:PT}

%XYXYXYX

We are now ready to conduct a thorough investigation of the partial transformation category $\PT$.  
We begin, in Section \ref{sect:Green_PT}, by characterising Green's relations and preorders in $\PT$, and classifying the $\R$- and/or~$\L$-stable elements.  
In Section \ref{sect:PTXYa}, we describe the sets $P_1^a,P_2^a,P^a,P_3^a$, using these to describe Green's relations and preorders on the sandwich semigroups $\PTXYa$ and characterise the regular elements; we also classify the regular $\D^a$-classes and maximal $\J^a$-classes in $\PTXYa$; some of the preliminary results of this section have been proved (sometimes in a very different form) in \cite{MGS2013,MS1975}.  
Section \ref{sect:structurePT} gives a structure theorem for the regular subsemigroup $\Reg(\PTXYa)$, and explores connections with certain non-sandwich semigroups,~$\PT(X,\B)$ and $\PT(Y,\si)$, of partial transformations of restricted range or kernel.  
Section \ref{sect:RegPTXYa} further explores $\Reg(\PTXYa)$, giving detailed structural and combinatorial information, including formulae for the size and rank of $\Reg(\PTXYa)$.  
In Section \ref{sect:EaPTXYa}, we describe and enumerate the idempotents of $\PTXYa$, and study the idempotent-generated subsemigroup $\E_a(\PTXYa)$, characterising the elements of this subsemigroup and calculating its rank and idempotent rank, which turn out to be equal.
Section \ref{sect:rank_PTXYa} calculates the rank of an arbitrary sandwich semigroup $\PTXYa$; the formulae given depend on whether the sandwich element is full and/or injective and/or surjective, and as one special case, we deduce a result from \cite{FS2014}.
Finally, Section \ref{sect:eggbox} gives egg-box diagrams of several sandwich semigroups $\PTXYa$; these may be used to visualise many of the results of the preceding sections.  
Throughout Section~\ref{sect:PT}, we will comment on various corollaries our general results have for the above-mentioned non-sandwich semigroups~$\PT(X,\B)$ and $\PT(Y,\si)$.

\subsectiontitle{Green's relations and stability in $\PT$}\label{sect:Green_PT}

Our first priority is to describe Green's relations and preorders.  We first collect some basic results about composition of partial transformations.  
Let $X,Y\in\Set$ with $X\sub Y$, and let $\si$ be an equivalence relation on $Y$.  We write $\si|_X=\si\cap(X\times X)$ for the \emph{restriction} of $\si$ to $X$.  We say that $X$ \emph{saturates} $\si$ if each $\si$-class contains at least one element of $X$.  We say $\si$ \emph{separates} $X$ if each $\si$-class contains at most one element of~$X$.  We say $X$ is a \emph{cross-section} of $\si$ if $X$ saturates and is separated by $\si$.
The next result is easily proved.

\begin{lemma}\label{lem:PT_prelim}
Let $A,B,C\in\Set$, and let $f\in\bPT_{AB}$ and $g\in\bPT_{BC}$.  Then
\bit
\itemit{i} $\dom(fg)\sub\dom(f)$, with equality if and only if $\im(f)\sub\dom(g)$,
\itemit{ii} $\im(fg)\sub\im(g)$, with equality if and only if $\im(f)$ saturates $\ker(g)$,
\itemit{iii} $\ker(fg)\supseteq\ker(f)|_{\dom(fg)}$, with equality if and only if $\ker(g)$ separates $\im(f)\cap\dom(g)$,
\itemit{iv} $\rank(fg)\leq\min(\rank(f),\rank(g))$. \epfres
%\itemit{v} $\dom(fg)=\dom(f) \iff \im(f)\sub\dom(g)$,
%\itemit{vi} $\im(fg)=\im(g) \iff \im(f)$ saturates $\ker(g)$,
%\itemit{vii} $\ker(fg)=\ker(f)|_{\dom(fg)} \iff \ker(g)$ separates $\im(f)\cap\dom(g)$.
\eit
\end{lemma}

In what follows, we will often use Lemma \ref{lem:PT_prelim} without explicit reference.

\newpage

\begin{prop}\label{prop:GreenPT}
Let $(A,f,B),(C,g,D)\in\PT$.  Then
\bit
\itemit{i} $(A,f,B)\leqR(C,g,D) \iff A=C$, $\dom(f)\sub\dom(g)$ and $\ker(f)\supseteq\ker(g)|_{\dom(f)}$,
\itemit{ii} $(A,f,B)\leqL(C,g,D) \iff B=D$ and $\im(f)\sub\im(g)$,
\itemit{iii} $(A,f,B)\leqJ(C,g,D) \iff \rank(f)\leq\rank(g)$,
\itemit{iv} $(A,f,B)\R(C,g,D) \iff A=C$, $\dom(f)=\dom(g)$ and $\ker(f)=\ker(g)$,
\itemit{v} $(A,f,B)\L(C,g,D) \iff B=D$ and $\im(f)=\im(g)$,
\itemit{vi} $(A,f,B)\J(C,g,D) \iff (A,f,B)\D(C,g,D) \iff \rank(f)=\rank(g)$.
\eit
\end{prop}

\pf (i).  Suppose $(A,f,B)\leqR(C,g,D)$, so that $(A,f,B)=(C,g,D)\cdot(E,h,F)$ for some $(E,h,F)\in\PT$.  So $D=E$ and $(A,f,B)=(C,gh,F)$, which gives $A=C$ and $f=gh$ (and $B=F$).  From $f=gh$, we deduce that $\dom(f)=\dom(gh)\sub\dom(g)$, and that $\ker(f)=\ker(gh)\supseteq\ker(g)|_{\dom(gh)} = \ker(g)|_{\dom(f)}$.

%.  Also, if $(x,y)\in\ker(g)|_{\dom(f)}$, then $x,y\in\dom(f)$ and $xf=(xg)h=(yg)h=yf$, so that $(x,y)\in\ker(f)$.

Conversely, suppose $A=C$, $\dom(f)\sub\dom(g)$ and $\ker(g)|_{\dom(f)}\sub\ker(f)$.  We define $h\in\bPT_{DB}$ as follows.  First, define $\dom(h)=\dom(f)g=\set{xg}{x\in\dom(f)}$.  For $x\in\dom(f)$, we define $(xg)h=xf$.  Note that $h$ is well defined, since $\ker(g)|_{\dom(f)}\sub\ker(f)$.  One then easily checks that $(A,f,B)=(C,g,D)\cdot(D,h,B)$.

\pfitem{ii}  The forwards implication is easy.  Conversely, suppose $B=D$ and $\im(f)\sub\im(g)$, and write $f=\binom{F_i}{f_i}$.  For each $i$, choose some $g_i\in f_ig^{-1}$, and put $h=\binom{F_i}{g_i}\in\bPT_{AC}$.  Then $(A,f,B)=(A,h,C)\cdot(C,g,D)$.

\pfitem{iii}  The forwards implication is again easy to check.  Conversely, suppose $\rank(f)\leq\rank(g)$, and write $f=\binom{F_i}{f_i}_{i\in I}$ and $g=\binom{G_j}{g_j}_{j\in J}$.  Without loss of generality, we may assume that $I\sub J$.  Choose some $e_i\in G_i$ for each $i$.  Put $h_1=\binom{F_i}{e_i}\in\bPT_{AC}$ and $h_2=\binom{g_i}{f_i}\in\bPT_{DB}$.  Then $(A,f,B)=(A,h_1,C)\cdot(C,g,D)\cdot(D,h_2,B)$.

\pfitem{iv) and (v}  These follow immediately from (i) and (ii), respectively, noting that $\ker(g)|_{\dom(g)}=\ker(g)$.

\pfitem{vi}  By (iii), and since ${\D}\sub{\J}$, it suffices to show that $\rank(f)=\rank(g)\implies(A,f,B)\D(C,g,D)$.  So suppose $\rank(f)=\rank(g)$, and write $f=\binom{F_i}{f_i}_{i\in I}$ and $g=\binom{G_i}{g_i}_{i\in I}$.  Then parts (iv) and (v) give $(A,f,B)\L(C,h,B)\R(C,g,D)$, where $h=\binom{G_i}{f_i}\in\bPT_{CB}$. \epf

The next result follows immediately from parts (iii) and (vi) of Proposition \ref{prop:GreenPT}.

\begin{cor}\label{cor:Jclasses_PT}
Let $A,B\in\Set$.  The ${\J}={\D}$-classes of $\PT_{AB}$ are the sets
\[
D_\mu = D_\mu(\PT_{AB}) = \set{(A,f,B)}{f\in\bPT_{AB},\ \rank(f)=\mu} \qquad\text{for each cardinal \ $0\leq\mu\leq\min(|A|,|B|)$.}
\]
These $\J$-classes form a chain: $D_\mu\leqJ D_\nu \iff \mu\leq\nu$. \epfres
\end{cor}

We will also need to know the sizes of certain $\K$-classes in $\PT$.  For this, we require the following notation.
For cardinals $\ka,\mu$ with $\mu\leq\ka$, we write $\binom{\ka}{\mu}$ for the number of subsets of size $\mu$, and $S(\ka,\mu)$ for the number of equivalence relations with $\mu$ equivalence classes, in a set of size $\ka$; we also write $\ka!$ for the size of the symmetric group on a set of size $\ka$.  When $\ka$ is finite, these are just the ordinary binomial coefficients, Stirling numbers of the second kind, and factorials, respectively.  When $\ka$ is infinite, $\binom\ka\mu=\ka^\mu$, $S(\ka,1)=1$, $S(\ka,\mu)=2^\ka$ for $\mu\geq2$, and $\ka!=2^\ka$.
If $\ka<\mu$, then we define $\binom\ka\mu=S(\ka,\mu)=0$.  %We also write $\mu!$ for the size of the symmetric group on a set of size $\mu$, for any 
%(possibly infinite) 
%cardinal $\mu$.  

The next result also follows quickly from Proposition \ref{prop:GreenPT}; note that for (i), if $A\in\Set$, then a pair $(D,K)$ where $D\sub A$ and $K$ is an equivalence on $D$ may be identified with an equivalence $K'$ on $A\cup\{\infty\}$, where~$\infty$ is a symbol that does not belong to $A$
% (if $|D/K|=\mu$, then $|(A\cup\{\infty\})/K'|=\mu+1$).
 (the $K'$-class containing $\infty$ is $\{\infty\}\cup(A\sm D)$).

\begin{cor}\label{cor:Green_sizes_PT}
Let $A,B\in\Set$, write $\al=|A|$ and $\be=|B|$, and let $0\leq\mu\leq\min(\al,\be)$.  Then
\begin{itemize}\begin{multicols}{2}
%\itemit{i} $|D_\mu/{\R}| = \sum_{\ka=\mu}^\al\binom\al\ka S(\ka,\mu)$,
\itemit{i} $|D_\mu/{\R}| = S(\al+1,\mu+1)$,
\itemit{ii} $|D_\mu/{\L}| = \binom\be\mu$,
%\itemit{iii} $|D_\mu/{\H}| = \binom\be\mu\sum_{\ka=\mu}^\al\binom\al\ka S(\ka,\mu)$,
\itemit{iii} $|D_\mu/{\H}| = \binom\be\mu S(\al+1,\mu+1)$,
\itemit{iv} each $\H$-class in $D_\mu$ has size $\mu!$,
%\itemit{v} $|D_\mu| = \mu! \binom\be\mu\sum_{\ka=\mu}^\al\binom\al\ka S(\ka,\mu)$.
\itemit{v} $|D_\mu| = \mu! \binom\be\mu S(\al+1,\mu+1)$.
\item[] ~ \epfres
\end{multicols}\eitmc
\end{cor}

\begin{rem}
By considering the size of $\bPT_{AB}$, Corollary \ref{cor:Green_sizes_PT}(v) gives rise to the identity
\[
(\be+1)^\al=\sum_{\mu=0}^{\min(\al,\be)}\mu! \tbinom\be\mu S(\al+1,\mu+1).
\]
Note that the sum is over all cardinals $\mu$ satisfying $0\leq\mu\leq\min(\al,\be)$.
\end{rem}

We have already seen that every element of $\PT$ is regular.  Next we wish to characterise the $\R$- and/or~$\L$-stable elements.  First we prove a preliminary result.  For a set $X$, we write
\[
\bPTfr_X=\set{f\in\bPT_X}{\rank(f)<\aleph_0}
\]
for the set of all finite-rank elements of $\bPT_X$.  Note that $\bPTfr_X$ is a subsemigroup of $\bPT_X$, by Lemma \ref{lem:PT_prelim}(iv).  Recall that a semigroup $T$ is \emph{periodic} if for each $x\in T$, some power of $x$ is an idempotent.

\begin{lemma}\label{lem:PTfr}
If $X$ is any set, then $\bPTfr_X$ is a periodic semigroup.
\end{lemma}

\pf Let $f\in\bPTfr_X$.  
%If $f^k=\emptyset$ for some $k\geq1$, then we are done, so suppose otherwise.  
The sequence $\rank(f),\rank(f^2),\rank(f^3),\ldots$ is non-increasing, and since its first term is finite, it must eventually become constant.  Suppose $m\geq1$ and $r\geq0$ are such that $\rank(f^k)=r$ for all $k\geq m$, and write $f^m=\binom{F_i}{f_i}_{i\in I}$.  Then every element of the set $\Om=\set{f^{k}}{k\geq m}$ is of the form $\binom{F_i}{f_{i\pi}}$, for some permutation $\pi$ of $I$.  Since $|I|=\rank(f^m)<\aleph_0$, it follows that $\Om$ is a finite semigroup and, hence, contains an idempotent. \epf

We say $f\in\bPT_{AB}$ is \emph{full} if $\dom(f)=A$: that is, if $f\in\bT_{AB}$.

\begin{lemma}\label{lem:stable_PT}
If $(A,f,B)\in\PT$, then
\bit
\itemit{i} $(A,f,B)$ is $\R$-stable $\iff$ $[\rank(f)<\aleph_0$ or $f$ is full and injective$]$,
\itemit{ii} $(A,f,B)$ is $\L$-stable $\iff$ $[\rank(f)<\aleph_0$ or $f$ is surjective$]$,
\itemit{iii} $(A,f,B)$ is stable $\iff$ $[\rank(f)<\aleph_0$ or $f$ is bijective$]$.
\eit
\end{lemma}

\pf First, suppose $\rank(f)<\aleph_0$.  Note that 
\[
(A,f,B)\cdot\PT_{BA}\sub\set{(A,g,A)}{g\in\bPTfr_A} \AND \PT_{BA}\cdot(A,f,B)\sub\set{(B,g,B)}{g\in\bPTfr_B}.
\]
%$(A,f,B)\cdot\PT_{BA}\sub\set{(A,g,A)}{g\in\bPTfr_A}$ and $\PT_{BA}\cdot(A,f,B)\sub\set{(B,g,B)}{g\in\bPTfr_B}$.  
Since $\bPTfr_A$ and $\bPTfr_B$ are both periodic, by Lemma \ref{lem:PTfr}, it follows from \cite[Lemma 1.3]{Sandwiches1} that $(A,f,B)$ is stable.

Next note that if $f$ is full and injective,
%$\dom(f)=A$ and $f$ is injective, 
then $\dom(gf)=\dom(g)$ and $\ker(gf)=\ker(g)$ for any $C\in\Set$ and $g\in\bPT_{CA}$, so that $(C,g,D)$ is $\R$-related to $(C,g,D)\cdot(A,f,B)$ whenever the latter product is defined: that is, $(A,f,B)$ is $\R$-stable.  Similarly, if $f$ is surjective, then $\im(fg)=\im(g)$ for any $D\in\Set$ and $g\in\bPT_{BD}$, and $\L$-stability of $(A,f,B)$ quickly follows.

Since (iii) clearly follows from (i) and (ii), it now just remains to prove the forward implications in each of~(i) and (ii).  In both cases, we do this by proving the contrapositive.  For the remainder of the proof, suppose $\rank(f)\geq\aleph_0$, and write $f=\binom{F_i}{f_i}_{i\in I}$.  Choose some $g_i\in F_i$ for each $i$.  Since $\rank(f)\geq\aleph_0$, it also follows that $|A|,|B|\geq\aleph_0$.

\pfitem{i}  Suppose that either (a) $f$ is not full, or (b) $f$ is not injective.
%
%one of the following holds:
%\begin{itemize}\begin{multicols}{2}
%\item[(a)] $\dom(f)\not=A$, or 
%\item[(b)] $f$ is not surjective.
%\end{multicols}\eitmc
In case (a), fix some ${a\in A\sm\dom(f)}$.  In case (b), let $i\in I$ be such that $|F_i|\geq2$, and fix some $a\in F_i\sm\{g_i\}$.  In either case, put $g=\binom{a\ g_i}{a\ g_i}_{i\in I}\in\bPT_{AA}$.  Then $\rank(gf)=\rank(g)$ in both cases, but $\dom(gf)\not=\dom(g)$ in case (a), while $\ker(gf)\not=\ker(g)$ in case~(b).  In both cases, it follows that $(A,g,A)\cdot(A,f,B)$ and $(A,g,A)$ are $\J$-related but not $\R$-related.  So $(A,f,B)$ is not $\R$-stable.

%Next, suppose $f$ is not injective.  Let $i\in I$ be such that $|F_i|\geq2$, choose some $a\in F_i\sm\{g_i\}$, and put $g=\binom{a\ g_i}{a\ g_i}_{i\in I}\in\bPT_{AA}$.  Then $\rank(gf)=\rank(g)$, but $\ker(gf)\not=\ker(g)$.  Again, it quickly follows that $(A,f,B)$ is not $\R$-stable.
%
%For (i), suppose first that $\dom(f)\not=A$.  Choose some $a\in A\sm\dom(f)$, and put $g=\binom{a\ g_i}{a\ g_i}_{i\in I}\in\bPT_{AA}$.  Then $\rank(gf)=\rank(g)$, but $\dom(gf)\not=\dom(g)$.  It follows that $(A,g,A)\cdot(A,f,B)$ and $(A,g,A)$ are $\J$-related but not $\R$-related.  So $(A,f,B)$ is not $\R$-stable.
%
%Next, suppose $f$ is not injective.  Let $i\in I$ be such that $|F_i|\geq2$, choose some $a\in F_i\sm\{g_i\}$, and put $g=\binom{a\ g_i}{a\ g_i}_{i\in I}\in\bPT_{AA}$.  Then $\rank(gf)=\rank(g)$, but $\ker(gf)\not=\ker(g)$.  Again, it quickly follows that $(A,f,B)$ is not $\R$-stable.

\pfitem{ii}  Suppose $f$ is not surjective.  Fix some $b\in B\sm\im(f)$, and put $g=\binom{b\ f_i}{b\ f_i}_{i\in I}\in\bPT_{BB}$.  Then $\rank(fg)=\rank(g)$ but $\im(fg)\not=\im(g)$, so that $(A,f,B)\cdot(B,g,B)$ and $(B,g,B)$ are $\J$-related but not $\L$-related.  So $(A,f,B)$ is not $\L$-stable.  \epf

\subsectiontitle{Green's relations, regularity and stability in $\PTXYa$}\label{sect:PTXYa}

Now that we have gathered the required preliminary material on the category $\PT$, we are now ready to study sandwich semigroups in $\PT$.  For the rest of Section \ref{sect:PT}, we fix two sets $X,Y\in\Set$.
In order to simplify notation, %in this section, 
we will identify $\bPT_{ZW}$ with $\PT_{ZW}$, where $Z,W\in\{X,Y\}$, via $f\equiv(Z,f,W)$.
For the rest of Section \ref{sect:PT}, we also fix a partial transformation $a\in\PT_{YX}$, with the aim of studying the sandwich semigroup~$\PTXYa$.  
We write
\begin{equation}\label{eq:aPT}
a=\tbinom{A_i}{a_i}_{i\in I} \COMMA
\B=\im(a)=\set{a_i}{i\in I} \COMMA
\A=\dom(a)=\bigcup_{i\in I}A_i \COMMA
\si=\ker(a) \COMMA
\al=\rank(a),
\end{equation}
so that $\al=|I|=|\B|=|\A/\si|$ and $\A/\si=\set{A_i}{i\in I}$.  
For each $i\in I$, we fix some $b_i\in A_i$, and write
\begin{equation}\label{eq:bPT}
b=\tbinom{a_i}{b_i}\in\PT_{XY},
\end{equation}
so that $a=aba$ and $b=bab$.  (More generally, we could pick any $b=\binom{B_i}{b_i}$, where $\set{B_i}{i\in I}$ is any partition of any subset of $X$ with $a_i\in B_i$ for each $i$.)  A number of other parameters will play a role in later calculations, but it will be convenient to define them here, so all notation is fixed at the beginning:
\begin{equation}\label{eq:cPT}
\be=|X\sm\im(a)| \COMMA \xi=\min(|X|,|Y|) \COMMA \lam_i=|A_i| \text{ for } i\in I \COMMA \Lam_J=\prod_{j\in J}\lam_j \text{ for } J\sub I.
\end{equation}
%an element $f\in\bPT_{ZW}$ with the corresponding element $(Z,f,W)\in\PT_{ZW}$, where $Z,W\in\{X,Y\}$.  Thus, the main object of our study in this section, is a sandwich semigroup $\PT_{XY}^a$.  
%By Theorem \ref{thm:green_Sij}, Green's relations on $\PTXYa$ are governed by the sets $P_1^a,P_2^a,P_3^a,P^a$.  
Our first result describes the sets 
\[
P_1^a = \set{f\in\PTXY}{fa\R f} \COMma
P_2^a = \set{f\in\PTXY}{af\L f} \COMma
P_3^a = \set{f\in\PTXY}{afa\J f} \COMma
P^a=P_1^a\cap P_2^a.
\]
%$P_1^a,P_2^a,P_3^a,P^a$; it 
It is an easy consequence of Lemma \ref{lem:PT_prelim} and Proposition~\ref{prop:GreenPT}.

\begin{prop}\label{prop:P_sets_PT}
We have
\bit
\itemit{i} $P_1^a = \set{f\in\PTXY}{\dom(fa)=\dom(f),\ \ker(fa)=\ker(f)}$
\item[] $\phantom{P_1^a} = \set{f\in\PTXY}{\im(f)\sub\dom(a),\ \ker(a)\text{ separates }\im(f)}$,
\itemit{ii} $P_2^a = \set{f\in\PTXY}{\im(af)=\im(f)} = \set{f\in\PTXY}{\im(a)\text{ saturates }\ker(f)}$,
%\itemit{ii} $P_2^a = \set{f\in\PTXY}{\im(af)=\im(f)}$
%\item[] $\phantom{P_2^a} = \set{f\in\PTXY}{\im(a)\text{ saturates }\ker(f)}$,
\itemit{iii} $P^a = \set{f\in\PTXY}{\dom(fa)=\dom(f),\ \ker(fa)=\ker(f),\ \im(af)=\im(f)}$
\item[] $\phantom{P^a} = \set{f\in\PTXY}{\im(f)\sub\dom(a),\ \ker(a)\text{ separates }\im(f),\ \im(a)\text{ saturates }\ker(f)}$,
\itemit{iv} $P_3^a = \set{f\in\PTXY}{\rank(afa)=\rank(f)}$. \epfres
\eitres
\end{prop}

\begin{rem}
Some simplifications arise in special cases.  For example, if $a$ is full, then $\im(f)\sub\dom(a)$ is automatically satisfied by any $f\in\PTXY$, and so $P_1^a=\set{f\in\PTXY}{\ker(a)\text{ separates }\im(f)}$.  Similarly, if $a$ is injective, then $P_1^a=\set{f\in\PTXY}{\im(f)\sub\dom(a)}$.
If $a$ is full \emph{and} injective, then $P_1^a=\PTXY$; cf.~\cite[Lemma 1.2]{Sandwiches1} and Lemma \ref{lem:al=xi}.
%
%or injective, then $P_1^a=\set{f\in\PTXY}{\ker(a)\text{ separates }\im(f)}$ or $\set{f\in\PTXY}{\im(f)\sub\dom(a)}$, respectively.
\end{rem}

Together with \cite[Theorem 1.1]{Sandwiches1}, Proposition \ref{prop:P_sets_PT} yields a description of Green's relations on $\PTXYa$: % the sandwich semigroup~$\PTXYa$: 

\begin{thm}\label{thm:green_PTXYa}
If $f\in \PTXY$, then   
\begin{itemize}\begin{multicols}{2}
\itemit{i} $R_f^a = \begin{cases}
R_f\cap P_1^a &\text{if $f\in P_1^a$}\\
\{f\} &\text{if $f\not\in P_1^a$,}
\end{cases}$
%\phantom{\begin{cases}a\\b\\c\\d\end{cases}}
\itemit{ii} $L_f^a = \begin{cases}
L_f\cap P_2^a &\hspace{0.7mm}\text{if $f\in P_2^a$}\\
\{f\} &\hspace{0.7mm}\text{if $f\not\in P_2^a$,}
\end{cases}
%\phantom{\begin{cases}a\\b\\c\\d\end{cases}}
$
\itemit{iii} $H_f^a = \begin{cases}
H_f &\hspace{7.4mm}\text{if $f\in P^a$}\\
\{f\} &\hspace{7.4mm}\text{if $f\not\in P^a$,}
\end{cases}$
\itemit{iv} $D_f^a = \begin{cases}
D_f\cap P^a &\text{if $f\in P^a$}\\
L_f^a &\text{if $f\in P_2^a\sm P_1^a$}\\
R_f^a &\text{if $f\in P_1^a\sm P_2^a$}\\
\{f\} &\text{if $f\not\in P_1^a\cup P_2^a$,}
\end{cases}$
\itemit{v} $J_f^a = \begin{cases}
J_f\cap P_3^a &\hspace{2.2mm}\text{if $f\in P_3^a$}\\
D_f^a &\hspace{2.2mm}\text{if $f\not\in P_3^a$.}
\end{cases}$
\end{multicols}\end{itemize}
Further, if $f\not\in P^a$, then $H_f^a=\{f\}$ is a non-group $\gHa$-class of $\PTXYa$.  \epfres
\end{thm}

\begin{rem}
Figures \ref{fig:PT_2}--\ref{fig:PT_6} give the so-called egg-box diagrams for various sandwich semigroups $\PTXYa$; as explained in Section \ref{sect:eggbox}, these display the structure of $\PTXYa$ as determined by Green's relations.
Green's relations on~$\PTXYa$ were also characterised in \cite[Theorems 2.6, 2.7 and 2.8]{MGS2013}, although the presentations of the results from \cite{MGS2013} are quite different from Theorem \ref{thm:green_PTXYa}, with the exception of the $\R$ and $\L$ relations.
\end{rem}

Recall that ${\D}={\J}$ in the category $\PT$ (see Proposition \ref{prop:GreenPT}(vi)).  Our next main result (Proposition~\ref{prop:JaDaPT}) shows that the $\gJa$ and $\gDa$ relations on the sandwich semigroup $\PTXYa$ need not coincide.  But first we prove a technical result that will be useful on a number of occasions.

\begin{lemma}\label{lem:JaDaPT}
Suppose $\mu$ is a cardinal with $\aleph_0\leq\mu\leq\al=\rank(a)$.
\bit
\itemit{i} If $a$ is not $\R$-stable, then there exists some $f\in P_3^a\sm P_1^a$ with $\rank(f)=\mu$.
\itemit{ii} If $a$ is not $\L$-stable, then there exists some $f\in P_3^a\sm P_2^a$ with $\rank(f)=\mu$.
\itemit{iii} If $a$ is not stable, then there exists some $f\in P_3^a\sm P^a$ with $\rank(f)=\mu$.
\eit
\end{lemma}

\pf By assumption, $\al=\rank(a)$ is infinite, so we may fix some proper subset $J\subsetneq I$ with $|J|=\mu$, and some $k\in I\sm J$.  It suffices to prove (i) and (ii), since $P_3^a\sm P_q^a\sub P_3^a\sm P^a$ for $q=1,2$.

\pfitem{i}  Suppose $a$ is not $\R$-stable.  By Lemma \ref{lem:stable_PT}(i), either (a) $a$ is not full, or (b) $a$ is not injective.  
In case~(a), choose some $y\in Y\sm \dom(a)$.  In case (b), without loss of generality, we may assume that $|A_l|\geq2$ for some $l\in J$, and we then choose some $y\in A_l\sm\{b_l\}$.  In both cases, put $f=\tran{a_k&a_j}{y&b_j}_{j\in J}\in\PTXY$.  Then $\mu=\rank(f)=\rank(afa)$, so that $f\in P_3^a$.  But $\dom(fa)\not=\dom(f)$ in case~(a), and $\ker(fa)\not=\ker(f)$ in case~(b), so that $f\not\in P_1^a$.
%, whence $g\not\in P^a$.  
%Thus, $g\in J_f\cap(P_3^a\sm P_1^a)\sub J_f\cap(P_3^a\sm P^a)$.

\pfitem{ii}  Suppose $a$ is not $\L$-stable.  Then $a$ is not surjective, by Lemma \ref{lem:stable_PT}(ii).  Choose some $x\in X\sm\im(a)$, and put $f=\tran{x&a_j}{b_k&b_j}_{j\in J}\in\PTXY$.  Again, we quickly obtain $\mu=\rank(f)=\rank(afa)$, so that $f\in P_3^a$, but $\im(af)\not=\im(f)$, so that $f\not\in P_2^a$. \epf

%\pfitem{iii}  This follows imediately from (i) and (ii), since $P_3^a\sm P_q^a\sub P_3^a\sm P^a$ for $q=1,2$.  \epf

%
\begin{prop}\label{prop:JaDaPT}
In $\PTXYa$, we have ${\gJa}={\gDa} \iff a$ is stable.
\end{prop}

\pf %Suppose first that $a$ is stable, and let $f\in\PTXY$.  We must show that $J_f^a=D_f^a$.  This follows immediately from Theorem \ref{thm:green_PTXYa}(v) if $f\not\in P_3^a$, so suppose $f\in P_3^a$.  Proposition \ref{prop:Reg(Sija)}(iii) gives $P_3^a=P^a$ (since $a$ is stable), while  Proposition \ref{prop:GreenPT}(vi) gives $J_f=D_f$.  Together with Theorem \ref{thm:green_PTXYa}(iv) and (v), it then follows that $J_f^a = J_f\cap P_3^a = D_f\cap P^a = D_f^a$.
If $a$ is stable, then \cite[Corollary 1.5]{Sandwiches1} and Proposition \ref{prop:GreenPT}(vi) gives ${\gJa}={\gDa}$.
Conversely, suppose $a$ is not stable.  We will show that $J_b^a\not=D_b^a$.  Since $a$ is not stable, Lemma \ref{lem:stable_PT}(iii) says that $\al=\rank(a)\geq\aleph_0$.  By Lemma~\ref{lem:JaDaPT}(iii), there exists 
%Choose any cardinal $\aleph_0\leq\mu\leq\al$, and let 
$f\in P_3^a\sm P^a$ with $\rank(f)=\al=\rank(b)$.  
%Also, let $J\sub I$ with $|J|=\mu$, and put $g=\binom{a_j}{b_j}_{j\in J}\in\PTXY$, so that $g\in P^a$ by Proposition \ref{prop:P_sets_PT}(iii).  
By Proposition \ref{prop:GreenPT}(vi) and Theorem \ref{thm:green_PTXYa}(iv) and (v), since $b\in P^a\sub P_3^a$, we have $D_b^a=D_b\cap P^a=J_b\cap P^a$ and $J_b^a=J_b\cap P_3^a$.  So $f\in J_b\cap(P_3^a\sm P^a)=J_b^a\sm D_b^a$. \epf

The next result shows how \cite[Proposition 1.4]{Sandwiches1} may be strengthened in the category $\PT$.  Among other things, it characterises the regular elements of $\PTXYa$, and shows that the inclusion $P^a\sub P_3^a$ (which holds in any partial semigroup) can sometimes be strict.

\begin{prop}\label{prop:Reg(PYXYa)}
We have $\Reg(\PTXYa) = P^a \sub P_3^a$.  Further,
\begin{itemize}\begin{multicols}{3}
\itemit{i} $a$ is $\R$-stable $\iff$ $P_3^a\sub P_1^a$,
\itemit{ii} $a$ is $\L$-stable $\iff$ $P_3^a\sub P_2^a$,
\itemit{iii} $a$ is stable $\iff$ $P_3^a=P^a$.
\end{multicols}\eitmc
\end{prop}

\pf By \cite[Proposition 1.7(v)]{Sandwiches1} and Proposition \ref{prop:reg_PTTI}, we have $\Reg(\PTXYa) = P^a$.  
By \cite[Proposition 1.4]{Sandwiches1}, 
%and since (iii) follows from (i) and (ii), 
it remains 
%suffices to show that $P^a\sub\Reg(\PTXYa)$, and 
to prove the reverse implication in (i)--(iii).  
%
%For the first, suppose $f=\binom{F_j}{f_j}\in P^a=P_1^a\cap P_2^a$, so that $\dom(fa)=\dom(f)$, $\ker(fa)=\ker(f)$ and $\im(af)=\im(f)$.  So we may write $fa=\binom{F_j}{g_j}$ and $af=\binom{G_j}{f_j}$.  For each $j$, choose some $h_j\in G_j$, and put $g=\binom{g_j}{h_j}\in\PTXY$.  Then $f=fagaf=f\star_ag\star_af$, so that $f\in\Reg(\PTXYa)$, as required.
%
For (i), we prove the contrapositive.  Suppose $a$ is not $\R$-stable.  By Lemma \ref{lem:stable_PT}(i), $\rank(a)\geq\aleph_0$, and Lemma~\ref{lem:JaDaPT}(i) says that $P_3^a\sm P_1^a$ is non-empty.  Part~(ii) is treated in similar fashion.  Part (iii) follows from~(i) and (ii). \epf

\begin{rem}
The characterisation $\Reg(\PTXYa)=P^a$ may be deduced from \cite[Theorem 5.3]{MS1975}.
\end{rem}

We now prove a number of results concerning the ordering $\leq_{\J^a}$ on the  $\gJa$-classes of $\PTXYa$.  For simplicity, we abbreviate $\leq_{\J^a}$ to $\leq$.

\newpage

\begin{prop}\label{prop:Ja_order_PT}
Let $f,g\in\PT_{XY}$.  Then $J_f^a\leq J_g^a$ in $\PTXYa$ if and only if one of the following holds:
\begin{itemize}\begin{multicols}{2}
\itemit{i} $f=g$,
\itemit{ii} $\rank(f)\leq\rank(aga)$,
\itemit{iii} $\im(f)\sub\im(ag)$,
\itemit{iv} $\dom(f)\sub\dom(ga)$ and $\ker(f)\supseteq\ker(ga)|_{\dom(f)}$.
\end{multicols}\eitmc
%The maximal $\gDa$-classes are those of the form $D_f^a=\{f\}$ where $\rank(f)>r$.
%{\red What are the maximal $\gJa$-classes?}
\end{prop}

\pf Note that $J_f^a\leq J_g^a$ if and only if one of the following holds:
\begin{itemize}\begin{multicols}{2}
\item[(a)] $f=g$,
\item[(b)] $f=uagav$ for some $u,v\in\PT_{XY}$,
\item[(c)] $f=uag$ for some $u\in \PT_{XY}$,
\item[(d)] $f=gav$ for some $v\in\PT_{XY}$.
\end{multicols}\eitmc
Then (b) $\iff$ $f\leqJ aga$ (in $\PT$) $\iff$ (ii), by Proposition \ref{prop:GreenPT}(iii).  Other parts of Proposition \ref{prop:GreenPT} may be used to show that (c) $\iff$ (iii), and (d) $\iff$ (iv). \epf

The next result shows how Proposition \ref{prop:Ja_order_PT} may be simplified in the case that one (or both) of $f,g$ belongs to one (or more) of the sets $P_1^a,P_2^a,P_3^a$.

\begin{prop}\label{prop:Ja_order_PT_P}
Let $f,g\in\PT_{XY}$.  
\bit
\itemit{i} If $f\in P_1^a$, then $J_f^a\leq J_g^a$ $\iff$ $[\rank(f)\leq \rank(aga) \text{ or } [\dom(f)\sub\dom(ga)$ and $\ker(f)\supseteq\ker(ga)|_{\dom(f)}]]$.
\itemit{ii} If $f\in P_2^a$, then $J_f^a\leq J_g^a \iff [\rank(f)\leq \rank(aga) \text{ or } \im(f)\sub\im(ag)]$.
\itemit{iii} If $f\in P_3^a$, then $J_f^a\leq J_g^a \iff \rank(f)\leq \rank(aga)$.
\itemit{iv} If $g\in P_1^a$, then $J_f^a\leq J_g^a \iff [\rank(f)\leq \rank(ag) \text{ or } [\dom(f)\sub\dom(g)$ and $\ker(f)\supseteq\ker(g)|_{\dom(f)}]]$.
\itemit{v} If $g\in P_2^a$, then $J_f^a\leq J_g^a \iff [\rank(f)\leq \rank(ga) \text{ or } \im(f)\sub\im(g)]$.
\itemit{vi} If $g\in P_3^a$, then $J_f^a\leq J_g^a \iff \rank(f)\leq \rank(g)$.
\eit
\end{prop}

\pf Again, $J_f^a\leq J_g^a$ if and only if one of (a)--(d) holds, as enumerated in the proof of Proposition \ref{prop:Ja_order_PT}.
%, $J_f^a\leq J_g^a$ if and only if 
%\begin{itemize}\begin{multicols}{2}
%\item[(a)] $f=g$,
%\item[(b)] $f=uagav$ for some $u,v\in\PT_{XY}$,
%\item[(c)] $f=uag$ for some $u\in \PT_{XY}$,
%\item[(d)] $f=gav$ for some $v\in\PT_{XY}$.
%\end{multicols}\eitmc

\pfitem{i} Suppose $f\in P_1^a$, so that $f=fah$ for some $h\in\PT_{XY}$.  Now,
\begin{align*}
\text{(b) } & \text{$\implies$ $f\leqJ aga$ $\implies$ $\rank(f)\leq\rank(aga)$, }\\
\text{(d) }& \text{$\implies$ $f\leqR ga$ $\implies$ $[\dom(f)\sub\dom(ga)$ and $\ker(f)\supseteq\ker(ga)|_{\dom(f)}]$,}\\
\text{(a) }& \text{$\implies$ $f=fah=gah$ $\implies$ (d), }\\
\text{(c) }& \text{$\implies f=fah=uagah \implies$ (b).}
\end{align*}
%Also,
%\[
%\text{(a) $\implies$ $f=fah=gah$ $\implies$ (d) \AND (c) $\implies f=fah=uagah \implies$ (b).}
%\]
This establishes the forwards implication.  The converse follows from Proposition \ref{prop:Ja_order_PT}.

\pfitem{ii} Suppose $f\in P_2^a$, so that $f=haf$ for some $h\in\PT_{XY}$.  This time,
\[
\text{(b) } \implies \rank(f)\leq\rank(aga) \COMMA
\text{(c) } \implies \im(f)\sub\im(ag) \COMMA
\text{(a) $\implies$ (c)} \COMMA
\text{(d) $\implies$ (b)} .
\]
Again, the converse follows from Proposition \ref{prop:Ja_order_PT}.

\pfitem{iii} If $f\in P_3^a$, then $f=h_1afah_2$ for some $h_1,h_2\in\PTXY$, and it is then easy to see that any of (a)--(d) implies $\rank(f)\leq\rank(aga)$: for example, (c) $\implies f=h_1afah_2=h_1au(aga)h_2$.  The converse is again clear.

%Suppose $f\in P^a=P_1^a\cap P_2^a$.  By the previous proofs, (b)$\implies \rank(f)\leq\rank(aga)$, while (c)$\implies$(b), (d)$\implies$(b), and (a)$\implies$(d)$\implies$(b).

%By (i) and (ii), any of (b), (c) or (d) implies $\rank(f)\leq\rank(aga)$.  Since $f\in P^a\sub P_3^a$, we have $f=h_1afah_2$ for some $h_1,h_2\in\PT_{XY}$, so (a) $\implies f=h_1afah_2 = h_1agah_2 \implies$ (b).  Again, Proposition \ref{prop:Ja_order_PT} completes the proof.

\pfitem{iv}  Suppose $g\in P_1^a$.  Clearly, (a) and (d) each imply $[\dom(f)\sub\dom(g)$ and $\ker(f)\supseteq\ker(g)|_{\dom(f)}]$, while~(b) and (c) each imply $\rank(f)\leq\rank(ag)$.  To prove the converse, suppose first that $\dom(f)\sub\dom(g)$ and $\ker(f)\supseteq\ker(g)|_{\dom(f)}$.  Since $g\in P_1^a$, Proposition \ref{prop:P_sets_PT}(i) gives $\dom(g)=\dom(ga)$ and $\ker(g)=\ker(ga)$, so it follows that $\dom(f)\sub\dom(ga)$ and $\ker(f)\supseteq\ker(ga)|_{\dom(f)}$, and we obtain $J_f^a\leq J_g^a$ from Proposition~\ref{prop:Ja_order_PT}.  Finally, suppose $\rank(f)\leq\rank(ag)$.  Since $g\in P_1^a$, we have $g=gah$ for some $h\in\PTXY$, so that $\rank(f)\leq\rank(ag)=\rank(agah)\leq\rank(aga)$,
%Now, 
%\[
%g\in P_1^a \implies g\R ga \implies ag\R aga \implies ag \J aga \implies \rank(ag)=\rank(aga),
%\]
%%where we have used Corollary \ref{cor:Green_relations_T}(iii), and the facts that $\R$ is a left-congruence and ${\R}\sub{\J}$.  In particular, 
%so $\rank(f)\leq\rank(ag)=\rank(aga)$, 
and $J_f^a\leq J_g^a$ again follows from Proposition \ref{prop:Ja_order_PT}.

\pfitem{v}  Suppose $g\in P_2^a$.  Clearly, (a) and (c) each imply~$\im(f)\sub\im(g)$, while (b) and (d) each imply~${\rank(f)\leq\rank(ga)}$.  To prove the converse, suppose first that $\im(f)\sub\im(g)$.  Since $g\in P_2^a$, Proposition~\ref{prop:P_sets_PT}(ii) gives $\im(g)=\im(ag)$, so $\im(f)\sub\im(ag)$, and we obtain $J_f^a\leq J_g^a$ from Proposition~\ref{prop:Ja_order_PT}.  Finally, suppose $\rank(f)\leq\rank(ga)$.  Since $g\in P_2^a$, we have $g=hag$ for some $h\in\PTXY$, so that ${\rank(f)\leq\rank(ga)=\rank(haga)\leq\rank(aga)}$, and $J_f^a\leq J_g^a$ again follows from Proposition \ref{prop:Ja_order_PT}.

%\pfitem{vi}  Finally, suppose $g\in P^a$.  Clearly, any of (a)--(d) implies $\rank(f)\leq\rank(g)$.  Conversely, suppose $\rank(f)\leq\rank(g)$.  Since $g\in P^a\sub P_3^a$, $\rank(g)=\rank(aga)$, so $\rank(f)\leq\rank(aga)$, giving $J_f^a\leq J_g^a$.~\epf

\pfitem{vi}  Finally, suppose $g\in P_3^a$.  Clearly, any of (a)--(d) implies $\rank(f)\leq\rank(g)$.  Conversely, suppose $\rank(f)\leq\rank(g)$.  Since $g\in P_3^a$, $\rank(g)=\rank(aga)$, so $\rank(f)\leq\rank(aga)$, giving $J_f^a\leq J_g^a$.~\epf

\begin{rem}
Since $P^a\sub P_3^a$, parts (iii) and (vi) of Proposition \ref{prop:Ja_order_PT_P} apply to elements of $P^a$.
\end{rem}

 % \implies g\R ga \implies g\J ga \implies \rank(g)=\rank(ga)$, so we deduce $J_f^a\leq J_g^a$ from part (v). \epf

Recall that the elements of a $\D$-class of a semigroup are either all regular or all non-regular; see for example \cite[Proposition 2.3.1]{Howie}.  The next result characterises the regular $\gDa$-classes of $\PTXYa$.
Recall from Proposition~\ref{prop:Reg(PYXYa)} that $P^a=\Reg(\PTXYa)$.

\begin{prop}\label{prop:regular_Da_classes}
The regular $\gDa$-classes of $\PTXYa$ are precisely the sets
\[
D_\mu^a = \set{f\in P^a}{\rank(f)=\mu} \qquad\text{for each cardinal $0\leq\mu\leq\al=\rank(a)$.}
\]
Further, if $f\in P^a$, then $D_f^a=J_f^a$ if and only if $\rank(f)<\aleph_0$ or $a$ is stable.
\end{prop}

\pf Let $f\in P^a$, and note that $\rank(f)=\rank(afa)\leq\rank(a)=\al$, since $f\in P^a\sub P_3^a$.  By Theorem~\ref{thm:green_PTXYa}(iv) and Proposition \ref{prop:GreenPT}(vi), $D_f^a=D_f\cap P^a = J_f\cap P^a = D_\mu^a$, where $\mu=\rank(f)$, showing that all regular $\gDa$-classes are of the specified form.  Since there exist regular elements of $\PTXYa$ of any rank from $0$ to $\al$ (for example, $\binom{a_j}{b_j}_{j\in J}$ is regular, for any subset $J\sub I$), this completes the proof of the first assertion.

Next, note that Theorem \ref{thm:green_PTXYa}(v) gives $J_f^a=J_f\cap P_3^a$, since $f\in P^a\sub P_3^a$.
If $a$ is stable, then Proposition~\ref{prop:JaDaPT} gives $D_f^a=J_f^a$.
Now suppose $\rank(f)<\aleph_0$.  To prove that $D_f^a=J_f^a$, it suffices to show that $J_f^a\sub D_f^a$.  With this in mind, let $g\in J_f^a=J_f\cap P_3^a$.  So $\rank(aga)=\rank(g)=\rank(f)<\aleph_0$.  Since $D_f^a=J_f\cap P^a$, we must show that $g\in P^a=P_1^a\cap P_2^a$.  Now,
\[
\rank(g)=\rank(aga)\leq\left\{ \begin{array}{c} \rank(ag) \\ \rank(ga) \end{array}\right\} \leq\rank(g),
\]
so that $\rank(ag)=\rank(ga)=\rank(g)$.  Now, $\im(ag)\sub\im(g)$, and since these are finite sets of the same size (equal to $\rank(ag)=\rank(g)$), it follows that $\im(g)=\im(ag)$, whence $g\L ag$, giving $g\in P_2^a$.
%
%To show that $g\in P_1^a$, write $g=\binom{G_j}{g_j}$.  
Since $\rank(ga)=\rank(g)<\aleph_0$, it follows that $\im(g)\sub\dom(a)$, and that $\ker(a)$ separates $\im(g)$, so that $g\in P_1^a$, by Proposition \ref{prop:P_sets_PT}(i), as required.  This proves the backwards implication of the second assertion.

For the forwards implication, we prove the contrapositive.  Suppose that $\rank(f)\geq\aleph_0$ and that $a$ is not stable.  By Lemma \ref{lem:JaDaPT}(iii), there exists $g\in P_3^a\sm P^a$ with $\rank(g)=\rank(f)$.  Together with Theorem~\ref{thm:green_PTXYa}(iv) and (v), this gives $g\in J_f\cap(P_3^a\sm P^a) = (J_f\cap P_3^a)\sm(J_f\cap P^a) = J_f^a\sm D_f^a$. \epf

We will have more to say about these regular $\D^a$-classes in Section \ref{sect:RegPTXYa}.  
We close this section with a description of the maximal $\J^a$-classes of $\PTXYa$.  Recall that we write $\xi=\min(|X|,|Y|)$.

\begin{prop}\label{prop:maximal_J_PT}
\begin{itemize}
\itemit{i} If $\al<\xi$, then the maximal $\J^a$-classes of $\PTXYa$ are precisely the singleton sets~$\{f\}$, for $f\in\PTXY$ with $\rank(f)>\al$.
\itemit{ii} If $\al=\xi$, then the set $J_b^a=\set{f\in P_3^a}{\rank(f)=\al}$ is a maximum $\J^a$-class of $\PTXYa$.
%\itemit{iii} If $\al=\xi<\aleph_0$, then $J_b^a=D_\al^a=\Hh_b^a$.
\end{itemize}
\end{prop}

\pf (i).  Suppose $\al<\xi$, and let $f\in\PTXY$ with $\rank(f)>\al$.  Now,
\[
\rank(af),\rank(fa),\rank(afa)\leq\rank(a)=\al<\rank(f),
\]
so it follows that $f\in\PTXY\sm(P_1^a\cup P_2^a\cup P_3^a)$.  Theorem \ref{thm:green_PTXYa}(iv) and (v) then give $J_f^a=D_f^a=\{f\}$.  To show that $J_f^a=\{f\}$ is maximal, suppose $J_f^a\leq J_g^a$ for some $g\in\PTXY$.  Then one of (a)--(d) from the proof of Lemma \ref{prop:Ja_order_PT} hold; any of (b)--(d) would imply $\rank(f)\leq\rank(a)=\al$, so it follows that $f=g$, as required.

To show that these are the only maximal $\J^a$-classes, suppose $g\in\PTXY$ is such that $\rank(g)\leq\al$, and write $g=\binom{G_j}{g_j}_{j\in J}$ where $J\sub I$.  Put $h_1=\binom{G_j}{b_j}$.  Since $|J|=\rank(g)\leq\al<\xi$, there is an element $h_2\in\I_{XY}$ extending $\binom{a_j}{g_j}$ and with $\rank(h_2)>\al$.  It is easy to check that $g=h_1\star_ah_2$, so that $J_g^a\leq J_{h_2}^a$.  But $\rank(g)\leq\al<\rank(h_2)$, so $(g,h_2)\not\in{\J}$ and $J_g^a\not=J_{h_2}^a$, showing that $J_g^a$ is not maximal.  %This completes the proof of (i).

\pfitem{ii}  Suppose $\al=\xi$.  Since $b\in P^a\sub P_3^a$, Theorem \ref{thm:green_PTXYa}(v) says that $J_b^a=J_b\cap P_3^a$ is the stated set.  If $g\in\PTXY$ is arbitrary, then $\rank(g)\leq\xi=\al=\rank(b)$; since $b\in P^a$, Proposition \ref{prop:Ja_order_PT_P}(vi) gives $J_g^a\leq J_b^a$.~\epf

\begin{rem}
The reader may locate the maximal $\J^a$-classes in various sandwich semigroups $\PTXYa$ in Figures \ref{fig:PT_2}--\ref{fig:PT_6}.  Note that $\al<\xi$ in Figures \ref{fig:PT_2}--\ref{fig:PT_4} (in which case there are many (singleton) maximal $\J^a$-classes), while $\al=\xi$ in Figures \ref{fig:PT_5} and \ref{fig:PT_6} (in which case there is a unique maximal $\J^a$-class).
\end{rem}

%$A=\dom(a)\not=Y$ or (b) $a$ is not injective.  Without loss of generality, we may assume that $\{1,2,3,\ldots\}\sub I$.  In case (a), choose some $y\in Y\sm A$.  In case (b), relabelling if necessary, choose some $y\in A_2\sm\{b_2\}$.  In either case, define $f=\left(\begin{smallmatrix} a_1&a_2&a_3&\cdots\\y&b_2&b_3&\cdots \end{smallmatrix}\right)\in\PTXY$.  Then $\rank(afa)=\rank(f)$ in both cases, so that $f\in P_3^a$.  However, $\dom(fa)\not=\dom(f)$ in case (a), while $\ker(fa)\not=\ker(f)$ in case (b), so that~$f\not\in P_1^a$.
%
%For (ii), suppose $a$ is not $\L$-stable, so that $|I|\geq\aleph_0$ and $a$ is not surjective.  Again, we assume that $\{1,2,3,\ldots\}\sub I$.  Choose some $y\in X\sm B$, and define $f=\left(\begin{smallmatrix} x&a_2&a_3&\cdots\\b_1&b_2&b_3&\cdots \end{smallmatrix}\right)\in\PTXY$.  Then $\rank(afa)=\rank(f)$ but $\im(af)\not=\im(f)$, so that $f\in P_3^a\sm P_2^a$. \epf

%In case (b), relabelling if necessary, choose some 
%
%suppose $P_3^a\sub P_1^a$.  To show that $a$ is $\R$-stable, suppose $f=\binom{F_j}{f_j}_{j\in J}\in\PT$ is such that $fa\J f$.  It suffices to show that $f\in P_3^a$.  

%Suppose first that (a) holds.  Choose some $y\in Y\sm A$, and put $g=\tran{a_j&a_i}{y&b_i}_{i\in I'}\in\PTXY$.  Then $\rank(aga)=\rank(g)$, so that $g\in P_3^a$.  However, $\dom(ga)\not=\dom(g)$, so that $g\not\in P_1^a$, whence $g\not\in P^a$.  Since also $\rank(g)=\rank(f)$, it follows that $g\in J_f$.
%
%Next, suppose (b) holds.  Choose some 

\subsectiontitle{A structure theorem for $\Reg(\PTXYa)$ and connections to (non-sandwich) semigroups of partial transformations}\label{sect:structurePT}

Again, we fix $a=\binom{A_i}{a_i}\in\PTYX$, and we keep the notation of \eqref{eq:aPT}--\eqref{eq:cPT}.  
From \cite[diagrams (2.1) and~(2.7)]{Sandwiches1}, we have the following two commutative diagrams, with all maps being semigroup epimorphisms:
\begin{equation}\label{eq:CD_PT}
\includegraphics{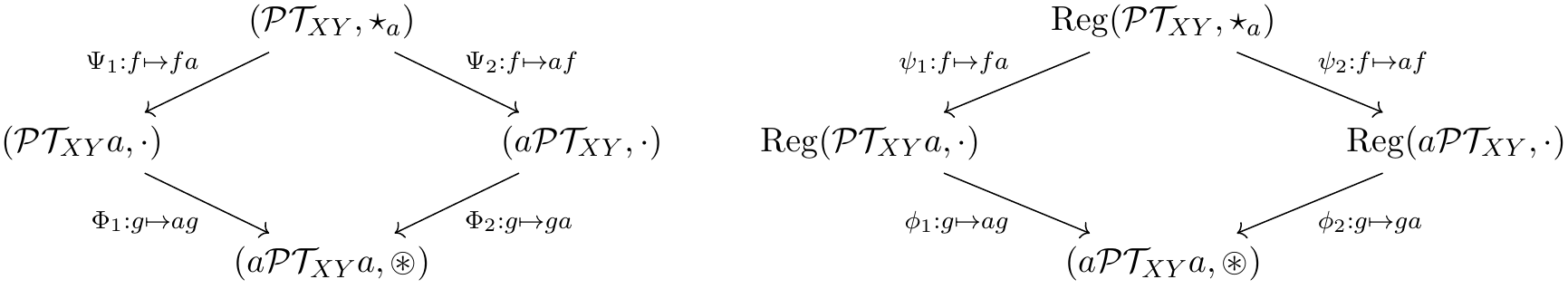}
%\begin{tikzcd} [column sep=small]
%~ & (\PTXY,\star_a) \arrow[swap]{dl}{\Psi_1:f\mt fa} \arrow{dr}{\Psi_2:f\mt af} & \\
%(\PTXY a,\cdot) \arrow[swap]{dr}{\Phi_1:g\mt ag} & & (a\PTXY ,\cdot) \arrow{dl}{\Phi_2:g\mt ga}\\
%& (a\PTXY a,\starb) & 
%\end{tikzcd}
%\qquad %\qquad
%\begin{tikzcd} [column sep=small]
%~ & \Reg(\PTXY ,\star_a) \arrow[swap]{dl}{\psi_1:f\mt fa} \arrow{dr}{\psi_2:f\mt af} & \\
%\Reg(\PTXY a,\cdot) \arrow[swap]{dr}{\phi_1:g\mt ag} & & \Reg(a\PTXY ,\cdot) \arrow{dl}{\phi_2:g\mt ga}\\
%& (a\PTXY a,\starb) & 
%\end{tikzcd}
\end{equation}
Recall that $\starb$ denotes the restriction of the $\star_b$ operation to $a\PTXY a\sub\PTYX$.  At the top of these diagrams, we have the sandwich semigroup $\PTXYa$, left, and its regular subsemigroup $P^a=\Reg(\PTXYa)$, right.  We now describe the other semigroups appearing in these diagrams.  Recall that $\B=\im(a)$, $\A=\dom(a)$ and $\si=\ker(a)$; in particular, $\si$ is an equivalence on $\A\sub Y$.  We write
\[
\PT(X,\B) = \set{f\in\PT_X}{\im(f)\sub \B} 
\]
for the set of partial transformations on $X$ with \emph{image restricted by $\B$}, and
\[
\PT(Y,\si)=\set{f\in\PT_Y}{
%\dom(f)\sub A,\ 
\text{every $\ker(f)$-class is a union of $\si$-classes}}
\]
for the set of partial transformations on $Y$ with \emph{kernel restricted by $\si$}; note that $\dom(f)\sub\A$ for all $f\in\PT(Y,\si)$.  It is easy to see that
\[
\PT(X,\B)=\PTXY a=\PT_Xba \AND \PT(Y,\si)=a\PTXY=ab\PT_Y.
\]
In particular, $\PT(X,\B)$ is a subsemigroup (indeed, a principal left ideal) of $\PT_X$, and $\PT(Y,\si)$ a subsemigroup (indeed, a principal right ideal) of $\PT_Y$; these semigroups, and their regular subsemigroups, make up the middle rows of the diagrams in \eqref{eq:CD_PT}.  The semigroups $\PT(X,\B)$ have been studied in \cite{FS2014}, the main results being a classification of the regular elements, a description of Green's relations, and the calculation of $\rank(\PT(X,\B))$ in the case of finite $X$.  
%
%Since the semigroups $\PT(X,\B)$ arise as special cases of the $\PTXYa$ construction, the results of the current article yield corresponding results concerning $\PT(X,\B)$, including (but not limited to) all the results from \cite{FS2014}.  
%
To the authors' knowledge, no systematic study of the semigroups $\PT(Y,\si)$ has been carried out.

Finally, we recall from \cite[Remark 2.6]{Sandwiches1} that the regular monoid $(a\PTXY a,\starb)$ appearing at the bottom of both diagrams in~\eqref{eq:CD_PT} is isomorphic to the local monoid $ba\PTXY a=ba\PT_Xba$ of $\PT_X$ with respect to the idempotent $ba=\binom{a_i}{a_i}\in\PT_X$, which is the partial identity map on $\B=\im(a)=\set{a_i}{i\in I}$.  
We will denote by $\eta:(a\PTXY a,\starb)\to(ba\PT_Xba,\cdot):f\mt bf$ this isomorphism.
The local monoid $ba\PT_Xba$ is clearly isomorphic to $\PT_\B$, and will therefore be identified with $\PT_\B$.
Thus, taking into account the above discussion, we see that the diagrams from \eqref{eq:CD_PT} become:
\begin{equation}\label{eq:CD_PT_2}
\includegraphics{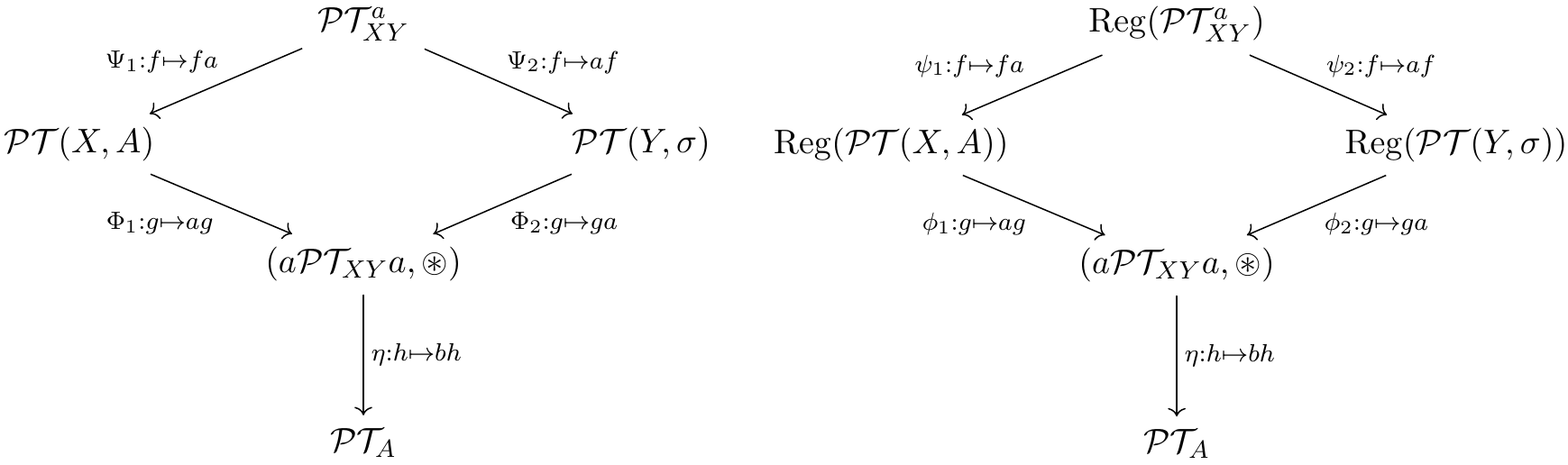}
%\begin{tikzcd} %[column sep=small]
%~ & \PTXYa \arrow[swap]{dl}{\Psi_1:f\mt fa} \arrow{dr}{\Psi_2:f\mt af} & \\
%\PT(X,\B) \arrow[swap]{dr}{\Phi_1:g\mt ag} & & \PT(Y,\si) \arrow{dl}{\Phi_2:g\mt ga}\\
%& (a\PTXY a,\starb)\arrow{dd}{\eta:h\mt bh} & \\
%\\
%& \PT_\B & 
%\end{tikzcd}
%\quad %\qquad
%\begin{tikzcd} [column sep=small]
%~ & \Reg(\PTXYa) \arrow[swap]{dl}{\psi_1:f\mt fa} \arrow{dr}{\psi_2:f\mt af} & \\
%\Reg(\PT(X,\B)) \arrow[swap]{dr}{\phi_1:g\mt ag} & & \Reg(\PT(Y,\si)) \arrow{dl}{\phi_2:g\mt ga}\\
%& (a\PTXY a,\starb)\arrow{dd}{\eta:h\mt bh}  & \\
%\\
%& \PT_\B & 
%\end{tikzcd}
\end{equation}
We will say more about the diagrams in \eqref{eq:CD_PT_2} shortly, but first we make some comments about certain special cases:~
\bit
\item If $a$ is full and injective, then $a$ is right-invertible (indeed, $ab=\id_Y$ in this case).  As discussed near the beginning of \cite[Section 2]{Sandwiches1}, this implies that the right cancellation law \cite[equation (2.2)]{Sandwiches1} holds, and so 
%$\Psi_1:\PTXYa \to (\PTXY a,\cdot) = \PT(X,B)$
${\Psi_1:\PTXYa \to \PT(X,\B):f\mt fa}$
is an isomorphism (as is $\psi_1$).  Thus, Figure \ref{fig:PT_5} shows an egg-box diagram for $\PT(X,\B)$, where $X=\{1,2,3,4\}$ and $\B=\{1,2,3\}$.
\item If $a$ is surjective, then $a$ is left-invertible (indeed, $ba=\id_X$ in this case), so \cite[equation (2.3)]{Sandwiches1} holds, and so
%$\Psi_2:\PTXYa \to (a\PTXY ,\cdot) = \PT(Y,\si)$
$\Psi_2:\PTXYa \to \PT(Y,\si):f\mt af$
is an isomorphism (as is $\psi_2$).  Thus, Figure \ref{fig:PT_6} shows egg-box diagrams for $\PT(Y,\si_1)$, $\PT(Y,\si_2)$ and $\PT(Y,\si_3)$, where $Y=\{1,2,3,4,5\}$, and $\si_1=(\ 1\mid2\mid3,4\ )$, $\si_2=(\ 1,2\mid3,4\mid5\ )$ and $\si_3=(\ 1\mid2\mid3\ )$, respectively, using an obvious notation for equivalences on subsets of $Y$.
\eit
Consequently, the (non-sandwich) semigroups $\PT(X,\B)$ and $\PT(Y,\si)$ arise as special cases of the $\PTXYa$ construction.  All the results we prove concerning arbitrary sandwich semigroups $\PTXYa$ therefore have corollaries for the semigroups $\PT(X,\B)$ and $\PT(Y,\si)$; in the former case, these include (but are not limited to) the above-mentioned results from \cite{FS2014}.  As a sub-case of the second special case above, we have the following:
\bit
\item If $a$ is injective and surjective (but not necessarily full), then $\PTXYa$ is isomorphic to $\PT(Y,\si)$, where $\si=\set{(y,y)}{y\in \A}$ is the diagonal relation on $\A$.  In this case, $\PT(Y,\si)=\set{f\in\PT_Y}{\dom(f)\sub\A}$.  Again, as far as the authors are aware, these semigroups have not been studied.  Figure \ref{fig:PT_6} (right) pictures an egg-box diagram of this semigroup in the case that $Y=\{1,2,3,4,5\}$ and $\A=\{1,2,3\}$.
\eit
The last special case combines the first two:
\bit
\item If $a$ is a (full) bijection, then \cite[equations (2.2) and (2.3)]{Sandwiches1} both hold, so that all maps in \eqref{eq:CD_PT_2} are isomorphisms; in particular, $\PTXYa \cong\PT(X,\B)=\PT(X,X)=\PT_X$ in this case.  Figure \ref{fig:PT_1} shows egg-box diagrams for $\PT_X$ when $|X|=2$ or $3$.
%
%so that 
%%$\PTXYa \cong (a\PTXY a,\starb) \cong \PT_B$.
%$\Psi_1\Phi_1\xi=\Psi_2\Phi_2\xi:\PTXYa \to \PT_B:f\mt bafa$ is an isomorphism.  Note that $B=X$, and so $\PTXYa$ is isomorphic to $\PT_X$, in this case.
\eit
Since most of the problems we consider for the sandwich semigroups $\PTXYa$ have been solved for the semigroups $\PT_X$ themselves, we generally exclude the case in which $a$ is a (full) bijection.  However, we allow for the possibility that $a$ is full and injective, \emph{or} that $a$ is surjective.

Since $\PT$ is a regular category, every element is sandwich-regular.  So, together with the above discussion, we immediately obtain the following structural result concerning $P^a=\Reg(\PTXYa)$ from \cite[Theorem 2.10]{Sandwiches1}.

\begin{thm}\label{thm:psi_PT}
Consider the map
\[
%\psi=(\psi_1,\psi_2):P^a=\Reg(\PTXYa) \to \Reg(\PT(X,B))\times\Reg(\PT(Y,\si)):f\mt(fa,af).
\psi:\Reg(\PTXYa) \to \Reg(\PT(X,\B))\times\Reg(\PT(Y,\si)):f\mt(fa,af).
\]
Then $\psi$ is injective, and
\[
\im(\psi) = \set{(g,h)\in \Reg(\PT(X,\B))\times\Reg(\PT(Y,\si))}{ag=ha}.
\]
%\begin{align*}
%\im(\psi) &= \set{(g,h)\in \Reg(\PT(X,B))\times\Reg(\PT(Y,\si))}{ag=ha} \\
%&= \set{(g,h)\in \Reg(\PT(X,B))\times\Reg(\PT(Y,\si))}{g\phi_1=h\phi_2}.
%\end{align*}
%
In particular, $\Reg(\PTXYa)$ is a pullback product of $\Reg(\PT(X,\B))$ and $\Reg(\PT(Y,\si))$ with respect to~$\PT_\B$.~\epfres
\end{thm}

Note that Theorem \ref{thm:psi_PT} makes reference to the regular semigroups $\Reg(\PT(X,\B))$ and $\Reg(\PT(Y,\si))$.  As a consequence of \cite[Proposition 2.5]{Sandwiches1} and Proposition \ref{prop:P_sets_PT}, we may describe these semigroups.

\begin{prop}\label{prop:RegPTXBYsi}
The sets $\Reg(\PT(X,\B))$ and $\Reg(\PT(Y,\si))$ are (regular) subsemigroups of $\PT(X,\B)$ and $\PT(Y,\si)$, and
\begin{align*}
\Reg(\PT(X,\B)) = \Reg(\PTXYa) a &= \set{f\in \PT(X,\B)}{\text{$\B$ saturates $\ker(f)$}},\\
\epfreseq \Reg(\PT(Y,\si)) = a\Reg(\PTXYa) &= \set{f\in \PT(Y,\si)}{\text{$\im(f)\sub \A$, $\si$ separates $\im(f)$}}.
\end{align*}
\end{prop}

\begin{rem}
The statement in Proposition \ref{prop:RegPTXBYsi} concerning $\Reg(\PT(X,\B))$ was proved in \cite[Theorem~1.2]{FS2014}, with the condition ``$\B$ saturates $\ker(f)$'' replaced with the equivalent condition ``$Xf=\B f$''; here, $Zf$ denotes the set $\set{zf}{z\in Z\cap\dom(f)}$ for any subset $Z\sub X$.
\end{rem}

\begin{rem}\label{rem:PTYsiPTA}
In the case that $\si$ is the diagonal relation on $\A\sub Y$, Proposition \ref{prop:RegPTXBYsi} gives
\[
\Reg(\PT(Y,\si)) = \set{f\in \PT(Y,\si)}{\im(f)\sub \A} = \set{f\in\PT_Y}{\dom(f),\im(f)\sub \A} \equiv\PT_{\A}.
\]
Compare the right-most diagrams in Figures \ref{fig:PT_1}, \ref{fig:PT_6} and \ref{fig:PT_9}.
\end{rem}

\begin{rem}
Since the semigroups $\PT(X,\B)$ and $\PT(Y,\si)$ arise as special cases of the $\PTXYa$ construction (as discussed above), one could also obtain Proposition \ref{prop:RegPTXBYsi} directly from Proposition \ref{prop:P_sets_PT}(iii), with suitable choices of $a$.
\end{rem}

\subsectiontitle{The regular subsemigroup $P^a=\Reg(\PTXYa)$}\label{sect:RegPTXYa}

Again, we fix $a=\binom{A_i}{a_i}\in\PTYX$, and keep the notation of \eqref{eq:aPT}--\eqref{eq:cPT}.  In this section, we further analyse the structure of the regular subsemigroup $P^a=\Reg(\PTXYa)$ of the sandwich semigroup $\PTXYa$.  In particular, we prove combinatorial results concerning $\K$- and $\gKh$-classes (Theorem \ref{thm:inflation_PT}), and also give the size and rank of $P^a$ (Corollary \ref{cor:size_P_PT}, Proposition \ref{prop:size_P_PT} and Theorem \ref{thm:rank_P_PT}).

%The proof of the next lemma is simple, but included for convenience.  {\red repetition now.}
%
%
%\begin{lemma}\label{lem:PTDJ}
%Suppose $T$ is a semigroup for which $P=\Reg(T)$ is a subsemigroup of $T$.  If $\K$ is any of Green's relations on $T$ other than $\J$, then ${\K^P}={\K}|_P$.
%\end{lemma}
%
%\pf This is well known if $\K$ is any of $\R,\L,\H$; see Lemma \ref{lem:Tgreen}.  Clearly ${\D^P}\sub{\D}|_P$.  Conversely, suppose $(x,y)\in{\D}|_P$.  So $x,y\in P$, and there exists $u\in T$ such that $x\R u\L y$ (in $T$).  Since any element in the $\R$-class of a regular element is regular, it follows that $u\in P$.  But then $(x,u)\in{\R}|_P={\R^P}$ and similarly $(u,y)\in{\L^P}$, so that $x\R^Pu\L^Py$, giving $(x,y)\in{\D^P}$. \epf

Because of \cite[Lemma 2.8]{Sandwiches1}, Green's $\R$, $\L$, $\H$ and $\D$ relations on $P^a=\Reg(\PTXYa)$ are simply the restrictions of the corresponding relations on $\PTXYa$.  So we will continue to denote these relations on $P^a$ by $\R^a$, $\L^a$,~$\H^a$ and $\D^a$.

\begin{lemma}
We have ${\J^{P^a}}={\D^a}$.
\end{lemma}

\pf Since ${\D}\sub{\J}$ in any semigroup, it is enough to show that ${\J^{P^a}}\sub{\D^a}$, so suppose $(f,g)\in{\J^{P^a}}$.  In particular, it follows that $(f,g)\in{\J}$, and so $\rank(f)=\rank(g)$, by Proposition \ref{prop:GreenPT}(vi): write $\mu$ for this common cardinality.  Then $f,g$ both belong to the (regular) $\D^a$-class~$D_\mu^a$, as defined in Proposition \ref{prop:regular_Da_classes}, so that $(f,g)\in{\D^a}$, as required. \epf

Thus, we will never need to refer to the ${\J^{P^a}}$ relation.  If $\K$ is any of $\R,\L,\H,\D$, then any $\K^a$-class in $\PTXYa$ contains only regular elements or only non-regular elements; it follows that every $\K^a$-class of $P^a$ is a $\K^a$-class of $\PTXYa$.  Thus, with no possibility of confusion, for $f\in P^a$, we will continue to write $K_f^a$ for the $\K^a$-class of $f$ in $P^a$.  We therefore obtain the following from Theorem \ref{thm:green_PTXYa} and Propositions \ref{prop:GreenPT} and~\ref{prop:regular_Da_classes}:

\begin{prop}\label{prop:green_Pa_PT}
If $f\in P^a=\Reg(\PTXYa)$, then
\bit
\itemit{i} $R_f^a = R_f\cap P^a = \set{g\in P^a}{\dom(f)=\dom(g),\ \ker(f)=\ker(g)}$,
\itemit{ii} $L_f^a = L_f\cap P^a = \set{g\in P^a}{\im(f)=\im(g)}$,
\itemit{iii} $H_f^a = H_f\cap P^a = \set{g\in P^a}{\dom(f)=\dom(g),\ \ker(f)=\ker(g),\ \im(f)=\im(g)}$,
\itemit{iv} $D_f^a = D_f\cap P^a = \set{g\in P^a}{\rank(f)=\rank(g)}$.
\eit
The ${\J^{P^a}}={\D^a}$-classes of $P^a$ are the sets
\[
D_\mu^a = \set{g\in P^a}{\rank(g)=\mu} \qquad\text{for each cardinal $0\leq\mu\leq\al=\rank(a)$,}
\]
and these form a chain under the $\J^{P^a}$-ordering: $D_\mu^a\leq D_\nu^a \iff \mu\leq\nu$. \epfres
\end{prop}

\begin{rem}
The chain of ${\J^{P^a}}={\D^a}$-classes in $P^a=\Reg(\PTXYa)$ has minimum element $D_0^a=\{\emptyset\}$, and maximum element $D_\al^a = \set{g\in P^a}{\rank(g)=\al}$.  Note that $b\in D_\al^a$.  Figure \ref{fig:PT_9} shows egg-box diagrams for several regular sandwich semigroups $\Reg(\PTXYa)$.  Parts (i)--(iv) of Proposition \ref{prop:green_Pa_PT} were also proved in \cite[Theorem 5.7]{MS1975}.
\end{rem}

Now that we know the ${\J^{P^a}}={\D^a}$-classes of $P^a$ form a chain, we wish to describe the internal structure of these classes.  To do so,
%In order to describe the internal structure of the $\D^a$-classes of $P^a=\Reg(\PTXYa)$, 
we must use the $\gKh$ relations defined in \cite[Section 2]{Sandwiches1}.  Recall that we have an epimorphism 
\[
\phi=\psi_1\phi_1=\psi_2\phi_2:P^a\to(a\PTXY a,\starb):f\mt afa.
\]
Here it will be convenient to compose $\phi$ with the isomorphism $\eta:(a\PTXY a,\starb)\to\PT_\B:f\mt bf$ discussed in Section \ref{sect:structurePT}, and instead work with the equivalent epimorphism
\[
\varphi = \psi_1\phi_1\eta=\psi_2\phi_2\eta:P^a\to\PT_\B:f\mt bafa.
\]
For this reason, we will slightly abuse previous notation and, for $f\in P^a$, we will write $\fb=f\varphi=bafa\in\PT_\B$.
%
%Combining this with the isomorphism $\Xi:(a\PTXY a,\starb)\to\PT_B:f\mt (fb)|_B$ yields an equivalent epimorphism
%\[
%\varphi=\phi\circ\Xi : P^a \to \PT_B : f \mt (afab)|_B.
%\]
%Since we intend to work with $\varphi$ rather than $\phi$, we will write $\fb=f\varphi=(afab)|_B$ for each $f\in P^a$, rather than $\fb=f\phi=afa$, as we did in Sections \ref{sect:RegSija} and \ref{sect:LMC}.  

If $\K$ is one of $\R$, $\L$, $\H$ or $\D$, and if $f,g\in P^a$, we write $f\gKh^a g$ if $\fb\K\gb$ in $\PT_\B$.  So, for example,~$f\gRh^ag \iff \dom(\fb)=\dom(\gb)$ and $\ker(\fb)=\ker(\gb)$.  For $f\in P^a$, we will write $\Kh^a_f$ for the $\gKh^a$-class of $f$ in~$P^a$.  Recall from \cite[Lemma 2.11]{Sandwiches1} that ${\gDh^a}={\D^a}$.

%Consider $f=\trans{F_j}{f_j}_{j\in J}\in P^a$.  By Proposition \ref{prop:P_sets_PT}, $\im(f)\sub\dom(a)$ and $\ker(a)$ separates $\im(f)$, so there is an injective map $J\to I:j\mt i_j$ such that $f_j\in A_{i_j}$ for each $j\in J$.

Consider a regular element $f\in P^a$.  By Proposition \ref{prop:P_sets_PT}, $\im(f)\sub\dom(a)$ and $\ker(a)$ separates $\im(f)$, so we may write $f=\binom{F_j}{f_j}_{j\in J}$, where $J\sub I$ and $f_j\in A_j$ for each $j\in J$.
%$f=\trans{F_j}{f_j}_{j\in J}\in P^a$. By Proposition \ref{prop:P_sets_PT}, $\im(f)\sub\dom(a)$ and $\ker(a)$ separates $\im(f)$, so there is an injective map $J\to I:j\mt i_j$ such that $f_j\in A_{i_j}$ for each $j\in J$.  
Since $\B=\im(a)$ saturates $\ker(f)$, it follows that $F_j\cap \B\not=\emptyset$ for all $j\in J$, so we may write $F_j\cap \B=\set{a_i}{i\in I_j}$ for some non-empty subset~$I_j\sub I$; note that the sets $I_j$ are pairwise disjoint, though their union is not necessarily all of $I$.  It is then easy to check that
\[
\fb = bafa = \tbinom{F_j\cap \B}{a_j}_{j\in J}.
\]
In fact, it also follows from this that $\fb = (fa)|_\B$; in particular, $\dom(\fb)=\dom(f)\cap \B$.
For the next statement and proof, recall the parameters $\be$, $\lam_i$, $\Lam_J$, from \eqref{eq:cPT}.

%Recall that $a=\binom{A_i}{a_i}_{i\in I}$.  
%To simplify the following statement, we introduce some further notation concerning the map $a=\binom{A_i}{a_i}$.  For each $i\in I$, we write $\lam_i=|A_i|$.  For a subset $J\sub I$, we write $\Lam_J=\prod_{j\in J}\lam_j$.  
%If $X$ is a set and $\mu\leq|X|$ a cardinal, we write $\stirlingii X\mu = \set{A\sub X}{|A|=\mu}$.

%\newpage

\begin{thm}\label{thm:inflation_PT}
Let $f=\tbinom{F_j}{f_j}_{j\in J}\in P^a$, keep the notation of the previous paragraph, and write $\mu=\rank(f)$.~  %Also let $\mu=|J|=\rank(f)$ and $\be=|X\sm B|$.
%$\ga=\big|X\sm(\im(a)\cap\dom(f))\big|$.
\bit
\itemit{i} $\Rh_f^a$ is the union of $(\mu+1)^\be$\ $\gRa$-classes of $P^a$.
\itemit{ii} $\Lh_f^a$ is the union of $\Lam_J$ $\gLa$-classes of $P^a$.
\itemit{iii} $\Hh_f^a$ is the union of $(\mu+1)^\be\Lam_J$ $\gHa$-classes of $P^a$, each of which has size $\mu!$.  
%The map $\phi:P^a\to\T_A$ is injective when restricted to any $\gHa$-class of $P^a$.
\itemit{iv} If $H_{\fb}$ is a non-group $\gH$-class of $\PT_\B$, then each $\gHa$-class of $P^a$ contained in $\Hh_f^a$ is a non-group.
\itemit{v} If $H_{\fb}$ is a group $\gH$-class of $\PT_\B$, then each $\gHa$-class of $P^a$ contained in $\Hh_f^a$ is a group isomorphic to the symmetric group $\S_\mu$.  Further, $\Hh_f^a$ is a $(\mu+1)^\be\times\Lam_J$ rectangular group over $\S_\mu$.
\itemit{vi} $\Dh_f^a=D_f^a=D_\mu^a=\set{g\in P^a}{\rank(g)=\mu}$ is the union of:
\begin{itemize}
\itemit{a} $(\mu+1)^\be S(\al+1,\mu+1)$ $\gRa$-classes of $P^a$,
\itemit{b} $\displaystyle\sum_{K\sub I \atop |K|=\mu}\Lam_K$ $\gLa$-classes of $P^a$,
\itemit{c} $\displaystyle(\mu+1)^\be S(\al+1,\mu+1)\sum_{K\sub I \atop |K|=\mu}\Lam_K$ $\gHa$-classes of $P^a$.
\eit
\eit
\end{thm}

\pf (i).  An $\R^a$-class $R_g^a$ contained in $\Rh_f^a$ is determined by the common domain and kernel of all its members---namely, $D=\dom(g)$ and $K=\ker(g)$---and, since $\gb\R\fb$ in $\PT_\B$, these are constrained so that $\dom(\gb)=\dom(\fb)$ and $\ker(\gb)=\ker(\fb)$.  The former tells us that $\dom(g)\cap \B=\dom(f)\cap \B$, and the latter that the $\ker(\gb)$-classes are precisely $F_j\cap \B$ ($j\in J$).  Because of this, we may write $g=\binom{G_j}{g_j}_{j\in J}$, where $G_j\cap \B=F_j\cap \B=\set{a_i}{i\in I_j}$ for each $j\in J$ (but note that we do not necessarily have $g_j\in A_j$, although~$g_j$ does belong to some set $A_k$).  So, to specify such a domain-kernel pair $(D,K)$, we note that 
\bit
\item $D\cap \B=\bigcup_{j\in J}(F_j\cap \B)=\dom(f)\cap \B$, and that 
\item each $K$-class contains precisely one of the sets $F_j\cap \B$.  
\eit
To complete the description of $(D,K)$, each element of $X\sm \B$ can either belong to $X\sm D$ or else it must be assigned to the $K$-class containing $F_j\cap \B$ for some $j\in J$.  So there are $(\mu+1)^\be$ such pairs $(D,K)$, and the proof of (i) is complete.

\pfitem{ii}  An $\L^a$-class $L_g^a$ contained in $\Lh_f^a$ is determined by the common image of all its members---namely, $\im(g)$---and, since $\gb\L\fb$ in $\PT_\B$, these are constrained so that $\im(\gb)=\im(\fb)$.  Because of this, we may write $g=\binom{G_j}{g_j}_{j\in J}$, where $g_j\in A_j$ for each $j\in J$.  In other words, $\im(g)$ must be a cross-section of the partition $\set{A_j}{j\in J}$, and there are $\Lam_J=\prod_{j\in J}\lam_j$ such cross-sections.

\pfitem{iii}  The number of $\gHa$-classes in $\Hh_f^a$ follows from (i) and (ii).  Each such $\gHa$-class has the same size as~$H_f^a$; but $H_f^a=H_f$ by Theorem \ref{thm:green_PTXYa}(iii), and $|H_f|=\mu!$ by Corollary \ref{cor:Green_sizes_PT}(iv).

\pfitem{iv) and (v}   These follow from \cite[Theorem 2.14(ii) and (iii)]{Sandwiches1}, and the fact that the $\H$-class of a rank $\mu$ idempotent from $\PT_\B$ is isomorphic to $\S_\mu$.

\pfitem{vi}  The number of $\gRh^a$-classes in $D_\mu^a$ is equal to the number of $\R$-classes in the $\D$-class $D_\mu(\PT_\B)={\set{g\in\PT_\B}{\rank(g)=\mu}}$ of $\PT_{\B}$, which is equal to $S(\al+1,\mu+1)$, by Corollary \ref{cor:Green_sizes_PT}(i); (a) now follows from (i).  Part (b) follows from (ii) and the fact that the $\L$-classes contained in $D_\mu(\PT_\B)$ (and hence the $\gLh^a$-classes contained in $D_\mu^a$) are indexed by the subsets of~$\B$ (equivalently, of $I$) of size $\mu$.  Finally, (c) follows immediately from (a) and~(b).~\epf

\begin{rem}
The reader may compare the values given in Theorem \ref{thm:inflation_PT} with the egg-box diagrams in Figure~\ref{fig:PT_9}.
\end{rem}

\begin{rem}
Again, some simplifications occur in the special cases enumerated in Section \ref{sect:structurePT}.  If $a$ is injective, then $\Lam_J=1$ for all $J\sub I$, so the rectangular groups in Theorem \ref{thm:inflation_PT}(v) are $(\mu+1)^\be\times1$ dimensional (see the third, fourth and seventh diagrams in Figure~\ref{fig:PT_9}).  If $a$ is surjective, then $\be=0$, so $(\mu+1)^\be=1$ for all $\mu$, and so the rectangular groups are $1\times\Lam_J$ dimensional (see the fifth, sixth and seventh diagrams in Figure~\ref{fig:PT_9}).  If $a$ is injective and surjective (but not necessarily full), then the rectangular groups are $1\times1$ dimensional: that is, they are groups (see the seventh diagram in Figure~\ref{fig:PT_9}).
\end{rem}

As an immediate consequence of Proposition \ref{prop:regular_Da_classes} and Theorem \ref{thm:inflation_PT}(vi)(c), we may deduce a formula for the size of $P^a=\Reg(\PTXYa)$.

\begin{cor}\label{cor:size_P_PT}
We have $\displaystyle |P^a| = \sum_{\mu=0}^\al \mu!(\mu+1)^\be S(\al+1,\mu+1)\sum_{K\sub I \atop |K|=\mu}\Lam_K$.  \epfres
%Moreover, if $\al\geq\aleph_0$, then $|P^a|=2^{|X|}\Lam_I$.
%\begin{itemize}
%\itemit{i} $\displaystyle \big|\Reg(\PTXYa)\big| = \sum_{\mu=0}^\al \mu!(\mu+1)^\be S(\al+1,\mu+1)\sum_{K\sub I \atop |K|=\mu}\Lam_K$. 
%%
%\itemit{ii} If $\al\geq\aleph_0$, then $\displaystyle \big|\Reg(\PTXYa)\big| =\al^{|X|}\Lam_I$.
%\end{itemize}
\end{cor}

The expression for $|P^a|$ in Corollary \ref{cor:size_P_PT} is valid whether $P^a$ is finite or infinite, but may be simplified substantially in the infinite case.  The next result gives such simplifications, and also provides necessary and sufficient conditions for $P^a$ to be finite, countably infinite or uncountable; these will be of use shortly, when we calculate the \emph{rank} of $P^a$.  The reader may refer to \cite[Chapter 5]{Jech2003} for some background on basic cardinal arithmetic.

%First, if $\al=\rank(a)=0$, then $|P^a|=1$; this follows from the fact that $P^a$ only contains the empty map, and also by applying Corollary~\ref{cor:size_P_PT} (note that $\Lam_\emptyset=1$, as an empty product).  Note also that if $|X|<\aleph_0$ (which implies $\al<\aleph_0$) and if $\lam_i<\aleph_0$ for all $i\in I$, then Corollary~\ref{cor:size_P_PT} gives $|P^a|<\aleph_0$, even if $|Y|\geq\aleph_0$; in this case, Corollary \ref{cor:size_P_PT} cannot be simplified (as far as we are aware).
%%
%The next result concerns the remaining situations.
%%If $\al\geq1$ and either $X$ or $Y$ is infinite, then we have the following:

\begin{prop}\label{prop:size_P_PT}
\begin{itemize}
%\itemit{i} If $\al=0$ or if $|X|<\aleph_0$ and $\lam_i<\aleph_0$ for all $i\in I$, then $|P^a|<\aleph_0$.
%\itemit{i} $|P^a|<\aleph_0 \iff [\al=0$ or $[|X|<\aleph_0$ and $\lam_i<\aleph_0$ for all $i\in I]]$.
%\itemit{ii} If $\al\geq1$ and $|X|\geq\aleph_0$, then $|P^a|=2^{|X|}\Lam_I$.
%\itemit{iii} If $\al\geq1$, $|X|<\aleph_0$, and $\lam_i\geq\aleph_0$ for some $i\in I$, then $|P^a|=\Lam_I=\displaystyle\max_{i\in I}\lam_i$.
\itemit{i} If $\al\geq1$ and $|X|\geq\aleph_0$, then $|P^a|=2^{|X|}\Lam_I=\max(2^{|X|},\Lam_I)$.
\itemit{ii} If $\al\geq1$, $|X|<\aleph_0$, and $\lam_i\geq\aleph_0$ for some $i\in I$, then $|P^a|=\Lam_I=\displaystyle\max_{i\in I}\lam_i$.
\itemit{iii} $|P^a|<\aleph_0 \iff \al=0$ or $[|X|<\aleph_0$ and $\lam_i<\aleph_0$ for all $i\in I]$.
\itemit{iv} $|P^a|=\aleph_0 \iff \al\geq1$, $|X|<\aleph_0$ and $\max_{i\in I}\lam_i=\aleph_0$.
\itemit{v} $|P^a|>\aleph_0 \iff \al\geq1$ and $[|X|\geq\aleph_0$ or $\lam_i>\aleph_0$ for some $i\in I]$.
%\itemit{v} $|P^a|>\aleph_0 \iff [\al\geq1$ and $|X|\geq\aleph_0]$ or $[\al\geq1$ and $\lam_i>\aleph_0$ for some $i\in I]$.
\end{itemize}
\end{prop}

%
%\begin{prop}\label{prop:size_P_PT}
%Suppose $\al=\rank(a)\geq1$.  Then
%\[
%|P^a| = \begin{cases}
%2^{|X|}\Lam_I &\text{if $|X|\geq\aleph_0$}\\
%\displaystyle\max_{i\in I}\lam_i &\text{if $|X|<\aleph_0$ and $\lam_i\geq\aleph_0$ for some $i\in I$.}
%\end{cases}
%\]
%\end{prop}

\pf %(i).  The reverse implication follows quickly from Corollary \ref{cor:size_P_PT}, noting that $|X|<\aleph_0 \implies\al<\aleph_0$.  The forwards implication follows from (ii) and (iii), which we now prove.  
%
%Part (i) follows immediately from Corollary \ref{cor:size_P_PT}, noting that $|X|<\aleph_0 \implies\al<\aleph_0$.  For the remaining parts, 
For any cardinal $\mu$ with $0\leq\mu\leq\al$, write $t_\mu$ for the $\mu$th term in the sum from Corollary~\ref{cor:size_P_PT}.

\pfitem{i}  Suppose $\al\geq1$ and $|X|\geq\aleph_0$.  For any $0\leq\mu\leq\al$, we have:
\bit
\item $\mu!\leq\al!\leq2^{|X|}$ (this is obvious if $\al<\aleph_0$, while if $\al\geq\aleph_0$, then $\al!=2^\al\leq2^{|X|}$), 
\item $(\mu+1)^\be\leq(|X|+1)^{|X|}=|X|^{|X|}=2^{|X|}$, 
\item $S(\al+1,\mu+1)\leq2^{|X|}$ (this is again obvious if $\al<\aleph_0$, while if $\al\geq\aleph_0$, then $S(\al+1,\mu+1)=1$ or $2^\al$), 
\item $\displaystyle\sum_{K\sub I \atop |K|=\mu}\Lam_K \leq \sum_{K\sub I}\Lam_I = 2^{|I|}\Lam_I=2^\al\Lam_I\leq2^{|X|}\Lam_I$.
\eit
Together, these give $t_\mu \leq 2^{|X|}\Lam_I$.  It follows that $|P^a|=\sum_{\mu=0}^\al t_\mu\leq2^{|X|}\Lam_I$, since there are fewer than $2^\al\leq2^{|X|}$ terms in the sum.  To prove the reverse inequality, note that
\[
|P^a|\geq t_\al = \al!(\al+1)^\be S(\al+1,\al+1) \sum_{K\sub I \atop |K|=\al}\Lam_K \geq \al!2^\be\Lam_I.
\]
If $\al<\aleph_0$, then $\be=|X|$ and $\al!<\aleph_0$, so $\al!2^\be\Lam_I=2^{|X|}\Lam_I$.  If $\al\geq\aleph_0$, then
\[
\al!2^\be\Lam_I = 2^\al2^\be\Lam_I = 2^{\al+\be}\Lam_I = 2^{|X|}\Lam_I.
\]

\pfitem{ii}  Suppose $\al\geq1$, $|X|<\aleph_0$, and $\lam_i\geq\aleph_0$ for some $i\in I$.  Since $|X|<\aleph_0$, it follows that $|I|=\al<\aleph_0$, and so $\Lam_I = \prod_{i\in I}\lam_i = \max_{i\in I}\lam_i$.  Now, we immediately have $|P^a|\geq t_\al \geq \Lam_I$.  Conversely, for any $0\leq\mu\leq\al$, $t_\mu$ is a finite multiple of 
\[
\sum_{K\sub I \atop |K|=\mu}\Lam_K \leq \sum_{K\sub I}\Lam_I = 2^\al \Lam_I=\Lam_I,
\]
so that $t_\mu\leq\Lam_I$ for all $\mu$, whence $|P^a|=\sum_{\mu=0}^\al t_\mu\leq (\al+1)\Lam_I=\Lam_I$. 

\pfitem{iii}  If $\al=0$, then $|P^a|=1$.  If $|X|<\aleph_0$ and $\lam_i<\aleph_0$ for all $i\in I$, then $\al,\be\leq|X|<\aleph_0$, so $|P^a|<\aleph_0$ follows from Corollary \ref{cor:size_P_PT}.  Conversely, suppose $|P^a|<\aleph_0$ and $\al\geq1$.  By (i), we must have $|X|<\aleph_0$.  By~(ii), it then follows that $\lam_i<\aleph_0$ for all $i\in I$.

\pfitem{iv}  The reverse implication follows immediately from (ii).  Conversely, suppose $|P^a|=\aleph_0$.  We obviously must have $\al\geq1$.  Then by (i), we have $|X|<\aleph_0$.  By (iii), it follows that $\lam_i\geq\aleph_0$ for some $i\in I$.  Then (ii) gives $\max_{i\in I}\lam_i=|P^a|=\aleph_0$.

\pfitem{v}  Suppose $\al\geq1$.  If $|X|\geq\aleph_0$, then (i) gives $|P^a|\geq2^{|X|}>\aleph_0$.  If $|X|<\aleph_0$ but $\lam_i>\aleph_0$ for some $i\in I$, then (ii) gives $|P^a|>\aleph_0$.  Conversely, suppose $|P^a|>\aleph_0$.  Again it is clear that $\al\geq1$.  Suppose $|X|<\aleph_0$; the proof will be complete if we can show that $\lam_i>\aleph_0$ for some $i\in I$.  But part (iii) gives $\lam_i\geq\aleph_0$ for some $i\in I$.  Together with (iv), which gives $\max_{i\in I}\lam_i\not=\aleph_0$, it follows that $\lam_i>\aleph_0$ for some $i\in I$. \epf

%(iii)--(v).  These all follow quickly from Corollary \ref{cor:size_P_PT} and parts (i) and (ii). \epf

%The reverse implication follows quickly from Corollary \ref{cor:size_P_PT}, noting that $|X|<\aleph_0 \implies\al<\aleph_0$, and that $\Lam_\emptyset=1$ (as an empty product).  The forwards implication follows from (i) and (ii).
%
%\epf

%%\[
%%\xi_\mu = \mu!(\mu+1)^\be S(\al+1,\mu+1)\sum_{K\sub I \atop |K|=\mu}\Lam_K.
%%\]
%Then for any such $\mu$, we have 
%\bit
%\item $\mu!\leq\al!=2^\al$, 
%\item $(\mu+1)^\be\leq(\al+1)^\be=\al^\be$, 
%\item $S(\al+1,\mu+1)\leq2^\al$ (the value is either $1$ or $2^\al$), and 
%\item $\displaystyle\sum_{K\sub I \atop |K|=\mu}\Lam_K \leq \sum_{K\sub I}\Lam_I = 2^{|I|}\Lam_I=2^\al\Lam_I$.
%\eit
%Together, these give $\xi_\mu \leq 2^\al\cdot\al^\be\cdot2^\al\cdot2^\al\Lam_I = \al^\al\cdot\al^\be\Lam_I = \al^{\al+\be}\Lam_I = \al^{|X|}\Lam_I = 2^{|X|}\Lam_I$.
%Since there are fewer than $2^\al$ terms in the sum \eqref{eq:size_P_PT}, it follows that $|P^a|=\sum_{\mu=0}^\al\xi_\mu\leq2^{|X|}\Lam_I$.
%%
%To establish the reverse inequality, one may show (using similar ideas to above) that $\xi_\al \geq 2^{|X|}\Lam_I$.  \epf

\begin{rem}\label{rem:size_P_PT}
The parameter $\Lam_I=\prod_{i\in I}\lam_i$ appears in the expressions for $|P^a|$ in parts (i) and (ii) of Proposition \ref{prop:size_P_PT}.  Note that $\Lam_I$ could be anything from $1$ (which occurs if and only if $a$ is injective) up to $|Y|^\al$ (which occurs, for example, when $\lam_i=|Y|$ for all $i\in I$).  Note also that if $\al\geq\aleph_0$, then the set $\set{\lam_i}{i\in I}$ may or may not have a maximum element.
\end{rem}

%
%\begin{rem}\label{rem:size_P_PT}
%Note that if $\al=\rank(a)\geq\aleph_0$, then $P^a=\Reg(\PTXYa)$ is uncountable, as may be seen by examining the $\mu=\al$ term in the sum from Corollary \ref{cor:size_P_PT}.  The assertion that $P^a$ is uncountable also follows from the fact that $\binom{a_j}{b_j}_{j\in J}$ is regular (indeed, an idempotent) for any subset $J\sub I$, so that $|P^a|\geq2^{|I|}=2^\al$.  From this, it follows that $\rank(P^a)=|P^a|$ if $\al\geq\aleph_0$.  

As mentioned above, we now wish to calculate the rank of $P^a=\Reg(\PTXYa)$.  If $|P^a|>\aleph_0$, then of course ${\rank(P^a)=|P^a|}$.  In order to deal with the case that $|P^a|\leq\aleph_0$, we first show that $P^a$ is MI-dominated (as defined in \cite[Section 3]{Sandwiches1}); this is true regardless of the values of the parameters $\al,\be,|X|,|Y|,\lam_i$, but we also show that RP-domination depends on $\al=\rank(a)$.  %The next lemma is the key step.

\begin{prop}\label{prop:LMC_PT}
\begin{itemize}
\itemit{i} The regular semigroup $P^a=\Reg(\PTXYa)$ is MI-dominated.
\itemit{ii} The regular semigroup $P^a=\Reg(\PTXYa)$ is RP-dominated if and only if $\rank(a)<\aleph_0$.
\end{itemize}
%If $f\in P^a$, then $f=g\star_af\star_ah$ for some $g,h\in \Hh_b^a$.  Consequently, $P^a$ is MI-dominated.
\end{prop}

\pf (i).  By \cite[Propositions 3.3(iii--iv) and 3.5(ii)]{Sandwiches1}, it suffices to show that
\[
P^a=\Hh_b^a\star_aP^a\star_a\Hh_b^a.
\]
With one containment being obvious, let $f\in P^a$.
As before, we may write $f=\binom{F_j}{f_j}_{j\in J}$, where $J\sub I$ and $f_j\in A_j$ for all $j$, and we also write $F_j\cap \B=\set{a_i}{i\in I_j}$, where $I_j\sub I$ is non-empty, noting that the $I_j$ are pairwise disjoint.  Put $K=I\sm\bigcup_{j\in J}I_j$.  For each $j\in J$, fix some partition $\set{F_{j,i}}{i\in I_j}$ of $F_j$, so that $a_i\in F_{j,i}$ for each $i\in I_j$.  Put
\[
g = \tran{F_{j,i}&a_k}{b_i&b_k}_{j\in J,\ i\in I_j,\ k\in K}
\AND
h = \tran{a_j&a_l}{f_j&b_l}_{j\in J,\ l\in I\sm J}.
\]
Then one easily checks that $f=gafah$ and $\gb=\hb=\id_\B$; the latter gives $g,h\in\Hh_b^a$. 

\pfitem{ii}  By \cite[Proposition 3.4]{Sandwiches1}, $P^a$ is RP-dominated if and only if the local monoid $e\star_aP^a\star_ae$ is factorisable for all $e\in\MI(P^a)$.  By \cite[Proposition 3.8]{Sandwiches1}, $e\star_aP^a\star_ae$ is isomorphic to $W=(a\PTXY a,\starb)$ for any $e\in\MI(P^a)$, and we observed in Section \ref{sect:structurePT} that $W$ is isomorphic to $\PT_\B$.  But $\PT_\B$ is factorisable if and only if $\B$ is finite, by \cite[Theorem 3.1]{Tirasupa1979}.  \epf

We may now give the rank of $P^a=\Reg(\PTXYa)$.

\begin{thm}\label{thm:rank_P_PT}
\begin{itemize}
\itemit{i} If $|P^a|\geq\aleph_0$, then $\rank(P^a)=|P^a|$.  (Proposition \ref{prop:size_P_PT} gives necessary and sufficient conditions under which $|P^a|\geq\aleph_0$, and formulae for $|P^a|$ in these cases.)
\itemit{ii} If $|P^a|<\aleph_0$ (so that $\al=0$ or $[|X|<\aleph_0$ and $\lam_i<\aleph_0$ for all $i\in I]$, by Proposition \ref{prop:size_P_PT}), then
\[
\rank(P^a) = 
\begin{cases}
1 &\text{if $\al=0$}\\
%x_\al+\max((\al+1)^\be,\Lam_I,y_\al) &\text{if $\al\geq1$,}
1+\max(2^\be,\Lam_I) &\text{if $\al=1$}\\
2+\max(3^\be,\Lam_I) &\text{if $\al=2$}\\
2+\max((\al+1)^\be,\Lam_I,2) &\text{if $\al\geq3$.}
\end{cases}
\]
%where $x_1=y_1=y_2=1$, and $x_j=y_k=2$ for $j\geq2$ and $k\geq3$.
\end{itemize}
%If $\al=\rank(a)<\aleph_0$, then the rank of $P^a=\Reg(\PTXYa)$ is given by
%\[
%\epfreseq
%\rank(P^a) = \rank(\PT_\al:\S_\al) + \max\big( (\al+1)^\be, \Lam_I, \rank(\S_\al) \big).
%\]
\end{thm}

\pf We have already noted that $\rank(P^a)=|P^a|$ if $|P^a|>\aleph_0$.  Also, if $\al=0$, then clearly $\rank(P^a)=|P^a|=1$.  For the remainder of the proof, we assume that $|P^a|\leq\aleph_0$ and $\al\geq1$.  It follows from Proposition~\ref{prop:size_P_PT}(v) that $|X|<\aleph_0$ (and so also $\al,\be<\aleph_0$) and that $\lam_i\leq\aleph_0$ for all $i\in I$.  In particular, $W=(a\PTXY a,\starb)\cong\PT_\B\cong\PT_\al$ is finite, and so $W\sm G_W\cong\PT_\al\sm\S_\al$ is an ideal of $W$.  By Theorem~\ref{thm:inflation_PT}(v), $\Hh_b^a$ is an $(\al+1)^\be\times\Lam_I$ rectangular group over the symmetric group $\S_\al$.  It then follows from \cite[Theorem 3.16]{Sandwiches1} and Proposition \ref{prop:LMC_PT}(i) that
\[
%\begin{equation}\label{eq:rank_P_PT}
\rank(P^a) = \rank(\PT_\al:\S_\al) + \max\big( (\al+1)^\be, \Lam_I, \rank(\S_\al) \big).
%\end{equation}
\]
If $\lam_i=\aleph_0$ for some $i\in I$, then it follows that $\rank(P^a)=\Lam_I=\aleph_0=|P^a|$, completing the proof of (i).  The proof of (ii) concludes with the observation that
\[
\rank(\PT_\al:\S_\al) = \begin{cases}
1 &\text{if $\al=1$}\\
2 &\text{if $2\leq\al<\aleph_0$}
\end{cases}
\AND
\rank(\S_\al) = \begin{cases}
1 &\text{if $\al\leq2$}\\
2 &\text{if $3\leq\al<\aleph_0$}.
\end{cases}
\]
\epf

\begin{rem}\label{rem:rank_Reg_PTXA}
We may deduce formulae for the ranks of $\Reg(\PT(X,\B))$ and $\Reg(\PT(Y,\si))$ as special cases of Theorem \ref{thm:rank_P_PT}.  Taking $a$ to be full, injective and non-surjective, we see that for finite $X$, and for a proper subset $\B\subsetneq X$,
\begin{align*}
\rank(\Reg(\PT(X,\B))) &= \begin{cases}
1 &\text{if $\B=\emptyset$}\\
1+2^{|X|-1} &\text{if $|\B|=1$}\\
2+(|\B|+1)^{|X|-|\B|} &\text{if $|\B|\geq2$.}
\end{cases}
\intertext{Taking $a$ to be surjective (and possibly full and/or injective), we see that for finite $Y$, and for an equivalence relation $\si$ on a subset $\A\sub Y$,}
\rank(\Reg(\PT(Y,\si))) &= \begin{cases}
1 &\hspace{7.25mm} \text{if $\A=\emptyset$}\\
1+\Lam_I &\hspace{7.25mm} \text{if $|\A/\si|=1$}\\
2+\Lam_I &\hspace{7.25mm} \text{if $|\A/\si|=2$}\\
2+\max(\Lam_I,2) &\hspace{7.25mm} \text{if $|\A/\si|\geq3$.}
\end{cases}
\end{align*}
Note that in the last case, $\max(\Lam_I,2)=2$ if $\si$ is the diagonal relation, or $\Lam_I$ otherwise.  %In particular, we have $\rank(\Reg(\PT(Y,\si)))=4$ if $\si$ is the diagonal relation and $|A|\geq3$; this is simply because 
These results are new, as far as the authors are aware.
\end{rem}

%
%\begin{thm}
%If $\al=\rank(a)<\aleph_0$, then the rank of $P^a=\Reg(\PTXYa)$ is given by
%\[
%\epfreseq
%\rank(P^a) = \rank(\PT_\al:\S_\al) + \max\big( (\al+1)^\be, \Lam, \rank(\S_\al) \big).
%\]
%\end{thm}
%
%
%\begin{rem}
%Note that
%\[
%\rank(\PT_\al:\S_\al) = \begin{cases}
%0 &\text{if $\al=0$}\\
%1 &\text{if $\al=1$}\\
%2 &\text{if $2\leq\al<\aleph_0$}
%\end{cases}
%\AND
%\rank(\S_\al) = \begin{cases}
%1 &\text{if $\al\leq2$}\\
%2 &\text{if $3\leq\al<\aleph_0$}.
%\end{cases}
%\]
%Beyond this, it is not necessarily possible to give a more precise value of $\rank(P^a)$.  Of course, if $\al<\aleph_0$ but $X$ and/or $Y$ is infinite, then it is possible for $P^a$ to have infinite rank.  
%\end{rem}

\subsectiontitle{Idempotents and idempotent-generation}\label{sect:EaPTXYa}
%\subsectiontitle{Idempotents and the idempotent-generated subsemigroup $\cE_{XY}^a=\E_a(\PTXYa)$}\label{sect:EaPTXYa}

Again, we fix $a=\binom{A_i}{a_i}\in\PTYX$, and keep the notation of \eqref{eq:aPT}--\eqref{eq:cPT}.  
Recall that
\[
E_a(\PTXYa)=E_a(P^a)=\set{f\in\PTXY}{f=f\star_af}
\]
denotes the set of $\star_a$-idempotents of $\PTXYa$.  We write $\cE_{XY}^a=\E_a(\PTXYa)=\E_a(P^a)$ for the idempotent-generated subsemigroup $\la E_a(\PTXYa)\ra_a$ of $\PTXYa$.  
In this section, we characterise the elements of this idempotent-generated subsemigroup, and calculate its rank and idempotent rank.  First we characterise and enumerate the idempotents themselves.

\begin{prop}\label{prop:EaPTXYa}
\begin{itemize}
\itemit{i} $E_a(\PTXYa) = \set{f\in\PTXY}{(af)|_{\im(f)}=\id_{\im(f)}}$.
\itemit{ii} If $|P^a|\geq\aleph_0$, then $|E_a(\PTXYa)|=|P^a|$.
\itemit{iii} If $|P^a|<\aleph_0$, then $\ds |E_a(\PTXYa)| = \sum_{\mu=0}^\al(\mu+1)^{|X|-\mu}\sum_{J\sub I\atop|J|=\mu}\Lam_J$.
\end{itemize}
\end{prop}

\pf (i).  This is easily checked. %; it also follows quickly from Lemma \ref{lem:EaEb}.

%Let $f\in E_a(\PTXYa)$, and suppose $y\in\im(f)$.  Then $y=xf$ for some $x\in X$, so $y=xf=xfaf=yaf$, giving $(af)|_{\im(f)}=\id_{\im(f)}$.  Conversely, suppose $(af)|_{\im(f)}=\id_{\im(f)}$, and let $x\in\dom(f)$.  Then $xf\in\im(f)$, so $xf=(xf)af$; since also $\dom(faf)\sub\dom(f)$, it follows that $f=faf$ and $f\in E_a(\PTXYa)$.

\pfitem{ii}  Suppose $|P^a|\geq\aleph_0$; in particular, $\al\geq1$.  Since $E_a(\PTXYa)\sub \Reg(\PTXYa)=P^a$, it suffices to show that $|E_a(\PTXYa)|\geq|P^a|$.  %By Proposition \ref{prop:size_P_PT}(i), we must consider two cases.
%
%Since $|P^a|\not=1$, we must also have $\al\geq1$.  
%
Now, $E_a(\PTXYa)\supseteq E_a(\Hh_b^a)=V(a)$, and the latter is an $(\al+1)^\be\times\Lam_I$ rectangular band, by Theorem \ref{thm:inflation_PT}(v).  Thus, $|E_a(\PTXYa)|\geq\Lam_I$.  By Proposition \ref{prop:size_P_PT}(ii), this completes the proof in the case that $|X|<\aleph_0$.  Now suppose $|X|\geq\aleph_0$.  By Proposition \ref{prop:size_P_PT}(i), it remains to show that $|E_a(\PTXYa)|\geq2^{|X|}$.  Fix some $i\in I$.  For any subset $C\sub X\sm\{a_i\}$, $\binom{C\cup\{a_i\}}{b_i}\in\PTXY$ is a $\star_a$-idempotent (of rank $1$).  Since there are $2^{|X\sm\{a_i\}|}=2^{|X|}$ such idempotents, we have $|E_a(\PTXYa)|\geq2^{|X|}$.

\pfitem{iii}  Suppose $|P^a|<\aleph_0$.  If $\al=0$, then the result is clear (noting that $\Lam_\emptyset=1$, as an empty product), so suppose $\al\geq1$.  By Proposition \ref{prop:size_P_PT}(iii), it follows that $|X|<\aleph_0$ and $\lam_i<\aleph_0$ for all $i\in I$.  To specify a $\star_a$-idempotent $f$ of rank $\mu$, we first note (as in the proof of Theorem \ref{thm:inflation_PT}(ii)) that $\im(f)$ is a cross-section of $\set{A_j}{j\in J}$ for some subset $J\sub I$ with $|J|=\mu$.  Once $J$ is chosen, there are then $\Lam_J$ such images.  Writing $\im(f)=\set{f_j}{j\in J}$, where $f_j\in A_j$ for each $j\in J$, the condition $(af)|_{\im(f)}=\id|_{\im(f)}$ gives $f_j=f_jaf=a_jf$; the remaining $|X|-\mu$ elements of $X\sm\set{a_j}{j\in J}$ may be either left out of $\dom(f)$ or else mapped to any of the $\mu$ elements of $\im(f)$.  So, once $\im(f)$ has been chosen, there are $(\mu+1)^{|X|-\mu}$ ways to complete the definition of $f$. \epf

\begin{rem}
If $a$ is surjective, then $\al=|X|=|\A/\si|$ in the formula from Proposition \ref{prop:EaPTXYa}(iii), so
\[
|E(\PT(Y,\si))| = \sum_{\mu=0}^{|\A/\si|}(\mu+1)^{|\A/\si|-\mu}\sum_{J\sub I\atop|J|=\mu}\Lam_J.
\]
If $a$ is injective, then $\ds\sum_{J\sub I\atop|J|=\mu}\Lam_J=\tbinom\al\mu$.  In particular, taking $a$ to be full and injective, we obtain\\[-4truemm]
\[
|E(\PT(X,\B))| = \sum_{\mu=0}^{|\B|}(\mu+1)^{|X|-\mu}\tbinom{|X|}\mu.
\]
Again, these two results appear to be new.  Of course, if $a$ is a full bijection, then both of the above formulae reduce to the well known formula (see for example \cite[Corollary 2.7.5]{GMbook} for the finite case):
\[
|E(\PT_X)| = |E(\PT_\al)| = \sum_{\mu=0}^{\al}(\mu+1)^{\al-\mu}\tbinom\al\mu.
\]
\end{rem}

Next, we need the following result concerning the idempotent-generated subsemigroup $\E(\PT_\B)=\la E(\PT_\B)\ra$ of $\PT_\B$.  For the statement, we need some definitions.  As in \cite{Howie1966}, we define the \emph{shift}, \emph{collapse}, \emph{defect} and \emph{codefect} of a partial transformation $f\in\PT_\B$ by
\[
\sh(f) = \big| \set{b\in\dom(f)}{bf\not=b} \big| \COmma
\col(f) = \sum_{b\in\im(f)}(|bf^{-1}|-1) \COmma
\defect(f)= |\B\sm\im(f)| \COmma
\codefect(f) = |\B\sm\dom(f)|,
\]
respectively.

\begin{prop}\label{prop:PTBSB}
\begin{itemize}
\itemit{i} If $|\B|<\aleph_0$, then $\E(\PT_\B) = \{\id_\B\} \cup (\PT_\B\sm\S_\B)$, and
\[
\rank(\E(\PT_\B))=\idrank(\E(\PT_\B))=\tbinom{\al+1}2+1.
\]
\itemit{ii} If $|\B|\geq\aleph_0$, then $\E(\PT_\B)=\{\id_\B\}\cup\set{f\in\PT_\B\sm\S_\B}{\sh(f)+\codefect(f)<\aleph_0}$
\item[] \hspace{4.8cm} ${}\cup\set{f\in\PT_\B}{\sh(f)+\codefect(f)=\col(f)+\codefect(f)=\defect(f)\geq\aleph_0}$,
\item[] and $\rank(\E(\PT_\B))=\idrank(\E(\PT_\B))=|\PT_\B|=2^{|\B|}$.
\end{itemize}
\end{prop}

\pf Part (i) was apparently first proved in \cite{EP1970}; see also \cite{Garba1990}.  Part (ii) follows quickly from \cite[Theorem~III]{Howie1966}, which gives a similar characterisation of the idempotent-generated subsemigroup of $\T_\B$. \epf

\begin{thm}\label{thm:EXYa}
We have $\EXYa=\E_a(\PTXYa)=\E(\PT_\B)\varphi^{-1}$; see Proposition \ref{prop:PTBSB} for a description of~$\E(\PT_\B)$.  Further,
\[
\rank(\EXYa)=\idrank(\EXYa)= \begin{cases}
|\EXYa|=|P^a| &\text{if $|P^a|\geq\aleph_0$}\\
\binom{\al+1}2+\max\big((\al+1)^\be,\Lam_I\big) &\text{if $|P^a|<\aleph_0$.}
\end{cases}
\]
%\begin{itemize}
%\itemit{i} We have $\EXYa=\E_a(\PTXYa)=\E(\PT_B)\varphi^{-1}$. %; see Proposition \ref{prop:PTBSB} for a description of~$\E(\PT_B)$.
%\itemit{ii} If $|P^a|\geq\aleph_0$, then $\rank(\EXYa)=\idrank(\EXYa)=|P^a|$. %; see 
%\itemit{iii} If $|P^a|<\aleph_0$, then $\rank(\EXYa)=\idrank(\EXYa)=\binom{\al+1}2+\max\big((\al+1)^\be,\Lam_I\big)$.
%\end{itemize}
\end{thm}

\pf The first assertion follows from \cite[Theorem 2.17]{Sandwiches1}.  Since $E_a(\PTXYa)\sub P^a$, and since $P^a$ is a subsemigroup of $\PTXYa$, it follows that $\EXYa\sub P^a$, and so $|E_a(\PTXYa)|\leq|\EXYa|\leq|P^a|$.  It follows from Proposition \ref{prop:EaPTXYa}(ii) that $|\EXYa|=|P^a|$ if $|P^a|\geq\aleph_0$.  It remains to prove the assertion concerning the rank and idempotent rank of $\EXYa$.  This clearly being true if $\al=0$, we assume $\al\geq1$ for the remainder of the proof. 

First suppose $|X|\geq\aleph_0$.  Proposition \ref{prop:size_P_PT}(v) gives $|P^a|>\aleph_0$, and the previous paragraph then gives $|\EXYa|=|P^a|>\aleph_0$.
%
%
%, and so Propositions \ref{prop:size_P_PT} and \ref{prop:EaPTXYa} give 
%\[
%|\EXYa| \geq |E_a(\PTXYa)| = |P_a| >\aleph_0.
%\]
%Since $P^a=\Reg(\PTXYa)$ is a semigroup, and since $E_a(\PTXYa)\sub P^a$, it follows that $\EXYa\sub P^a$, so that $|\EXYa|\leq|P^a|$, giving $|\EXYa|=|P^a|$.
%%
Since $\EXYa$ is uncountable, its rank and idempotent rank are equal to its size.

Now suppose $|X|<\aleph_0$, noting that $\al<\aleph_0$ also.  Together with the fact that ${W=(a\PTXY a,\starb)\cong\PT_\B}$,~\cite[Theorem~3.17]{Sandwiches1} and Propositions \ref{prop:LMC_PT}(i) and \ref{prop:PTBSB}(i) then give
\begin{align}
\tag*{} \rank(\EXYa)=\idrank(\EXYa)&=\rank(\E(\PT_\B))+\max\big((\al+1)^\be,\Lam_I\big)-1\\
\label{eq:EXYa}&=\tbinom{\al+1}2+\max\big((\al+1)^\be,\Lam_I\big).
\end{align}
%Proposition \ref{prop:PTBSB}(i) concludes 
This completes the proof in the case of $|P^a|<\aleph_0$.  If $|P^a|\geq\aleph_0$, then (keeping in mind that $|X|<\aleph_0$) Proposition \ref{prop:size_P_PT}(iv) and (v) give $\lam_i\geq\aleph_0$ for some $i\in I$; since all terms in \eqref{eq:EXYa} other than $\Lam_I$ are finite, it follows that $\rank(\EXYa)=\Lam_I=|P^a|$, using Proposition \ref{prop:size_P_PT}(ii) in the last step. \epf

In the case that $\al<\aleph_0$, $\EXYa$ has a very neat alternative description.  Before we can give it (in Theorem~\ref{thm:finite_EXYa_PT}), we require the following lemma, which will also be useful on a number of other occasions.  Recall that $\xi=\max(|X|,|Y|)$.

\begin{lemma}\label{lem:maximal_J_PT}
If $\al<\aleph_0$, then $J_b^a=D_\al^a=\Hh_b^a$.  If also $\al=\xi$, then $J_b^a$ is a maximum $\J^a$-class of $\PTXYa$.
\end{lemma}

\pf Since $\al<\aleph_0$, it follows that $W=(a\PTXY a,\starb)\cong\PT_\al$ is finite.  So \cite[Lemmas 3.6 and 3.7]{Sandwiches1} give ${\Hh_b^a=J_b^a=D_b^a=D_\al^a}$.  If also $\al=\xi$, then Proposition \ref{prop:maximal_J_PT}(ii) gives the second assertion. \epf

\begin{thm}\label{thm:finite_EXYa_PT}
If $\al=\rank(a)<\aleph_0$, then $\EXYa=\E_a(\PTXYa) = E_a(D_\al^a) \cup (P^a\sm D_\al^a)$.
\end{thm}

\pf This follows quickly from $\EXYa=\E(\PT_\B)\varphi^{-1}=\id_\B\varphi^{-1}\cup(\PT_\B\sm\S_\B)\varphi^{-1}$, and the fact that $D_\al^a=\Hh_b^a$ (which itself follows from Lemma \ref{lem:maximal_J_PT}, as $\al<\aleph_0$). \epf

\begin{rem}
Again, by considering various special cases, Theorems \ref{thm:EXYa} and \ref{thm:finite_EXYa_PT} yield results on the idempotent-generated subsemigroups of the semigroups $\PT(X,\B)$ and $\PT(Y,\si)$.  We will not state these here; the reader may supply the details if they wish.
\end{rem}

\subsectiontitle{The rank of a sandwich semigroup $\PTXYa$}\label{sect:rank_PTXYa}

Again, we fix $a=\binom{A_i}{a_i}\in\PTYX$, and keep the notation of \eqref{eq:aPT}--\eqref{eq:cPT}.  
In this section, we calculate the rank of a sandwich semigroup $\PTXYa$.  
%
%The value of $\rank(\PTXYa)$ depends crucially on whether the sandwich element $a$ is injective and/or full and/or surjective; see Theorems~\ref{thm:rankPTXYa1},~\ref{thm:rankPTXYa2} and \ref{thm:rankPTXYa3}.
%
Unlike the case with $P^a=\Reg(\PTXYa)$ and $\EXYa=\E_a(\PTXYa)$, treated in Sections~\ref{sect:RegPTXYa}  and \ref{sect:EaPTXYa}, there is no general result to aid us (apart from \cite[Lemma 5.1]{Sandwiches1}), so our results in this section depend very much on the structure of the category $\PT$; we will also see that they depend heavily on the nature of the sandwich element $a$.

We begin by eliminating some easy special cases:
\bit
\item If $X=\emptyset$ or $Y=\emptyset$, then $\rank(\PTXYa)=|\PTXYa|=1$.
\item If $X,Y\not=\emptyset$ and $\al=0$, then $f\star_ag=\emptyset$ for all $f,g\in\PTXY$, and so $\rank(\PTXYa)=|\PTXYa|-1$.
\item If $X,Y\not=\emptyset$, then since $|\PTXYa|=|\PTXY|=(|Y|+1)^{|X|}$, it follows that $|\PTXYa|>\aleph_0$ if and only if $|X|\geq\aleph_0$ or $|Y|>\aleph_0$.  In such cases, we have $\rank(\PTXYa)=|\PTXYa|$. %=|Y|^{|X|}$.
%\item We also see that $|\PTXYa|=\aleph_0$ if and only if $|X|<\aleph_0$ and $|Y|=\aleph_0$.
\eit

%\medskip

{\bf \boldmath Thus, for the remainder of Section \ref{sect:rank_PTXYa}, we assume that $X,Y$ are non-empty, that $\al\geq1$, and that $|X|<\aleph_0$ and $|Y|\leq\aleph_0$.}

\medskip

The value of $\rank(\PTXYa)$ depends crucially on whether the sandwich element $a$ is injective and/or full and/or surjective; see Table \ref{tab:rankPTXYa} and Theorems~\ref{thm:rankPTXYa1},~\ref{thm:rankPTXYa2} and \ref{thm:rankPTXYa3}.  Note that if $a$ is a (full) bijection, then $\PTXYa\cong\PT_X$ (as discussed in Section \ref{sect:structurePT}), and so (see \cite[Theorem 3.1.5]{GMbook}, for example)
\begin{equation}\label{eq:rankPTX}
\rank(\PTXYa)=\rank(\PT_X)=\begin{cases}
|X|+1 &\text{if $|X|\leq2$}\\
4 &\text{if $|X|\geq3$.}
\end{cases}
\end{equation}
Thus, we will also generally assume that $a$ is not a (full) bijection.

\begin{table}[H]
\begin{center}
{
\begin{tabular}{|c|c|c|c|c|}
\hline
$a$ full? & $a$ injective? & $a$ surjective? & Reference & Egg-box diagram\\
\hline
N&N&N&&Figure \ref{fig:PT_2}\\ %Theorem \ref{thm:rankPTXYa1}\\
Y&N&N&Theorem \ref{thm:rankPTXYa1}&Figure \ref{fig:PT_3}\\
N&Y&N&&Figure \ref{fig:PT_4}\\ %Theorem \ref{thm:rankPTXYa1}\\
\hline
Y&Y&N&Theorem \ref{thm:rankPTXYa2}&Figure \ref{fig:PT_5}\\
\hline
N&N&Y&&Figure \ref{fig:PT_6} (left)\\ %Theorem \ref{thm:rankPTXYa3}\\
Y&N&Y&Theorem \ref{thm:rankPTXYa3}&Figure \ref{fig:PT_6} (middle)\\
N&Y&Y&&Figure \ref{fig:PT_6} (right)\\ %Theorem \ref{thm:rankPTXYa3}\\
\hline
Y&Y&Y&classical; see \eqref{eq:rankPTX} &Figure \ref{fig:PT_1}\\
%\hline
%Y&Y&Y&classical (see above)\\
%Y&Y&N&Theorem \ref{thm:rankPTXYa2}\\
%Y&N&Y&Theorem \ref{thm:rankPTXYa3}\\
%Y&N&N&Theorem \ref{thm:rankPTXYa1}\\
%N&Y&Y&Theorem \ref{thm:rankPTXYa3}\\
%N&Y&N&Theorem \ref{thm:rankPTXYa1}\\
%N&N&Y&Theorem \ref{thm:rankPTXYa3}\\
%N&N&N&Theorem \ref{thm:rankPTXYa1}\\
\hline
\end{tabular}
}
\end{center}
\vspace{-5mm}
\caption{References for $\rank(\PTXYa)$ for all combinations of $a$ being full and/or injective and/or surjective.}
\label{tab:rankPTXYa}
\end{table}

Recall that we write $\xi=\min(|X|,|Y|)$, and note that $\xi\leq|X|<\aleph_0$.
Since $\al=\rank(a)\leq\xi<\aleph_0$, $a$ is stable (by Lemma \ref{lem:stable_PT}(iii)), and so ${\J^a}={\D^a}$ on $\PTXYa$ (by Proposition \ref{prop:JaDaPT}).
Throughout this section, 
%we will also write $\xi=\min(|X|,|Y|)$, noting that $\xi<\aleph_0$, and 
we denote the ${\J}={\D}$-classes of~$\PTXY$ (not to be confused with the ${\J^a}={\D^a}$-classes of $\PTXYa$) by
\[
D_\mu = D_\mu(\PTXY) = \set{f\in\PTXY}{\rank(f)=\mu} \qquad\text{for each $\mu\in\{0,1,\ldots,\xi\}$.}
\]
Note that $\PTXY=D_0\cup D_1\cup\cdots\cup D_\xi$.
Also, since $\al<\aleph_0$, we may assume that $I=\{1,\ldots,\al\}$, and so write
\[
a=\trans{A_1&\cdots&A_\al}{a_1&\cdots&a_\al}
\AND
b=\trans{a_1&\cdots&a_\al}{b_1&\cdots&b_\al}.
\]
%{\red We call $f\in\PT_{CD}$ \emph{full} if $\dom(f)=C$: that is, if $f\in\T_{CD}$.}  
The next two results play a key role in the calculation of $\rank(\PTXYa)$ in all the remaining cases.

\begin{lemma}\label{lem:Dmu_Dmu+1}
\begin{itemize}
\itemit{i} If $\mu\leq\al-2$, then $D_\mu\sub D_{\mu+1}\star_aD_{\mu+1}$.
\itemit{ii} If $a$ is not surjective, then $D_{\al-1}\sub D_\al\star_aD_\al$.
%$\al<\min(|X|,|Y|)$, then $D_{\al-1}\sub D_\al\star_aD_\al$.
\end{itemize}
%If $f\in D_\mu$ where $\mu\leq\al-2$, then $f=g\star_ah$ for some $g,h\in D_{\mu+1}$.
\end{lemma}

\pf We prove both parts together.  Suppose $\mu\leq\al-1$, and let $f=\trans{F_1&\cdots&F_\mu}{f_1&\cdots&f_\mu}\in D_\mu$.  Since $\mu<\al\leq|X|$, one of the following must hold:
\bit
\item[(a)] $f$ is not full, in which case we fix some $x\in X\sm\dom(f)$, or
\item[(b)] $f$ is not injective, in which case, renaming if necessary, we fix a partition $\{F_\mu',F_\mu''\}$ of $F_\mu$.
\eit
Since $\mu<\al\leq|Y|$, we may also choose some $y\in Y\sm\im(f)$.  If $\mu\leq\al-2$, or if $a$ is not surjective, we may choose some $z\in X\sm\{a_1,\ldots,a_{\mu+1}\}$.  
%
%
%
%
%
%(i).  Suppose $\mu\leq\al-2$, and let $f=\trans{F_1&\cdots&F_\mu}{f_1&\cdots&f_\mu}\in D_\mu$.  Then one of the following must hold:
%\bit
%\item[(a)] there exists $x\in X\sm\dom(f)$, or
%\item[(b)] there exists $1\leq j\leq\mu$ such that $|F_j|\geq2$; in this case, we may assume that $j=\mu$, and we split $F_\mu$ into two non-empty disjoint subsets $F_\mu'$ and $F_\mu''$.
%\eit
%Since $\mu\leq\al-2$, we may also choose some $y\in Y\sm\im(f)$.  
Then it is easy to check that $f=g\star_ah$, where
\begin{align*}
g= \begin{cases}
\trans{F_1&\cdots&F_\mu&x}{b_1&\cdots&b_\mu&b_{\mu+1}} &\text{in case (a)}\\[5truemm]
\trans{F_1&\cdots&F_{\mu-1}&F_\mu'&F_\mu''}{b_1&\cdots&b_{\mu-1}&b_\mu&b_{\mu+1}} &\text{in case (b)}
\end{cases}
\AND
h= \begin{cases}
\trans{a_1&\cdots&a_\mu&z}{f_1&\cdots&f_\mu&y} &\text{in case (a)}\\[5truemm]
\trans{a_1&\cdots&a_\mu&a_{\mu+1}&z}{f_1&\cdots&f_\mu&f_\mu&y} &\text{in case (b).}
%\trans{a_1&\cdots&a_\mu&a_{\mu+2}}{f_1&\cdots&f_\mu&y} &\text{in case (a)}\\~\\
%\trans{a_1&\cdots&a_\mu&a_{\mu+1}&a_{\mu+2}}{f_1&\cdots&f_\mu&f_\mu&y} &\text{in case (b).}
\end{cases}
\end{align*}
\epf

%We immediately deduce the following.

\begin{cor}\label{cor:Dmu_Dmu+1}
\begin{itemize}
\itemit{i} We have $D_0\cup D_1\cup\cdots\cup D_\al=\la D_\al\cup D_{\al-1}\ra_a$.
\itemit{ii} If $a$ is non-surjective, then $D_0\cup D_1\cup\cdots\cup D_\al=\la D_\al\ra_a$.
\end{itemize}
\end{cor}

\pf The forward set containments follow immediately from Lemma \ref{lem:Dmu_Dmu+1}, and the reverse containments from the fact that $D_0\cup D_1\cup\cdots\cup D_\al$ is a subsemigroup of $\PTXYa$, which itself follows from Lemma~\ref{lem:PT_prelim}(iv).~\epf

The cases in which $\al<\xi=\min(|X|,|Y|)$ and $\al=\xi$ are different in flavour, so we will treat them separately.  The next two results are required to treat the $\al<\xi$ case.
% (see Theorems \ref{thm:rankPTXYa1}, \ref{thm:rankPTXYa2} and \ref{thm:rankPTXYa3}).  We begin with the former.  

\begin{lemma}\label{lem:rankPTXYa11}
Suppose $\al<\xi$, and that $f\in D_\al$.  If $a,f$ are both non-injective, or both non-full, then $f\in D_{\al+1}\star_aD_{\al+1}$.
\end{lemma}

\pf Write $f=\trans{F_1&\cdots&F_\al}{f_1&\cdots&f_\al}$.  Because $\al<\xi$, we may choose some $x\in X\sm\im(a)$ and $y\in Y\sm\im(f)$.  If $a,f$ are both non-injective, then without loss of generality, we may choose some partition $\{F_\al',F_\al''\}$ of $F_\al$, and some $z\in A_\al\sm\{b_\al\}$.  If $a,f$ are both non-full, then choose some $u\in X\sm\dom(f)$ and $v\in Y\sm\dom(a)$.  Then $f=g\star_ah$, where $g=\trans{F_1&\cdots&F_{\al-1}&F_\al'&F_\al''}{b_1&\cdots&b_{\al-1}&b_\al&z}$ or $\trans{F_1&\cdots&F_\al&u}{b_1&\cdots&b_\al&v}$, respectively, and $h=\trans{a_1&\cdots&a_\al&x}{f_1&\cdots&f_\al&y}$. \epf

\newpage

\begin{lemma}\label{lem:rankPTXYa12}
Suppose $\al<\xi$, and that $f=g\star_ah$, where $f\in D_\al$ and $g,h\in\PTXY$.
\begin{itemize}\begin{multicols}{2}
\itemit{i} If $a$ is injective and $f$ full, then $f\R g$.
\itemit{ii} If $a$ is full and $f$ injective, then $f\R g$.
\end{multicols}\end{itemize}
\end{lemma}

\pf We begin with some observations relevant to both parts.  From $f=gah$, Proposition \ref{prop:GreenPT}(i) gives 
\begin{itemize}\begin{multicols}{2}
\item[(a)] $\dom(f)\sub\dom(g)$, and
\item[(b)] $\ker(f)\supseteq\ker(g)|_{\dom(f)}$.
\end{multicols}\eitmc
From $\al=\rank(f)=\rank(gah)\leq\rank(ah)\leq\rank(a)=\al$, it follows that $\rank(ah)=\rank(a)=\al<\aleph_0$, and so:
\begin{itemize}\begin{multicols}{2}
\item[(c)] $\im(a)\sub\dom(h)$, and
\item[(d)] $h$ acts injectively on $\im(a)$.
\end{multicols}\eitmc
(i).  Suppose $a$ is injective and $f$ full.  Then (a) gives $\dom(f)=\dom(g)=X$ (since $f$ is full), and (b) then gives $\ker(g)\sub\ker(f)$.  By Proposition \ref{prop:GreenPT}(iv), it remains to prove the reverse inclusion.  But by (d) and the fact that $a$ is injective, we obtain
\[
(x,y)\in\ker(f) \implies xf=yf \implies xgah=ygah \implies xga=yga \implies xg=yg \implies (x,y)\in\ker(g).
\]

\pfitem{ii} Suppose $a$ is full and $f$ injective.  By (c), and since $a$ is full,
\[
x\in\dom(g) \implies xg\in Y=\dom(a) \implies xga\in\im(a)\sub\dom(h) \implies x\in\dom(gah)=\dom(f).
\]
Together with (a), this gives $\dom(f)=\dom(g)$.  Combining this with (b) gives $\ker(g)\sub\ker(f)$; the reverse containment follows because $f$ is injective. \epf

We may now give the rank of $\PTXYa$ in the case that $\al<\xi$.

\begin{thm}\label{thm:rankPTXYa1}
Suppose $1\leq\al=\rank(a)<\xi=\min(|X|,|Y|)$, and that $|X|<\aleph_0$ and $|Y|\leq\aleph_0$.  Then
\[
\rank(\PTXYa) = \sum_{\mu=\al+1}^\xi \mu! \tbinom{|Y|}{\mu}S(|X|+1,\mu+1) + 
\begin{cases}
0 &\text{if $a$ is non-injective and non-full}\\[2truemm]
S(|X|,\al) &\text{if $a$ is injective (and non-full)}\\[2truemm]
\binom{|X|}{\al} &\text{if $a$ is full (and non-injective).}
\end{cases}
\]
\end{thm}

\pf Put $M=D_{\al+1}\cup\cdots\cup D_\xi=\set{f\in\PTXY}{\rank(f)>\al}$.  Proposition \ref{prop:maximal_J_PT}(i) says that $\{f\}$ is a maximal $\J^a$-class of $\PTXYa$ for all $f\in M$.  Consequently, $M$ is contained in any generating set for $\PTXYa$.  It follows that $\rank(\PTXYa)=|M|+\rank(\PTXYa:M)$.  Corollary \ref{cor:Green_sizes_PT}(v) gives
\[
|M|=\sum_{\mu=\al+1}^\xi \mu! \tbinom{|Y|}{\mu}S(|X|+1,\mu+1),
\]
so it remains to calculate $\rank(\PTXYa:M)$.  We consider three cases, as indicated by the statement of the theorem.  Since $\al<\xi$, note that $a$ is non-surjective, so Corollary \ref{cor:Dmu_Dmu+1}(ii) gives
\[
\rank(\PTXYa:M) = \min\bigset{|\Om|}{\Om\sub\PTXY,\ D_\al\sub\la M\cup\Om\ra_a}.
\]
{\bf Case 1.}  Suppose first that $a$ is non-injective and non-full.  Let $f\in D_\al$ be arbitrary.  Since $\al<\xi$, $f$ is non-injective or non-full, and it follows from Lemma \ref{lem:rankPTXYa11} that $f\in\la M\ra_a$.  So $\rank(\PTXYa:M)=0$.

\medskip \noindent {\bf Case 2.}  Next suppose that $a$ is injective (and non-full).  We show that $S(|X|,\al)$ is both a lower and upper bound for $\rank(\PTXYa:M)$.  Beginning with the former, suppose $D_\al\sub\la M\cup\Om\ra_a$; we must show that ${|\Om|\geq S(|X|,\al)}$.  Consider some full transformation $f\in D_\al$.  We claim that there exists $g\in\Om$ with $g\R f$.  Indeed, consider an expression $f=g_1\star_a\cdots\star_ag_k$, where $g_1,\ldots,g_k\in M\cup\Om$.  Then $g_1\R f$ follows trivially if $k=1$, or from Lemma \ref{lem:rankPTXYa12}(i) if $k\geq2$;
%If $k=1$, then $f=g_1$ so $g_1\R f$, while if $k\geq2$, then Lemma \ref{lem:rankPTXYa12}(i) gives $g_1\R f$; 
this establishes the claim, since we cannot have $g_1\in M$ (as ${g_1\R f \implies g_1\J f\implies \rank(g_1)=\al}$).  So $|\Om|$ is at least as large as the number of $\R$-classes in $D_\al$ containing full transformations; such an $\R$-class is uniquely determined by a partition of $X$ into $\al$ blocks, and there are $S(|X|,\al)$ of these.

To complete the proof in this case, it remains to show that there exists $\Om\sub\PTXY$ with $D_\al\sub\la M\cup\Om\ra_a$ and $|\Om|=S(|X|,\al)$.  Let $\sE$ be the set of all equivalences on $X$ with $\al$ blocks.  For each $\ve\in\sE$, choose some $f_\ve\in D_\al$ with $\ker(f_\ve)=\ve$ and $\im(f_\ve)=\dom(a)=\{b_1,\ldots,b_\al\}$, and put $\Om=\set{f_\ve}{\ve\in\sE}$.  Now let $g\in D_\al$ be arbitrary.  If $g$ is non-full, then Lemma \ref{lem:rankPTXYa11} gives $g\in\la D_{\al+1}\ra_a\sub\la M\cup\Om\ra_a$, so suppose $g$ is full.  Write $g=\trans{G_1&\cdots&G_\al}{g_1&\cdots&g_\al}$, and put $\ve=\ker(g)\in\sE$.  Relabelling if necessary, we may assume that $f_\ve=\trans{G_1&\cdots&G_\al}{b_1&\cdots&b_\al}$.  Since $\al<\xi$, we may choose some $x\in X\sm\im(a)$ and $y\in Y\sm\im(g)$.  Then $g=f_\ve\star_ah$, where $h=\trans{a_1&\cdots&a_\al&x}{g_1&\cdots&g_\al&y}\in M$.  Since $|\Om|=|\sE|=S(|X|,\al)$, we are done.

\medskip \noindent {\bf Case 3.}  Suppose now that $a$ is full (and non-injective).  The proof that $\rank(\PTXYa:M)=\binom{|X|}\al$ is very similar to Case 2, so we just give the outline.  First, one may show that if $D_\al\sub\la M\cup\Om\ra_a$, then $\Om$ must contain an element from every $\R$-class in $D_\al$ containing injective partial transformations; from this, it quickly follows that $\rank(\PTXYa:M)\geq\binom{|X|}\al$.  To establish the reverse inequality, write $\sF=\set{C\sub X}{|C|=\al}$.  For each $C\in\sF$, choose some $f_C\in D_\al$ with $\dom(f_C)=C$ and $\im(f_C)$ a cross-section of $\ker(a)$, and put $\Om=\set{f_C}{C\in\sF}$.  Then for any $g\in D_\al$, we have $g\in\la D_{\al+1}\ra_a$ if $g$ is non-injective, or else $g=f_C\star_ah$ for some $C\in\sF$ and $h\in M$. \epf

\begin{rem}\label{rem:rankPTXYa1}
If $1\leq\al<|X|<|Y|=\aleph_0$, then Theorem \ref{thm:rankPTXYa1} gives $\rank(\PTXYa)=|\PTXYa|=\aleph_0$.
\end{rem}

Having now covered the case in which $\al<\xi=\min(|X|,|Y|)$, we assume that $\al=\xi$ for the remainder of the section.  We still assume that $|X|<\aleph_0$ and $|Y|\leq\aleph_0$.  If we had $|X|=|Y|$, then $a$ would be a (full) bijection, and so $\PTXYa\cong\PT_X$, with $|X|<\aleph_0$.  Since the rank of finite $\PT_X$ is well known (see above),
% (it is $|X|+1$ for $|X|\leq2$, and $4$ for $|X|\geq3$), 
we assume that either $\al=|X|<|Y|$ or $\al=|Y|<|X|$.  The former implies that~$a$ is a surjection, and the latter that $a$ is a full injection, so that $\PTXYa\cong\PT(Y,\si)$ or $\PTXYa\cong\PT(X,\B)$, respectively, as discussed in Section~\ref{sect:structurePT}.
As noted earlier, the rank of $\PT(X,\B)$ has already been calculated \cite[Theorem~2.4]{FS2014}; however, we include a proof here (of the special case in which $\al=|Y|<|X|$) for completeness, and to show the versatility of our general techniques.
To the knowledge of the authors, the rank of $\PT(Y,\si)$ has not previously been calculated, even in the case that $\si$ is a diagonal relation (on a proper subset of $Y$).

\begin{lemma}\label{lem:al=xi}
\begin{itemize}
\itemit{i} If $\al=|Y|<\aleph_0$, then $P_1^a=\PTXYa$, $P_2^a=P^a$, and ${\R^a}={\R}$ on $\PTXYa$.
\itemit{ii} If $\al=|X|<\aleph_0$, then $P_2^a=\PTXYa$, $P_1^a=P^a$, and ${\L^a}={\L}$ on $\PTXYa$.
%
%then:
%\begin{itemize}\begin{multicols}{3}
%\itemit{a} $P_1^a=\PTXYa$,
%\itemit{b} $P_2^a=P^a$, % is a right ideal,
%\itemit{c} ${\R^a}={\R}$ on $\PTXYa$.
%\end{multicols}\end{itemize}
%\itemit{ii} If $\al=|X|<|Y|\leq\aleph_0$, then:
%\begin{itemize}\begin{multicols}{3}
%\itemit{a} $P_2^a=\PTXYa$,
%\itemit{b} $P_1^a=P^a$, % is a left ideal,
%\itemit{c} ${\L^a}={\L}$ on $\PTXYa$.
%\end{multicols}\end{itemize}
\end{itemize}
\end{lemma}

\pf If $\al=|Y|<\aleph_0$, then $a$ is a full injection, and so right-invertible.  If $\al=|X|<\aleph_0$, then $a$ is a surjection, and so left-invertible.
So the result follows immediately from \cite[Lemma 1.2]{Sandwiches1}. \epf

%We just prove (i), as (ii) is virtually identical.  Suppose $\al=|Y|<|X|\leq\aleph_0$.  Since $\al=|Y|<\aleph_0$, $a$ is a full injection, so (a) follows from Proposition \ref{prop:P_sets_PT}(i).  Part (b) immediately follows from (a), while (c) follows from (a) and Theorem \ref{thm:green_PTXYa}(i). \epf
%
%Then (c) is an immediate consequence of (a) and Theorem \ref{thm:green_PTXYa}(i).  It also follows from~(a) that $P^a=P_1^a\cap P_2^a=P_2^a$.  To show that $P_2^a$ is a right ideal, let $f\in P_2^a$ and $g\in\PTXY$.  The former gives $\im(f)=\im(af)$, by Proposition \ref{prop:P_sets_PT}(ii).  But then
%\[
%\im(fag)=[\im(f)\cap\dom(ag)]ag = [\im(af)\cap\dom(ag)]ag = \im(afag),
%\]
%so that $f\star_ag=fag\in P_2^a$, using Proposition \ref{prop:P_sets_PT}(ii) again.
%
%(ii).  In similar fashion to (i), we quickly obtain (a), (c) and $P_1^a=P^a$.  To show that $P_1^a$ is a left ideal, let $f\in P_1^a$ and $g\in\PTXY$.  The former gives $\dom(f)=\dom(fa)$ and $\ker(f)=\ker(fa)$, by Proposition \ref{prop:P_sets_PT}(i).  Then
%\[
%\dom(gaf)=[\im(ga)\cap\dom(f)](ga)^{-1}=[\im(ga)\cap\dom(fa)](ga)^{-1}=\dom(gafa).
%\]
%This also gives $\ker(gaf)\sub\ker(gafa)$, by Lemma \ref{lem:PT_prelim}(iii).  Conversely,
%\[
%(x,y)\in\ker(gafa) 
%%\implies xgafa=ygafa 
%\implies (xga,yga)\in\ker(fa)=\ker(f) 
%%\implies xgaf=ygaf 
%\implies (x,y)\in\ker(gaf),
%\]
%so that $\ker(gafa)=\ker(gaf)$.  Together with $\dom(gafa)=\dom(gaf)$, it follows that $g\star_af=gaf\in P_1^a$. \epf

\begin{lemma}\label{lem:rank_PTXYa_lower}
Suppose $\al=\xi<\aleph_0$, and that $D_\al\sub\la\Om\ra_a$, where $\Om\sub\PTXY$.
\bit
\itemit{i} $\Om$ contains a cross-section of $D_\al/{\R}$ and a cross-section of $D_\al/{\L}$.
\itemit{ii} If $\al=|Y|<|X|\leq\aleph_0$, then $|\Om\cap D_\al|\geq S(|X|+1,\al+1)$.
\itemit{iii} If $\al=|X|<|Y|\leq\aleph_0$, then $|\Om\cap D_\al|\geq\binom{|Y|}\al$.
\itemit{iv} If $3\leq\al=|X|<|Y|\leq\aleph_0$ and $a$ is injective, then $|\Om\cap D_\al|\geq\binom{|Y|}\al+1$.
%\itemit{iv} If $\al=|X|<|Y|\leq\aleph_0$ and $a$ is injective, then $|\Om\cap D_\al|\geq\binom{|Y|}\al+
%\begin{cases} 
%0 &\text{if $\al\leq2$}\\
%1 &\text{if $\al\geq3$.}
%\end{cases}
%$
\end{itemize}
\end{lemma}

\pf By Corollary \ref{cor:Green_sizes_PT}(i) and (ii), $|D_\al/{\R}|=S(|X|+1,\al+1)$ and $|D_\al/{\L}|=\binom{|Y|}\al$.

\pfitem{i}  This follows immediately from \cite[Lemma 5.1]{Sandwiches1}, since $D_\al=J_b$ is the maximum $\J$-class in $\PTXY$, and since every element of $a\PTXY\cup\PTXY a$ is stable; the latter follows from Lemma \ref{lem:stable_PT}(iii) as $\al=\rank(a)<\aleph_0$.

%Suppose first that $\al=|Y|<|X|$.  Since $|D_\al/{\L}|=\binom\al\al=1$, $\Om$ clearly contains a cross-section of $D_\al/{\L}$.  Now suppose $f\in D_\al$.  We will show that $\Om\cap R_f\not=\emptyset$.  Consider some expression $f=g_1\star_a\cdots\star_ag_k$, where $g_1,\ldots,g_k\in\Om$.  If $k=1$, then $f=g_1\in\Om\cap R_f$, so suppose $k\geq2$.  For simplicity, write $g=g_1$ and $h=g_2\star_a\cdots\star_ag_k$.  Now, $\al=\rank(f)=\rank(gah)\leq\rank(g)\leq\al$, so that $\rank(g)=\rank(f)$.  Thus, $g\J f=g(ah)$.  But $\rank(ah)\leq\rank(a)<\aleph_0$, so $ah$ is ($\R$-)stable, by Lemma \ref{lem:stable_PT}.  It follows that $g\R g(ah)=f$, whence $g=g_1\in\Om\cap R_f$, as required.

%The case in which $\al=|X|<|Y|$ is proved by a dual argument.

\pfitem{ii) and (iii}   These follow immediately from (i), together with the above-mentioned fact about the sizes of $D_\al/{\R}$ and~$D_\al/{\L}$.

\pfitem{iv}  Suppose $3\leq\al=|X|<|Y|\leq\aleph_0$ and $a$ is injective.  It is enough to show that $\Om$ must contain \emph{two} elements from some $\L$-class in $D_\al$.  By Proposition \ref{prop:maximal_J_PT}(ii), $J_b^a$ is a maximal $\J^a$-class of $\PTXYa$.  By Lemma~\ref{lem:maximal_J_PT} and Theorem \ref{thm:inflation_PT}(v), $J_b^a=\Hh_b^a$ is a subsemigroup of $\PTXYa$; thus, since $\Hh_b^a\sub D_\al\sub\la\Om\ra_a$, it follows that $\Om$ contains a generating set for $\Hh_b^a$.  By Theorem \ref{thm:inflation_PT}(v), $\Hh_b^a$ is an $(\al+1)^\be\times\Lam_I$ rectangular group over~$\S_\al$.  Since $\al=|X|$, $a$ is surjective, so $\be=0$ and $(\al+1)^\be=1$.  Since $a$ is injective, $\Lam_I=1$.  Consequently, $\Hh_b^a$ is in fact a group; thus, $\Hh_b^a=H_b^a$ is isomorphic to $\S_\al$.  Since $\Om$ contains a generating set for $\Hh_b^a$, as shown above, and since $\rank(\S_\al)=2$ (as $\al\geq3$), it follows that $\Om$ contains (at least) two elements of $H_b^a\sub L_b^a$, as required. \epf

\begin{rem}
In the proof of (iv), we showed that the maximum $\J^a$-class of $\PTXYa$ is a group if $a$ is injective and $3\leq\al=|X|<|Y|\leq\aleph_0$ (this is still true if $\al\leq2$); see Figure \ref{fig:PT_6} (right).  The identity element of this group (i.e.,~$b$) is not an identity element of $\PTXYa$, although it is a \emph{left} identity.
\end{rem}

Recall that $D_\al^a=D_\al\cap P^a$ is the set of all regular elements of $\PTXYa$ of maximum rank $\al$.  By Lemma~\ref{lem:maximal_J_PT}, $D_\al^a=J_b^a$ is the maximum $\J^a$-class of $\PTXYa$.

\begin{lemma}\label{lem:OmGa}
Suppose $\al=\xi<\aleph_0$, 
%that $D_\al^a=\la\Om\ra_a$, 
and that $\Ga$ is a cross-section of the non-regular $\D^a$-classes contained in~$D_\al$.  Then $D_\al\sub\la D_\al^a\cup\Ga\ra_a$.
%Suppose $\al=\xi<\aleph_0$, and that $D_\al^a=\la\Om\ra_a$.
%\bit
%\itemit{i} If $\Ga$ is a cross-section of the non-regular $\D^a$-classes contained in $D_\al$.  Then $D_\al\sub\la\Om\cup\Ga\ra_a$.
%\eit
\end{lemma}

\pf We must consider two cases, depending on whether $\al=|Y|$ or $|X|$.  We give the proof only in the former case, as the latter is nearly identical.  
%By assumption, we have $D_\al^a\sub\la D_\al^a\cup\Ga\ra_a$, so we just need 
Clearly it suffices to show that $D_\al\sm D_\al^a\sub\la D_\al^a\cup\Ga\ra_a$.  With this in mind, let $f\in D_\al\sm D_\al^a$; we must show that $f\in\la D_\al^a\cup\Ga\ra_a$.  By Lemma \ref{lem:al=xi}(i), $D_\al^a=D_\al\cap P^a=D_\al\cap P_2^a$.  Again using Lemma \ref{lem:al=xi}(i), it follows that $f\in D_\al\sm(D_\al\cap P_2^a)=D_\al\sm P_2^a\sub\PTXYa\sm P_2^a=P_1^a\sm P_2^a$.
%D_\al\cap(\PTXYa\sm P_2^a)=D_\al\cap(P_1^a\sm P_2^a)$.  
By Theorem \ref{thm:green_PTXYa}(iv), $D_f^a=R_f^a$.  Thus, by assumption, there exists some $g\in\Ga$ with $g\R^af$.  If $g=f$, then $f\in\Ga\sub\la D_\al^a\cup\Ga\ra_a$, as desired, so suppose $g\not=f$.  Then $f=g\star_ah$ for some $h\in\PTXY$.  Now, $\al=\rank(f)=\rank(gah)\leq\rank(ah)\leq\rank(h)\leq\al$, so that $\rank(ah)=\rank(h)=\al$.  Now, $\rank(h)=\al$ gives $h\in D_\al$, and $\rank(ah)=\rank(h)$ gives $ah\J h$; stability of $a$ (since $\al<\aleph_0$) then gives $ah\L h$, so that $h\in P_2^a=P^a$.  It follows that $h\in D_\al\cap P^a=D_\al^a$.  Thus, $f=g\star_ah\in\la D_\al^a\cup\Ga\ra_a$, as required. \epf

We now have all we need to calculate $\rank(\PTXYa)$ in the case that $\al=|Y|$.

\begin{thm}\label{thm:rankPTXYa2}
Suppose $1\leq\al=\rank(a)=|Y|<|X|<\aleph_0$.  Then $\rank(\PTXYa) = S(|X|+1,\al+1)$.
\end{thm}

\pf It follows immediately from Lemma \ref{lem:rank_PTXYa_lower}(ii) that $\rank(\PTXYa)\geq S(|X|+1,\al+1)$.  To complete the proof, it suffices to provide a generating set of the specified size.

By Theorem \ref{thm:inflation_PT}(v) and Lemma \ref{lem:maximal_J_PT}, $D_\al^a=\Hh_b^a$ is an $(\al+1)^\be\times\Lam_I$ rectangular group over $\S_\al$.  Since $\al=|Y|$, $a$ is injective, so $\Lam_I=1$.  Since $\al<|X|$, $a$ is not surjective, so $\be\geq1$; together with $\al\geq1$, it follows that $(\al+1)^\be\geq2\geq\rank(\S_\al)$.  Thus, by \cite[Proposition 3.11(i)]{Sandwiches1}, $\rank(D_\al^a)=\max\big((\al+1)^\be,1,\rank(\S_\al)\big)=(\al+1)^\be$.  From \cite[Proposition 3.11(iii)]{Sandwiches1}, we may fix a minimal-size generating set $\Om$ for $D_\al^a$ that is a cross-section of the $\R^a$-classes of $D_\al^a$.  Also, let $\Ga$ be a cross-section of the non-regular $\D^a$-classes in $D_\al$.  By Lemmas \ref{lem:Dmu_Dmu+1}(ii) and~\ref{lem:OmGa}, $\PTXYa=\la D_\al\ra_a\sub\la D_\al^a\cup\Ga\ra_a=\la\la\Om\ra_a\cup\Ga\ra_a=\la\Om\cup\Ga\ra_a$, so it remains to show that $\Om\cup\Ga$ has the desired size.  In fact, by Corollary \ref{cor:Green_sizes_PT}(i), it is enough to show that $\Om\cup\Ga$ is a cross-section of the $\R$-classes in $D_\al$.  By Lemma~\ref{lem:al=xi}(i), ${\R}={\R^a}$ on $\PTXYa$, so every $\R$-class contained in $D_\al$ is an $\R^a$-class.  Now, $\Om$ is a cross-section of the $\R^a$-classes in $D_\al^a$; these are precisely the regular $\R^a$-classes in $D_\al$.  As in the proof of Lemma \ref{lem:OmGa}, $\Ga$ is a cross-section of the non-regular $\R^a$-classes. \epf

\begin{rem}
If $1\leq\al=|Y|<|X|<\aleph_0$, then $a$ is injective, so that $\PTXYa\cong\PT(X,\B)$.  So Theorem~\ref{thm:rankPTXYa2} gives $\rank(\PT(X,\B))=S(|X|+1,|\B|+1)$ if $1\leq|\B|<|X|<\aleph_0$.  This was first proved in \cite[Theorem~2.4]{FS2014}. 
\end{rem}

%Since $\al<|X|$, $a$ is non-surjective, so it follows from Lemma \ref{lem:Dmu_Dmu+1} that any generating set for $\PTXYa$ has size at least $S(|X|+1,\al+1)$.  
%%\[
%%\rank(\PTXYa) = \min \set{|\Om|}{\Om\sub\PTXY,\ D_\al\sub\la\Om\ra_a}.
%%\]
%It follows from Lemma \ref{lem:rank_PTXYa_lower}(ii) that $\rank(\PTXYa)\geq S(|X|+1,\al+1)$.  
%To complete the proof, it suffices to give a generating set for $\PTXYa$ of the required size.  
%
%Let $\Om$ be a minimal-size generating set for $D_\al^a$.  Also, let $\Ga$ be a cross-section of the non-regular $\D^a$-classes in $D_\al$.  By Lemma \ref{lem:OmGa}, $D_\al\sub\la\Om\cup\Ga\ra_a$, and so $\PTXYa=\la\Om\cup\Ga\ra_a$.  So it remains to check that $\Om\cup\Ga$ has the desired size.
%
%Since $D_\al^a=\Hh_b^a$ is a $(\al+1)^\be\times1$ rectangular group over $\S_\al$, and since $\al\geq1$, \cite[Theorem~4.6]{Ruskuc1994} gives
%\[
%|\Om|=\rank(\Hh_b^a)=\max\big((\al+1)^\be,1,\rank(\S_\al)\big)=(\al+1)^\be.
%\]
%
%
%Now, $D_\al = (D_\al\cap P^a) \cup (D_\al\cap(\PTXY\sm P^a)) = D_\al^a \cup (D_\al\cap(\PTXY\sm P^a))$.  
%
%\cite[Theorem~4.6]{Ruskuc1994}

We now focus on the case in which $\al=|X|$.  Note that this implies that $a$ is surjective, that $a$ cannot be both full and injective, that every element of $D_\al$ is full and injective, and that every element of $D_{\al-1}$ is either full or injective (but not both).

\begin{lemma}\label{lem:X<Y1}
Suppose $1\leq\al=|X|<|Y|\leq\aleph_0$.  %, and that $a$ is non-injective.  
\bit
\itemit{i} If $a$ is non-injective, and if $f\in D_{\al-1}$ is full, then $f\in D_\al\star_aD_\al$.
\itemit{ii} If $a$ is non-full, and if $f\in D_{\al-1}$ is injective, then $f\in D_\al\star_aD_\al$.
\itemit{iii} If $a$ is non-injective and non-full, then $\PTXYa=\la D_\al\ra_a$.
%\itemit{iii} If $a$ is injective, if $f\in D_{\al-1}$ is full, and if $f=g\star_ah$ for some $g,h\in\PTXY$, then one of $g,h$ belongs to $D_{\al-1}$.
%\itemit{iv} If $a$ is full, if $f\in D_{\al-1}$ is injective, and if $f=g\star_ah$ for some $g,h\in\PTXY$, then one of $g,h$ belongs to $D_{\al-1}$.
\eit
\end{lemma}

\pf (i).  Write $f=\trans{F_1&\cdots&F_{\al-1}}{f_1&\cdots&f_{\al-1}}$.  Without loss of generality, we may assume that $|F_1|,|A_1|\geq2$.  Since $|X|=\al$, this forces $|F_1|=2$, and $|F_2|=\cdots=|F_{\al-1}|=1$.  Fix some $x\in A_1\sm\{b_1\}$, some $y\in Y\sm\im(f)$, and write $F_1=\{u,v\}$.  Then $f=g\star_ah$, where $g=\trans{u&v&F_2&\cdots&F_{\al-1}}{x&b_1&b_2&\cdots&b_{\al-1}}$ and $h=\trans{a_1&\cdots&a_{\al-1}&a_\al}{f_1&\cdots&f_{\al-1}&y}$. 

\pfitem{ii}  Write $f=\trans{f_1&\cdots&f_{\al-1}}{g_1&\cdots&g_{\al-1}}$, $X\sm\dom(f)=\{x\}$, and choose some $y\in Y\sm\dom(a)$ and some $z\in Y\sm\im(f)$.  Then $f=g\star_ah$, where $g=\trans{f_1&\cdots&f_{\al-1}&x}{b_1&\cdots&b_{\al-1}&y}$ and $h=\trans{a_1&\cdots&a_{\al-1}&a_\al}{g_1&\cdots&g_{\al-1}&z}$. 

\pfitem{iii}  Suppose $a$ is non-injective and non-full.  By Corollary \ref{cor:Dmu_Dmu+1}(i), it suffices to show that $D_{\al-1}\sub\la D_\al\ra_a$.  So suppose $f\in D_{\al-1}$.  If $f$ is injective, then part (ii) gives $f\in\la D_\al\ra_a$.  If $f$ is non-injective, then it is full (since $\rank(f)=\al-1$), so part (i) gives $f\in\la D_\al\ra_a$.
\epf

%\pfitem{iii}  Note that $g,h\in D_{\al-1}\cup D_\al$, by Lemma \ref{lem:PT_prelim}(iv), and that $f$ is non-injective, as $\rank(f)<|\dom(f)|$.  In order to obtain a contradiction, suppose $g,h\in D_\al$.  Since $|X|=\al$, it follows that $g,h$ are (full) injections.  But then $f=gah$ is injective, a contradiction.
%
%\pfitem{iv}  Note that $g,h\in D_{\al-1}\cup D_\al$, by Lemma \ref{lem:PT_prelim}(iv).  In order to obtain a contradiction, suppose $g,h\in D_\al$.  Then $a$ is surjective, and $h$ a full injection, so it follows that $\rank(ah)=\al$.  
%
%
%
%\epf

\begin{lemma}\label{lem:X<Y2}
Suppose $1\leq\al=|X|<|Y|\leq\aleph_0$.  %, and that $a$ is non-injective.  
\bit
\itemit{i} If $a$ is injective, then every element of $\la D_\al\ra_a$ is injective. % and if $f_1,\ldots,f_k\in D_\al$, then $f_1\star_a\cdots\star_af_k$ is injective.
\itemit{ii} If $a$ is full, then every element of $\la D_\al\ra_a$ is full.  %and if $f_1,\ldots,f_k\in D_\al$, then $f_1\star_a\cdots\star_af_k$ is full.
\eit
\end{lemma}

\pf Both statements follow from the facts that: every element of $D_\al$ is a full injection; the composite of injective partial transformations is injective; and the composite of full transformations is full. \epf

%(i).  Since $|X|=\al$, every element of $D_\al$ is injective, and the result quickly follows.  %So, if $a$ is injective, then $f_1\star_a\cdots\star_af_k=f_1af_2a\cdots af_k$ is a composite of injective partial transformations, and so injective itself.
%
%(ii). Since $|X|=\al$, every element of $D_\al$ is full.  So, if $a$ is full, then $f_1\star_a\cdots\star_af_k=f_1af_2a\cdots af_k$ is a composite of full transformations, and so full itself. \epf

\begin{lemma}\label{lem:X<Y3}
Suppose $1\leq\al=|X|<|Y|\leq\aleph_0$.  %^, and put $g=\trans{a_1&\cdots&a_{\al-1}}{b_1&\cdots&b_{\al-1}}$.  %If $f\in D_{\al-1}$ is injective or full, then $f\in D_\al\star_ag\star_aD_\al$.
\bit
\itemit{i} If $f\in D_{\al-1}$ is injective, then $f\in D_\al\star_ag\star_aD_\al$, where $g=\trans{a_1&\cdots&a_{\al-1}}{b_1&\cdots&b_{\al-1}}$.
\itemit{ii} If $f\in D_{\al-1}$ is full, then $f\in D_\al\star_ag\star_aD_\al$, where $g=\trans{a_1&\cdots&a_{\al-2}&\{a_{\al-1},a_\al\}}{b_1&\cdots&b_{\al-2}&b_{\al-1}}$.
\eitres
%If $f,g$ are both full, or both injective, then $f\in D_\al\star_ag\star_aD_\al$.
%%%, and that $a$ is non-injective.  
%%\bit
%%\itemit{i} If $a$ is injective, and if $f,g\in D_{\al-1}$ are full, then $f=h_1\star_ag\star_ah_2$ for some $h_1,h_2\in D_\al$.
%%\itemit{ii} If $a$ is full, and if $f,g\in D_{\al-1}$ are injective, then $f=h_1\star_ag\star_ah_2$ for some $h_1,h_2\in D_\al$.
%%\eit
\end{lemma}

%\pf We prove both parts together.  Let $f\in D_{\al-1}$.  If $f$ is injective, write $f=\trans{f_1&\cdots&f_{\al-1}}{g_1&\cdots&g_{\al-1}}$, and let $X\sm\dom(f)=\{x\}$.  If $f$ is surjective, write $f=\trans{f_1&\cdots&f_{\al-2}&\{f_{\al-1},x\}}{g_1&\cdots&g_{\al-2}&g_{\al-1}}$.

%\pf We may write $f=\trans{f_1&\cdots&f_{\al-1}}{g_1&\cdots&g_{\al-1}}$ or $f=\trans{f_1&\cdots&f_{\al-2}&\{f_{\al-1},x\}}{g_1&\cdots&g_{\al-2}&g_{\al-1}}$, according to whether $f$ is injective or full, respectively.  In the former case, write $X\sm\dom(f)=\{x\}$.  

\pf (i).  Write $f=\trans{f_1&\cdots&f_{\al-1}}{g_1&\cdots&g_{\al-1}}$, $X\sm\dom(f)=\{x\}$, and fix some $y\in Y\sm\im(f)$.  Then $f=h_1\star_ag\star_ah_2$, where $h_1=\trans{f_1&\cdots&f_{\al-1}&x}{b_1&\cdots&b_{\al-1}&b_\al}$ and $h_2=\trans{a_1&\cdots&a_{\al-1}&a_\al}{g_1&\cdots&g_{\al-1}&y}$.

\pfitem{ii}  Write $f=\trans{f_1&\cdots&f_{\al-2}&\{f_{\al-1},x\}}{g_1&\cdots&g_{\al-2}&g_{\al-1}}$.  Fix $y\in Y\sm\im(f)$.  Then we define $h_1,h_2$ exactly as in the previous case. \epf

\begin{lemma}\label{lem:X<Y4}
Suppose $1\leq\al=|X|<|Y|\leq\aleph_0$.
\bit
\itemit{i} If $a$ is full, then $\PTXYa=\la D_\al\cup\{g\}\ra_a$, where $g=\trans{a_1&\cdots&a_{\al-1}}{b_1&\cdots&b_{\al-1}}$.
\itemit{ii} If $a$ is injective, then $\PTXYa=\la D_\al\cup\{g\}\ra_a$, where $g=\trans{a_1&\cdots&a_{\al-2}&\{a_{\al-1},a_\al\}}{b_1&\cdots&b_{\al-2}&b_{\al-1}}$.
\end{itemize}
\end{lemma}

\pf We just prove (i), as the proof of (ii) is virtually identical.  Suppose $a$ is full, and note that this implies $a$ is non-injective.  By Corollary \ref{cor:Dmu_Dmu+1}(i), it suffices to show that $D_{\al-1}\sub\la D_\al\cup\{g\}\ra_a$, so suppose $f\in D_{\al-1}$.  If $f$ is full, then Lemma \ref{lem:X<Y1}(i) gives $f\in \la D_\al\ra_a\sub\la D_\al\cup\{g\}\ra_a$.  If $f$ is non-full, then it is injective (since $\rank(f)=\al-1$), and so Lemma \ref{lem:X<Y3}(i) gives $f\in\la D_\al\cup\{g\}\ra_a$. \epf

%
%
%%\pf (i).  Write $f=\trans{f_1&\cdots&f_{\al-2}&\{f_{\al-1},f_\al\}}{f_1'&\cdots&f_{\al-2}'&f_{\al-1}'}$.
%\pf (i).  Write $f=\trans{f_1&\cdots&f_{\al-2}&\{f_{\al-1},f_\al\}}{g_1&\cdots&g_{\al-2}&g_{\al-1}}$.  Now, $g$ has a similar form.  However, it will be convenient to write $g$ in an even more specific form.  Namely, 
%
%  and $g=\trans{g_1&\cdots&g_{\al-2}&\{g_{\al-1},g_\al\}}{g_1'&\cdots&g_{\al-2}'&g_{\al-1}'}$.  

\begin{thm}\label{thm:rankPTXYa3}
Suppose $1\leq\al=\rank(a)=|X|<|Y|\leq\aleph_0$.  Then
\[
\rank(\PTXYa) = \tbinom{|Y|}\al+
\begin{cases}
0 &\text{if $a$ is neither full nor injective}\\[2truemm]
1 &\text{if $a$ is full or $[a$ is injective and $\al\leq2]$}\\[2truemm]
2 &\text{if $a$ is injective and $\al\geq3$.}
\end{cases}
%\begin{cases}
%\binom{|Y|}\al &\text{if $a$ is neither full nor injective}\\[2truemm]
%\binom{|Y|}\al+1 &\text{if $a$ is full or $[a$ is injective and $\al\leq2]$}\\[2truemm]
%\binom{|Y|}\al+2 &\text{if $a$ is injective and $\al\geq3$.}
%\end{cases}
\]
\end{thm}

\pf First suppose $\PTXYa=\la\Om\ra_a$.  Lemma \ref{lem:rank_PTXYa_lower}(iii) gives $|\Om\cap D_\al| \geq \binom{|Y|}\al$, while Lemma \ref{lem:rank_PTXYa_lower}(iv) gives $|\Om\cap D_\al| \geq \binom{|Y|}\al+1$ if $a$ is injective and $\al\geq3$.  Lemma \ref{lem:X<Y2} gives $|\Om\sm D_\al|\geq1$ if $a$ is injective or full.  Together, these show that the claimed value of $\rank(\PTXYa)$ is a lower bound.  To complete the proof, it suffices to provide a generating set of the specified size.

By Theorem \ref{thm:inflation_PT}(v) and Lemma \ref{lem:maximal_J_PT}, $D_\al^a=\Hh_b^a$ is an $(\al+1)^\be\times\Lam_I$ rectangular group over $\S_\al$.  Since $\al=|X|$, $a$ is surjective, so $\be=0$ and $(\al+1)^\be=1$.  If $a$ is not injective, or if $\al\leq2$, then $\Lam_I\geq\rank(\S_\al)$, so \cite[Proposition 3.11(i)]{Sandwiches1} gives
\[
\rank(D_\al^a)=\max(1,\Lam_I,\rank(\S_\al)) = \begin{cases}
\Lam_I &\text{if $a$ is non-injective or $\al\leq2$}\\
2 &\text{if $a$ is injective and $\al\geq3$.}
\end{cases}
\]
Choose a minimal-size generating set $\Om$ for $D_\al=\Hh_b^a$.  In the first case (in which $a$ is non-injective or $\al\leq2$), by \cite[Proposition 3.11(iv)]{Sandwiches1}, we may assume that $\Om$ is a cross-section of the $\L^a$-classes of $D_\al^a$.  In the second case, as in the proof of Lemma \ref{lem:rank_PTXYa_lower}(iv),~$D_\al^a=\Hh_b^a$ is a single $\H^a$-class, which is a group isomorphic to $\S_\al$, and $\Om$ has two elements from this group.  Let $\Ga$ be a cross-section of the non-regular $\D^a$-classes contained in~$D_\al$.  By Lemma \ref{lem:OmGa}, $D_\al\sub\la D_\al^a\cup\Ga\ra_a=\la\Om\cup\Ga\ra_a$.  As in the proof of Theorem~\ref{thm:rankPTXYa2}, noting that ${\L^a}={\L}$ on $\PTXYa$ and that $|D_\al/{\L}|=\binom{|Y|}\al$, by Lemma \ref{lem:al=xi}(ii) and Corollary \ref{cor:Green_sizes_PT}, respectively, we have 
\[
|\Om\cup\Ga| = \tbinom{|Y|}\al+
\begin{cases}
0 &\text{if $a$ is non-injective or $\al\leq2$}\\
1 &\text{if $a$ is injective and $\al\geq3$.}
\end{cases}
\]
If $a$ is neither full nor injective, then Lemma \ref{lem:X<Y1}(iii) gives $\PTXYa=\la D_\al\ra_a=\la\Om\cup\Ga\ra_a$, and the proof is complete in this case.  If $a$ is full or injective, then Lemma \ref{lem:X<Y4} says that $\PTXYa=\la D_\al\cup\{g\}\ra_a=\la\Om\cup\Ga\cup\{g\}\ra_a$ for some $g\in D_{\al-1}$, completing the proof in these cases also. \epf

\begin{rem}
Since $\PTXYa\cong\PT(Y,\si)$ if $a$ is surjective, Theorem \ref{thm:rankPTXYa3} gives a formula for $\rank(\PT(Y,\si))$ in the case that $1\leq|Y|\leq\aleph_0$ and $|Y/\si|<\aleph_0$; for $|Y|=0$ or $|Y|>\aleph_0$ and/or $|Y/\si|\geq\aleph_0$, see the discussion at the beginning of this section.  To the authors' knowledge, $\rank(\PT(Y,\si))$ has not been previously calculated.
In the case that $\si$ is the diagonal relation, $\PT(Y,\si)=\set{f\in\PT_Y}{\dom(f)\sub \A}$, so Theorem \ref{thm:rankPTXYa3} gives the rank of these semigroups as well; this rank is $\binom{|Y|}{|\A|}+1$ or $\binom{|Y|}{|\A|}+2$, according to whether $|\A|\leq2$ or $|\A|\geq3$, respectively.
\end{rem}

\begin{rem}
As in Remark \ref{rem:rankPTXYa1}, again note that $\rank(\PTXYa)=|\PTXYa|=\aleph_0$ if $1\leq\al=|X|<|Y|=\aleph_0$.
Also note that the value of $\rank(\PTXYa)$ given in Theorem \ref{thm:rankPTXYa3} does not depend on the actual sizes, $\lam_1,\ldots,\lam_\al$, of the $\ker(a)$-classes, apart from some small effect in the case that $a$ is injective (i.e., when $\lam_i=1$ for all $i$).
%
%Finally, note that if $1\leq\al=|X|<|Y|=\aleph_0$, then Theorem \ref{thm:rankPTXYa3} gives $\rank(\PTXYa)=|\PTXYa|=\aleph_0$.  
\end{rem}

\subsectiontitle{Egg-box diagrams}\label{sect:eggbox}

In this section, we give computer-generated egg-box diagrams for a number of sandwich semigroups $\PTXYa$, in large part to illustrate various results in the preceding sections; these diagrams were all produced with the {\sc Semigroups} package in GAP \cite{GAP}, and we thank Dr Attila Egri-Nagy and Dr James Mitchell for writing the code for creating them.  In each such diagram, elements related under the relevant Green's~$\R$ (or $\L$ or $\H$) relation are in the same row (or column or cell, respectively).  The $\leq_{\J}$ order is indicated in the usual way (with a line segment from one $\J$-class up to another $\J$-class if the former is covered by the latter).  Group $\H$-classes are always coloured gray.  In all cases, we assume $X=\{1,\ldots,m\}$ and~$Y=\{1,\ldots,n\}$ for appropriate $m,n\in\N$, writing $\PT_{mn}$ for $\PT_{\{1,\ldots,m\},\{1,\ldots,n\}}$, and we denote a partial transformation in standard tableau form: for example, $f=\trans{1&2&3&4&5}{3&4&3&-&1}\in\PT_{54}$ denotes the partial transformation with ${\dom(f)=\{1,2,3,5\}}$ and $\im(f)=\{1,3,4\}$, under which $1\mt3$, $2\mt4$, $3\mt3$ and $5\mt1$.

\begin{figure}[ht]
\begin{center}
\includegraphics[height=4cm]{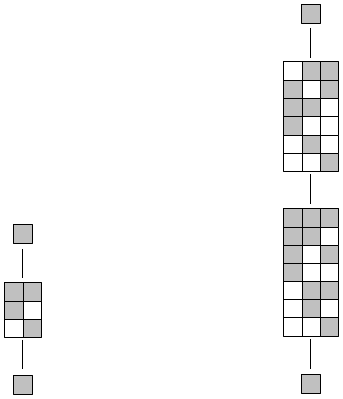} 
\end{center}
\vspace{-5mm}
\caption{Egg-box diagrams of the partial transformation semigroups $\PT_2$ (left) and $\PT_3$ (right).}
\label{fig:PT_1}
\end{figure}

\begin{figure}[ht]
\begin{center}
\includegraphics[width=\textwidth]{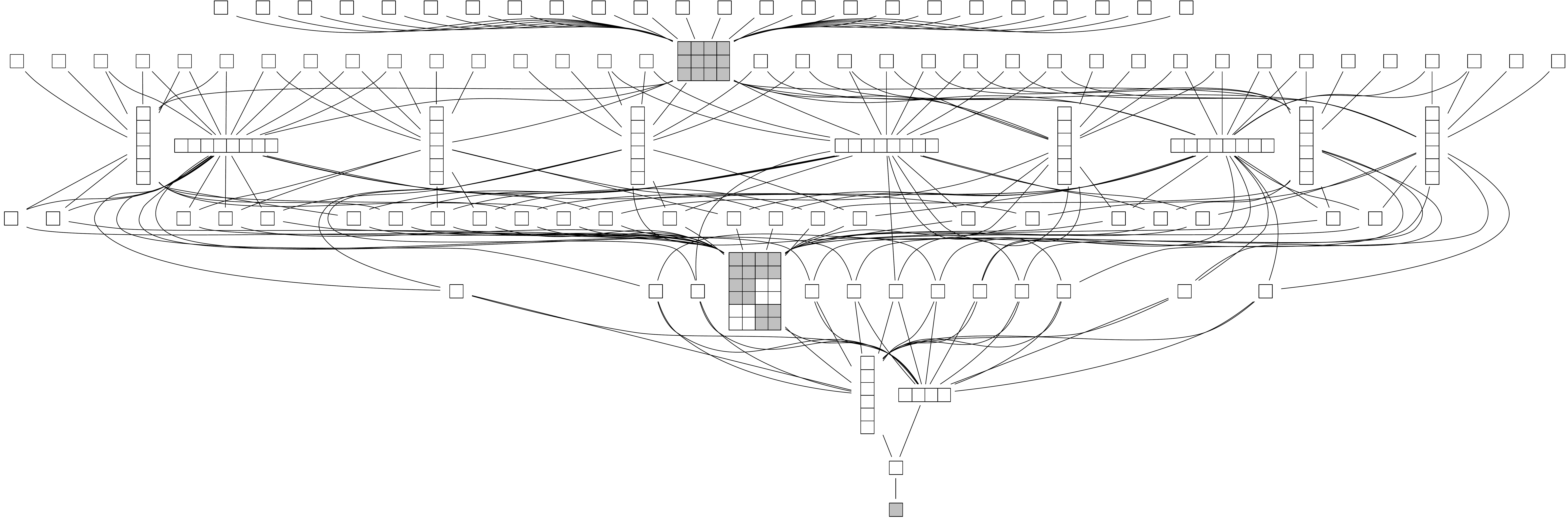} \\~\\
\vspace{-5mm}
\caption[blah]{Egg-box diagram of the sandwich semigroup $\PT_{35}^a$, where $a=\trans{1&2&3&4&5}{1&1&2&2&-}\in\PT_{53}$.  Note that $a$ is neither full, nor injective, nor surjective.}
\label{fig:PT_2}
\end{center}
\end{figure}

\begin{figure}[ht]
\begin{center}
\includegraphics[width=\textwidth]{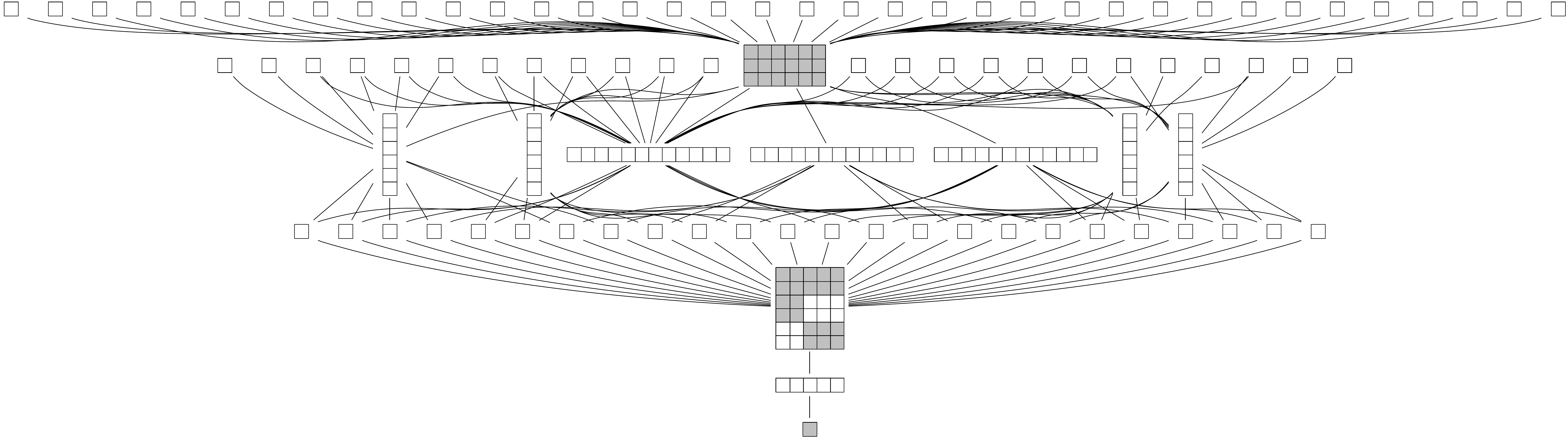} \\~\\
\vspace{-5mm}
\caption[blah]{Egg-box diagram of the sandwich semigroup $\PT_{35}^b$, where $b=\trans{1&2&3&4&5}{1&1&2&2&2}\in\PT_{53}$. Note that $b$ is full, but neither injective nor surjective.}
\label{fig:PT_3}
\end{center}
\end{figure}

\begin{figure}[ht]
\begin{center}
\includegraphics[width=\textwidth]{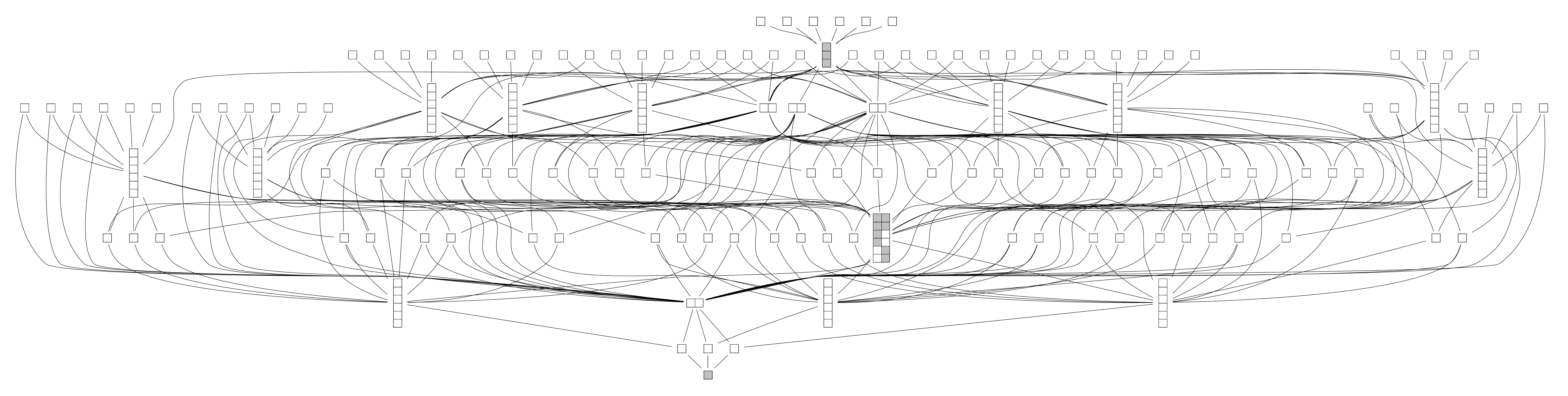} \\~\\
\vspace{-5mm}
\caption[blah]{Egg-box diagram of the sandwich semigroup $\PT_{35}^c$, where $c=\trans{1&2&3&4&5}{1&2&-&-&-}\in\PT_{53}$. Note that $c$ is injective, but neither full nor surjective.}
\label{fig:PT_4}
\end{center}
\end{figure}

\begin{figure}[ht]
\begin{center}
\includegraphics[width=0.5\textwidth]{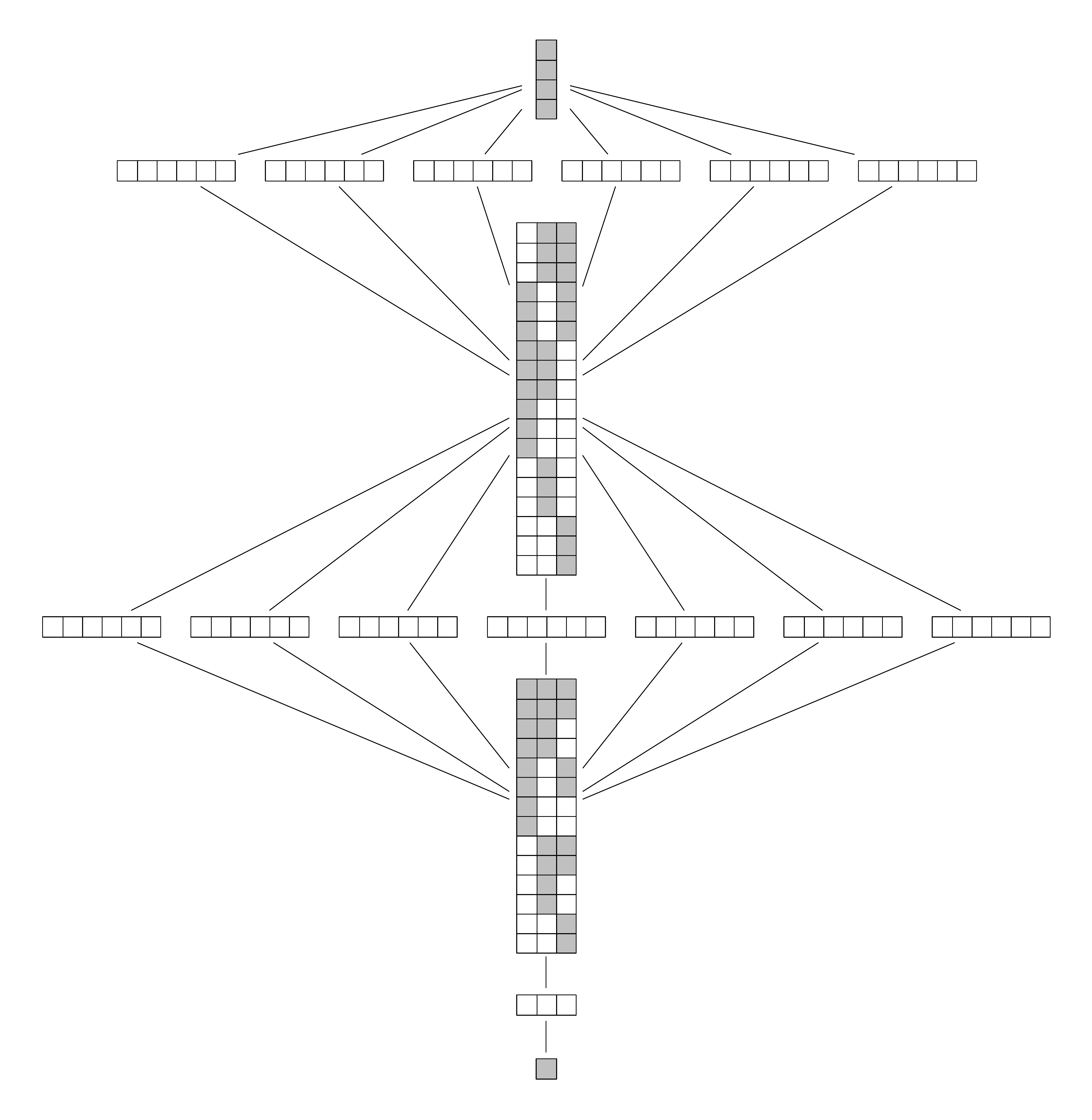} \\~\\
\vspace{-5mm}
\caption[blah]{Egg-box diagram of the sandwich semigroup $\PT_{43}^d$, where $d=\trans{1&2&3}{1&2&3}\in\PT_{34}$. Note that $d$ is full and injective, but not surjective.}
\label{fig:PT_5}
\end{center}
\end{figure}

\begin{figure}[ht]
\begin{center}
\includegraphics[height=10cm]{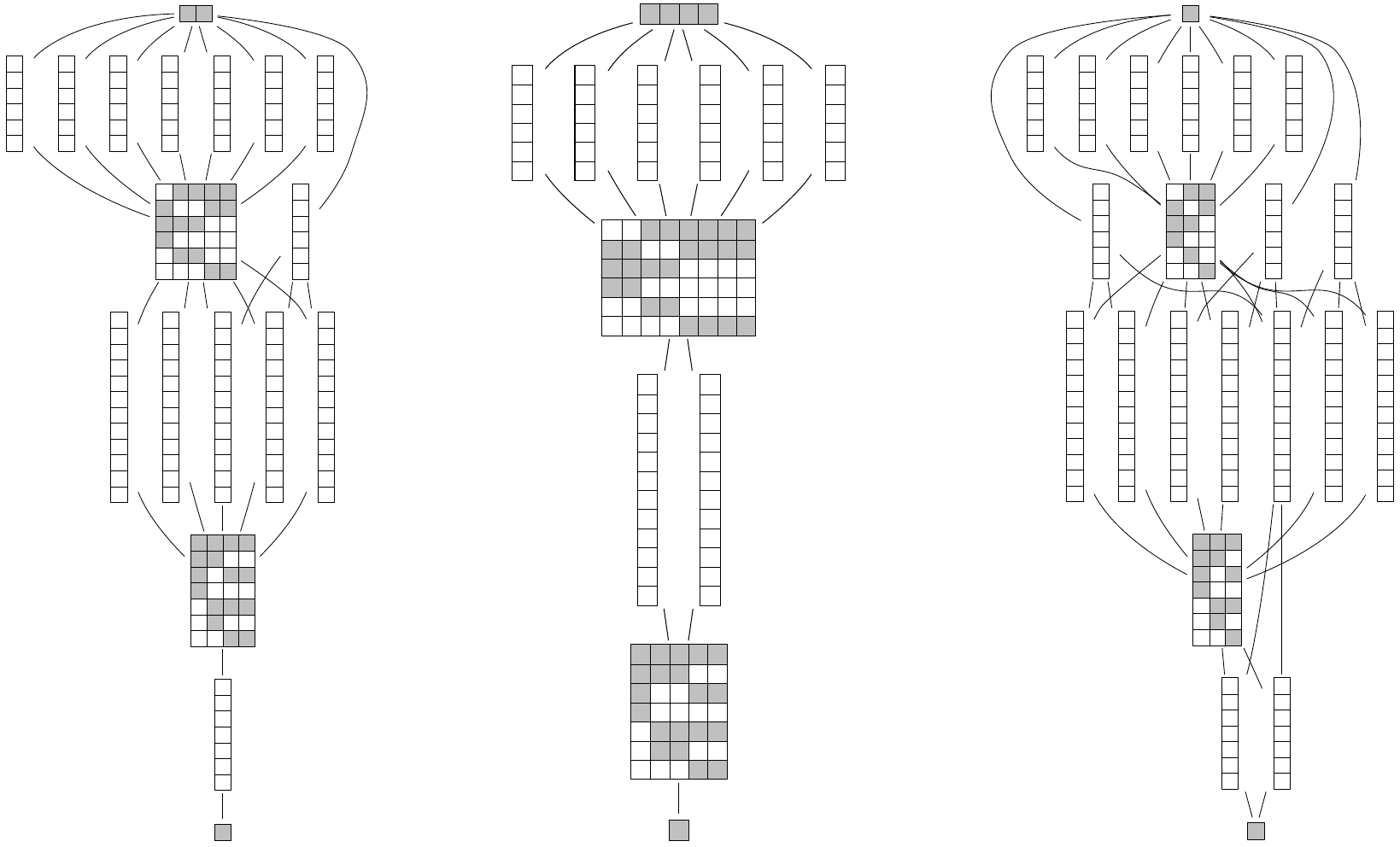}
\caption[blah]{Left to right: egg-box diagrams of the sandwich semigroups $\PT_{35}^e$, $\PT_{35}^f$ and $\PT_{35}^g$, where $e=\trans{1&2&3&4&5}{1&2&3&3&-},\ f=\trans{1&2&3&4&5}{1&1&2&2&3},\ g=\trans{1&2&3&4&5}{1&2&3&-&-}\in\PT_{53}$. Note that $e$ is surjective, but neither full nor injective; $f$ is full and surjective, but not injective; and $g$ is injective and surjective, but not full.}
\label{fig:PT_6}
\end{center}
\end{figure}

\begin{figure}[ht]
\begin{center}
\includegraphics[height=8cm]{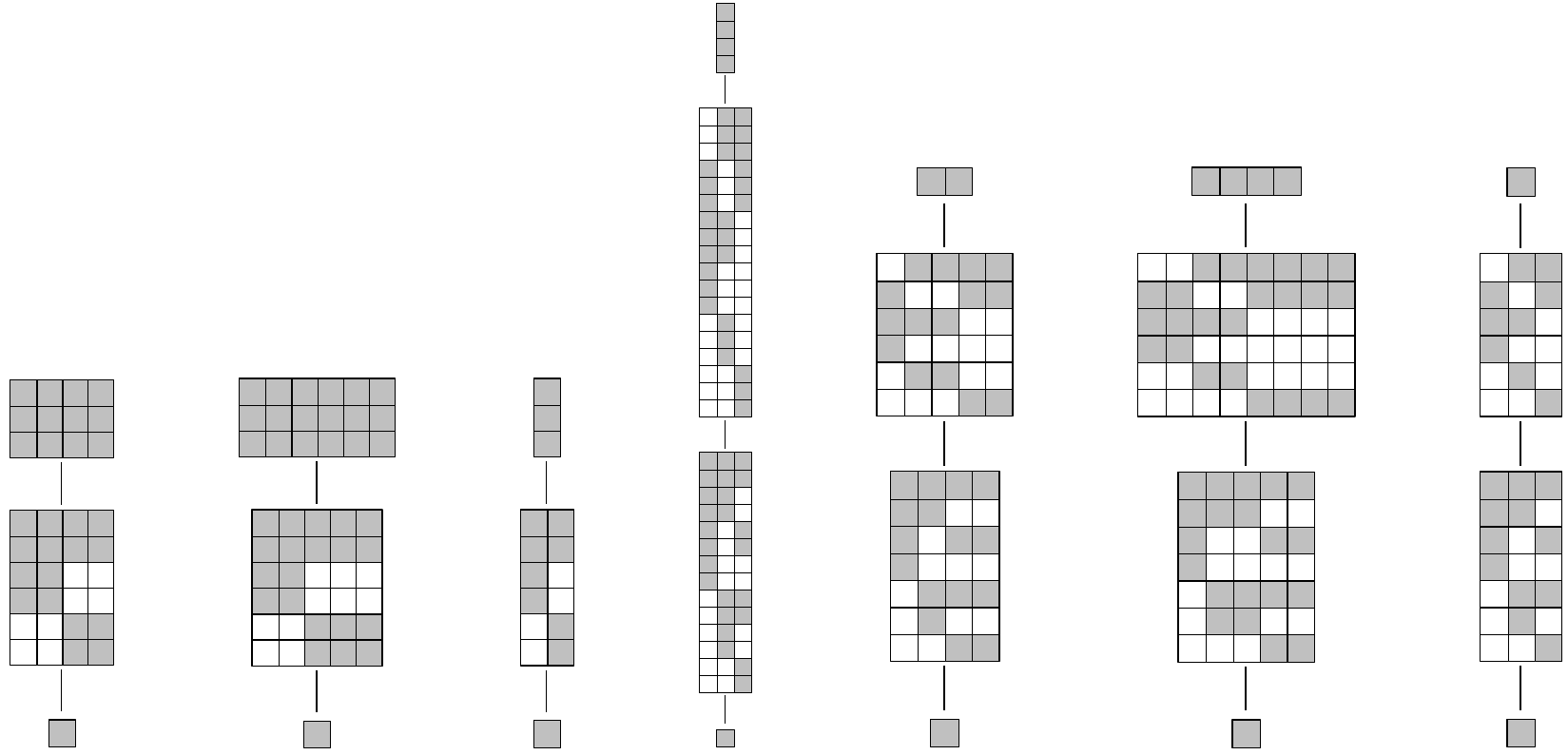}
\caption[blah]{Left to right: egg-box diagrams of the regular sandwich semigroups $\Reg(\PT_{35}^a)$, $\Reg(\PT_{35}^b)$, $\Reg(\PT_{35}^c)$, $\Reg(\PT_{43}^d)$, $\Reg(\PT_{35}^e)$, $\Reg(\PT_{35}^f)$ and $\Reg(\PT_{35}^g)$, where the sandwich elements $a,b,c,d,e,f,g$, are as in Figures \ref{fig:PT_2}--\ref{fig:PT_6}.  The first three are inflations (in the sense of \cite[Remark 2.15]{Sandwiches1}) of $\PT_2$, and the other four are inflations of $\PT_3$; cf Figure \ref{fig:PT_1}.  In fact, $\Reg(\PT_{35}^g)\cong\PT_3$, as in Remark \ref{rem:PTYsiPTA}.}
\label{fig:PT_9}
\end{center}
\end{figure}

\sectiontitle{The category $\T$}\label{sect:T}

Recall that the category $\T = \bigset{(A,f,B)}{A,B\in\Setp,\ f\in\bT_{AB}}$ is the (regular monoidal) partial subsemigroup of the category $\PT$ consisting of all \emph{full} transformations between \emph{non-empty} sets.  Our thorough investigation of the category $\PT$ in Section \ref{sect:PT} may be adapted to prove analogous results concerning the category $\T$.  Rather than giving the full details, we will just list the main results, and comment on similarities with $\PT$.  Many results in this section are either proved along similar lines to the corresponding results of Section \ref{sect:PT}, or else may be obtained as direct consequences of results in Section \ref{sect:PT}, by applying a suitable ``reduction'' result from \cite[Section 1.4]{Sandwiches1}.  We also note that many of the results of \cite{DEvariants} are special cases of the results in this section, taking $|X|=|Y|<\aleph_0$.

\subsectiontitle{Green's relations and stability in $\T$}\label{sect:Green_T}

As in Proposition \ref{prop:GreenPT}, we may describe Green's relations and preorders on $\T$.  If $(A,f,B),(C,g,D)\in\T$, then~
\bit
\itemnit{i} $(A,f,B)\leqR(C,g,D) \iff A=C$ and $\ker(f)\supseteq\ker(g)$,
\itemnit{ii} $(A,f,B)\leqL(C,g,D) \iff B=D$ and $\im(f)\sub\im(g)$,
\itemnit{iii} $(A,f,B)\leqJ(C,g,D) \iff \rank(f)\leq\rank(g)$,
\itemnit{iv} $(A,f,B)\R(C,g,D) \iff A=C$ and $\ker(f)=\ker(g)$,
\itemnit{v} $(A,f,B)\L(C,g,D) \iff B=D$ and $\im(f)=\im(g)$,
\itemnit{vi} $(A,f,B)\J(C,g,D) \iff (A,f,B)\D(C,g,D) \iff \rank(f)=\rank(g)$.
\eit
Indeed, parts (i), (ii), (iv) and (v) follow directly from Proposition \ref{prop:GreenPT} and \cite[Lemma 1.8]{Sandwiches1}; note that statements concerning domains in Proposition \ref{prop:GreenPT} become vacuous when considering $\T$, as all transformations are full.  Parts (iii) and (vi) are proved in similar fashion to the corresponding parts of Proposition \ref{prop:GreenPT}; we only note that in the proof of (iii), the partial transformation $h_2\in\bPT_{DB}$ needs to be arbitrarily extended to a full transformation from $\bT_{DB}$. %\epf

As in Corollary \ref{cor:Jclasses_PT}, but noting that the minimum rank of a full transformation is $1$, it follows that the~${\J}={\D}$-classes of $\T_{AB}$ are the sets
\[
D_\mu = D_\mu(\T_{AB}) = \set{(A,f,B)}{f\in\bT_{AB},\ \rank(f)=\mu} \qquad\text{for each cardinal \ $1\leq\mu\leq\min(|A|,|B|)$.}
\]
%(Note that the minimum rank of a full transformation is $1$.)
These $\J$-classes form a chain: $D_\mu\leqJ D_\nu \iff \mu\leq\nu$.  As in Corollary \ref{cor:Green_sizes_PT}, if $|A|=\al$ and $|B|=\be$, then
\[
|D_\mu/{\R}| = S(\al,\mu) \COMMA
|D_\mu/{\L}| = \tbinom\be\mu \COMMA
|D_\mu/{\H}| = \tbinom\be\mu S(\al,\mu) \COMMA
|D_\mu| = \mu! \tbinom\be\mu S(\al,\mu) ,
\]
and each $\H$-class in $D_\mu$ has size $\mu!$.  Consequently,
\[
\be^\al=|\T_{AB}|=\sum_{\mu=1}^{\min(\al,\be)}\mu! \tbinom\be\mu S(\al,\mu).
\]
As in Lemma \ref{lem:stable_PT} (and using a preliminary result analogous to Lemma \ref{lem:PTfr}), we may also prove that if $(A,f,B)\in\T$, then
\bit
\itemnit{i} $(A,f,B)$ is $\R$-stable $\iff$ $[\rank(f)<\aleph_0$ or $f$ is injective$]$,
\itemnit{ii} $(A,f,B)$ is $\L$-stable $\iff$ $[\rank(f)<\aleph_0$ or $f$ is surjective$]$,
\itemnit{iii} $(A,f,B)$ is stable $\iff$ $[\rank(f)<\aleph_0$ or $f$ is bijective$]$.
\eit
One could also use \cite[Lemma 1.9]{Sandwiches1} to prove the backwards implications in each of (i)--(iii).

\subsectiontitle{Green's relations, regularity and stability in $\TXYa$}\label{sect:TXYa}

In order to study sandwich semigroups in $\T$, for the rest of Section \ref{sect:T}, we fix two non-empty sets $X,Y\in\Setp$.
Again, we will identify $\bT_{ZW}$ with $\T_{ZW}$, where $Z,W\in\{X,Y\}$, via $f\equiv(Z,f,W)$.
For the rest of Section \ref{sect:T}, we also fix a transformation $a\in\T_{YX}$, and write
\begin{equation}\label{eq:aT}
a=\tbinom{A_i}{a_i}_{i\in I} \COMMA
A=\im(a)=\set{a_i}{i\in I} \COMMA
\si=\ker(a) \COMMA
\al=\rank(a),
\end{equation}
noting that $\dom(a)=\bigcup_{i\in I}A_i=Y$, $\si$ is an equivalence on $Y$, $Y/\si=\set{A_i}{i\in I}$, and $\al=|I|=|A|=|Y/\si|$.  
For each $i\in I$, we fix some $b_i\in A_i$.  We also fix a partition $\set{B_i}{i\in I}$ of $X$ such that $a_i\in B_i$ for each $i$, and we define
\begin{equation}\label{eq:bT}
b=\tbinom{B_i}{b_i}\in\T_{XY},
\end{equation}
so that $a=aba$ and $b=bab$.  We also define
\begin{equation}\label{eq:cT}
\be=|X\sm\im(a)| \COMMA \xi=\min(|X|,|Y|) \COMMA \lam_i=|A_i| \text{ for } i\in I \COMMA \Lam_J=\prod_{j\in J}\lam_j \text{ for } J\sub I,
\end{equation}
noting that $\sum_{i\in I}\lam_i=|Y|$.

The sets 
\[
P_1^a = \set{f\in\TXY}{fa\R f} \COMma
P_2^a = \set{f\in\TXY}{af\L f} \COMma
P_3^a = \set{f\in\TXY}{afa\J f} \COMma
P^a=P_1^a\cap P_2^a
\]
may also be easily described:
\bit
\itemnit{i} $P_1^a = \set{f\in\TXY}{\ker(fa)=\ker(f)} = \set{f\in\TXY}{\ker(a)\text{ separates }\im(f)}$,
\itemnit{ii} $P_2^a = \set{f\in\TXY}{\im(af)=\im(f)} = \set{f\in\TXY}{\im(a)\text{ saturates }\ker(f)}$,
\itemnit{iii} $P^a = \set{f\in\TXY}{\ker(fa)=\ker(f),\ \im(af)=\im(f)}$
\item[] $\phantom{P^a} = \set{f\in\TXY}{\ker(a)\text{ separates }\im(f),\ \im(a)\text{ saturates }\ker(f)}$,
\itemnit{iv} $P_3^a = \set{f\in\TXY}{\rank(afa)=\rank(f)}$. 
\eitres
From these, a characterisation of Green's relations on $\TXYa$ may be deduced; we won't state this explicitly, as it is exactly analogous to Theorem \ref{thm:green_PTXYa} (see also \cite[Remark 2.3]{MGS2013}).  Since $\T$ is regular, we again have $\Reg(\TXYa)=P^a$; this last statement may also be deduced from \cite[Theorem 5.3]{MS1975}.  Also, as in Propositions~\ref{prop:JaDaPT} and \ref{prop:Reg(PYXYa)}, one may show that
\begin{itemize}\begin{multicols}{3}
\itemnit{i} $a$ is $\R$-stable $\iff$ $P_3^a\sub P_1^a$,
\item[]
\itemnit{ii} $a$ is $\L$-stable $\iff$ $P_3^a\sub P_2^a$,
\item[]
\itemnit{iii} $a$ is stable $\iff$ $P_3^a=P^a$
\item[] {\white $a$ is stable} $\iff$ ${\J^a}={\D^a}$.
\end{multicols}\eitmc
As in Propositions \ref{prop:Ja_order_PT} and \ref{prop:Ja_order_PT_P}, we may characterise the ordering on $\J^a$-classes of $\TXYa$.  Since these are directly analogous to the just-mentioned results, we will not state these here.  Again the reader should note that simplifications occur; for example, condition (iv) of Proposition \ref{prop:Ja_order_PT} reduces to ``$\ker(f)\supseteq\ker(ga)$''.  However, it is worth noting that (as in Proposition \ref{prop:regular_Da_classes}) the regular $\gDa$-classes of $\TXYa$ are precisely the sets
\[
D_\mu^a = \set{f\in P^a}{\rank(f)=\mu} \qquad\text{for each cardinal $1\leq\mu\leq\al=\rank(a)$,}
\]
and that, if $f\in P^a$, then $D_f^a=J_f^a$ if and only if $\rank(f)<\aleph_0$ or $a$ is stable.  Also, as in Proposition \ref{prop:maximal_J_PT}:
\bit
\itemnit{i} If $\al<\xi$, then the maximal $\J^a$-classes of $\TXYa$ are the singleton sets~$\{f\}$, for $f\in\TXY$ with $\rank(f)>\al$.
\itemnit{ii} If $\al=\xi$, then the set $J_b^a=\set{f\in P_3^a}{\rank(f)=\al}$ is a maximum $\J^a$-class of $\TXYa$.
\eit
A subtle but crucial difference between sandwich semigroups in $\PT$ and $\T$ is that the minimum $\J^a$-class in $\TXYa$ is not a singleton in general, since $\T$ has no mappings of rank $0$.  Rather, this minimum $\J^a$-class is the set of all constant mappings $X\to Y$,
\[
D_1^a = \set{f\in P^a}{\rank(f)=1} = \set{f\in\T_{XY}}{\rank(f)=1},
\]
which has size $|Y|$.

\subsectiontitle{A structure theorem for $\Reg(\TXYa)$ and connections to (non-sandwich) semigroups of transformations}\label{sect:structureT}

As in Section \ref{sect:structurePT}, the commutative diagrams \cite[diagrams (2.1) and~(2.7)]{Sandwiches1} lead to relationships between the sandwich semigroup $\TXYa$ and certain (non-sandwich) transformation semigroups: namely,
\[
\TXY a = \T_X(ba) = \set{f\in\T_X}{\im(f)\sub A} 
\AND
a\TXY = (ab)\T_Y = \set{f\in\T_Y}{\ker(f)\supseteq\si}.
\]
These are subsemigroups (indeed, principal one-sided ideals) of $\T_X$ and $\T_Y$, respectively; they have been studied extensively in the literature, where they are commonly denoted by $\T(X,A)$ and $\T(Y,\si)$, respectively; see for example \cite{SS2013,FS2014,MGS2010}.  We also note that the regular monoid $(a\TXY a,\starb)$, which is isomorphic to~${ba\TXY a=(ba)\T_X(ba)}$ by \cite[Remark 2.6]{Sandwiches1}, is in fact isomorphic to $\T_A$.  Indeed, this is (only slightly) less obvious than the corresponding isomorphism $(a\PTXY a,\starb)\cong\PT_A$ discussed in Section \ref{sect:structurePT}.  Note that~$ba=\binom{B_i}{a_i}\in\T_X$, and since $a_i\in B_i$ for each $i$, every element of $(ba)\T_X(ba)$ is uniquely determined by its restriction to $A$.  It quickly follows that $(ba)\T_X(ba)\to\T_A:f\mt f|_A$ is an isomorphism.  Thus, \cite[diagrams~(2.1) and~(2.7)]{Sandwiches1} yield commutative diagrams of semigroup epimorphisms:
\begin{equation}\label{eq:CD_T_2}
\includegraphics{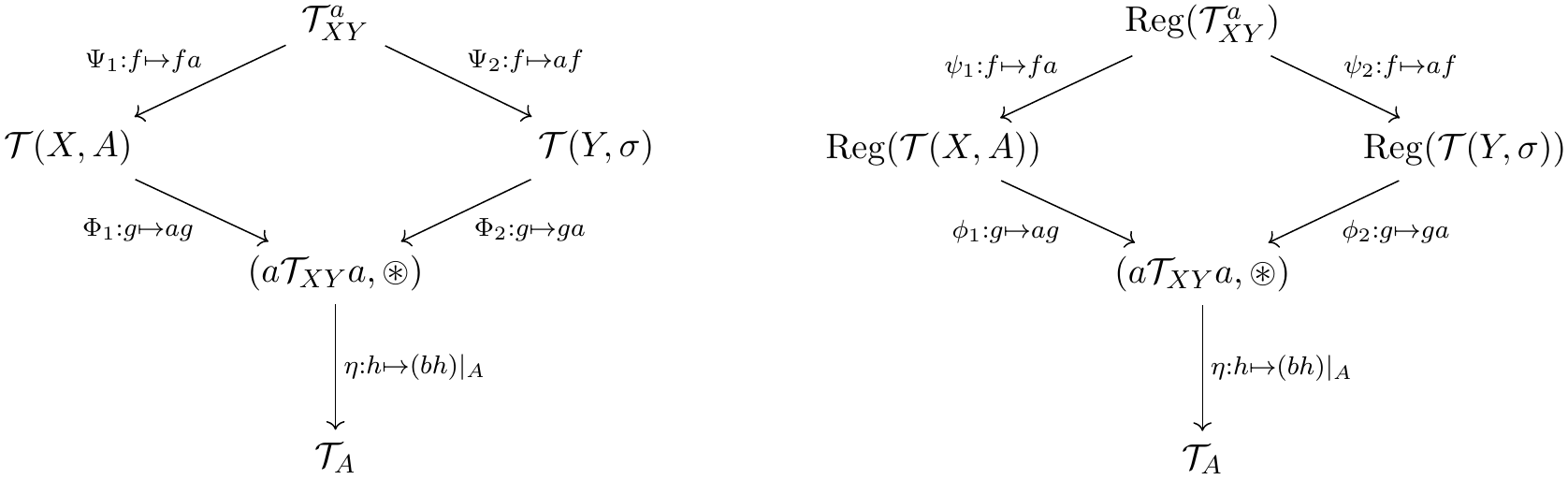}
%\begin{tikzcd} %[column sep=small]
%~ & \TXYa \arrow[swap]{dl}{\Psi_1:f\mt fa} \arrow{dr}{\Psi_2:f\mt af} & \\
%\T(X,\B) \arrow[swap]{dr}{\Phi_1:g\mt ag} & & \T(Y,\si) \arrow{dl}{\Phi_2:g\mt ga}\\
%& (a\TXY a,\starb)\arrow{dd}{\eta:h\mt (bh)|_A} & \\
%\\
%& \T_\B & 
%\end{tikzcd}
%\qquad \qquad
%\begin{tikzcd} [column sep=small]
%~ & \Reg(\TXYa) \arrow[swap]{dl}{\psi_1:f\mt fa} \arrow{dr}{\psi_2:f\mt af} & \\
%\Reg(\T(X,\B)) \arrow[swap]{dr}{\phi_1:g\mt ag} & & \Reg(\T(Y,\si)) \arrow{dl}{\phi_2:g\mt ga}\\
%& (a\TXY a,\starb)\arrow{dd}{\eta:h\mt (bh)|_A}  & \\
%\\
%& \T_\B & 
%\end{tikzcd}
\end{equation}
If $a$ is injective, then $\Psi_1$ and $\psi_1$ are isomorphisms.  If $a$ is surjective, then $\Psi_2$ and $\psi_2$ are isomorphisms.  If $a$ is bijective, then all maps in \eqref{eq:CD_T_2} are isomorphisms.  Consequently, $\T(X,A)$ and $\T(Y,\si)$ arise as special cases of the $\TXYa$ construction.  As with Theorem \ref{thm:psi_T}, the next result follows immediately from \cite[Theorem~2.10]{Sandwiches1}, and the above observations.

\begin{thm}\label{thm:psi_T}
Consider the map $\psi:\Reg(\TXYa) \to \Reg(\T(X,\B))\times\Reg(\T(Y,\si)):f\mt(fa,af)$.  Then $\psi$ is injective, and $\im(\psi) = \set{(g,h)\in \Reg(\T(X,\B))\times\Reg(\T(Y,\si))}{ag=ha}$.  In particular, $\Reg(\TXYa)$ is a pullback product of $\Reg(\T(X,\B))$ and $\Reg(\T(Y,\si))$ with respect to~$\PT_\B$.~\epfres
\end{thm}

It is also possible to quickly deduce simple descriptions of the regular subsemigroups $\Reg(\T(X,\B))$ and $\Reg(\T(Y,\si))$; we will not state these description explicitly, as the reader may easily modify Proposition \ref{prop:RegPTXBYsi} as appropriate, and since they have been previously given in \cite{MGS2010,SS2008}.

\subsectiontitle{The regular subsemigroup $P^a=\Reg(\TXYa)$}\label{sect:RegTXYa}

As in Section \ref{sect:RegPTXYa}, we may give a thorough structural description of the regular subsemigroup $P^a=\Reg(\TXYa)$.  First, we note that Green's $\R$, $\L$, $\H$ and $\D$ relations on $P^a$ are simply the restrictions of the corresponding relations $\R^a$, $\L^a$, $\H^a$ and $\D^a$ on $\TXYa$, and that ${\J^{P^a}}={\D^a}$.  As in Proposition \ref{prop:green_Pa_PT}, if $f\in P^a$, then $K_f^a=K_f\cap P^a$ if $\K$ is any of $\R,\L,\H,\D$.  The ${\J^{P^a}}={\D^a}$-classes of $P^a$ are the sets
\[
D_\mu^a = \set{g\in P^a}{\rank(g)=\mu} \qquad\text{for each cardinal $1\leq\mu\leq\al=\rank(a)$,}
\]
and these form a chain under the $\J^{P^a}$-ordering: $D_\mu^a\leq D_\nu^a \iff \mu\leq\nu$.

The internal structure of the ${\D^a}$-classes of $P^a=\Reg(\TXYa)$ may be described using the $\gKh^a$-relations; these are defined by means of the epimorphism
\[
\varphi=\psi_1\phi_1\eta=\psi_2\phi_2\eta:P^a\to\T_A:f\mt(bafa)|_A=(fa)|_A.
\]
As in Theorem \ref{thm:inflation_PT}, we have the following.

\begin{thm}\label{thm:inflation_T}
Let $f=\tbinom{F_j}{f_j}_{j\in J}\in P^a=\Reg(\TXYa)$, where $J\sub I$ and $f_j\in A_j$ for each $j$, and write $\mu=\rank(f)$.~  
\bit
\itemit{i} $\Rh_f^a$ is the union of $\mu^\be$\ $\gRa$-classes of $P^a$.
\itemit{ii} $\Lh_f^a$ is the union of $\Lam_J$ $\gLa$-classes of $P^a$.
\itemit{iii} $\Hh_f^a$ is the union of $\mu^\be\Lam_J$ $\gHa$-classes of $P^a$, each of which has size $\mu!$.  
\itemit{iv} If $H_{\fb}$ is a non-group $\gH$-class of $\T_\B$, then each $\gHa$-class of $P^a$ contained in $\Hh_f^a$ is a non-group.
\itemit{v} If $H_{\fb}$ is a group $\gH$-class of $\T_\B$, then each $\gHa$-class of $P^a$ contained in $\Hh_f^a$ is a group isomorphic to the symmetric group $\S_\mu$.  Further, $\Hh_f^a$ is a $\mu^\be\times\Lam_J$ rectangular group over $\S_\mu$.
\itemit{vi} $\Dh_f^a=D_f^a=D_\mu^a=\set{g\in P^a}{\rank(g)=\mu}$ is the union of:
\item[]
\emph{(a)} $\mu^\be S(\al,\mu)$ $\gRa$-classes, \qquad
\emph{(b)} $\displaystyle\sum_{K\sub I \atop |K|=\mu}\Lam_K$ $\gLa$-classes, \qquad
\emph{(c)} $\displaystyle\mu^\be S(\al,\mu)\sum_{K\sub I \atop |K|=\mu}\Lam_K$ $\gHa$-classes. \epfres
\eit
\end{thm}

As in Corollary \ref{cor:size_P_PT}, it follows that
\[
|P^a| = \sum_{\mu=1}^\al \mu!\mu^\be S(\al,\mu)\sum_{K\sub I \atop |K|=\mu}\Lam_K.
\]
If $\al=1$, then the above expression reduces to $|P^a|=\sum_{i\in I}\lam_i=|Y|$, agreeing with the previously-mentioned fact that $P^a=D_1^a=D_1(\TXY)$ if $\al=1$.  
As in Proposition \ref{prop:size_P_PT}, there is also some simplification possible in the case that $|P^a|$ is infinite (note some differences because mappings in $\T$ have minimum rank of $1$):
\bit
\itemnit{i} If $\al\geq2$ and $|X|\geq\aleph_0$, then $|P^a|=2^{|X|}\Lam_I=\max(2^{|X|},\Lam_I)$.
\itemnit{ii} If $\al\geq2$, $|X|<\aleph_0$, and $\lam_i\geq\aleph_0$ for some $i\in I$, then $|P^a|=\Lam_I=\displaystyle\max_{i\in I}\lam_i$.
\itemnit{iii} $|P^a|<\aleph_0 \iff [\al=1$ and $|Y|<\aleph_0]$ or $[\al\geq2$, $|X|<\aleph_0$ and $\lam_i<\aleph_0$ for all $i\in I]$.
\itemnit{iv} $|P^a|=\aleph_0 \iff [\al=1$ and $|Y|=\aleph_0]$ or $[\al\geq2$, $|X|<\aleph_0$ and $\max_{i\in I}\lam_i=\aleph_0]$.
\itemnit{v} $|P^a|>\aleph_0 \iff [\al=1$ and $|Y|>\aleph_0]$ or $[\al\geq2$ and $[|X|\geq\aleph_0$ or $\lam_i>\aleph_0$ for some $i\in I]]$.
\eit
The proof of Proposition \ref{prop:LMC_PT} works essentially unchanged (quoting Theorem 3.2 of \cite{Tirasupa1979}, instead of Theorem~3.1) to prove that $P^a=\Reg(\TXYa)$ is MI-dominated for any $\al$, while $P^a$ is RP-dominated if and only if $\al<\aleph_0$.  As with Theorem \ref{thm:rank_P_PT}, the MI-domination property may then be used to prove the following result, concerning the rank of $P^a=\Reg(\TXYa)$:

\begin{thm}\label{thm:rank_P_T}
Suppose $a$ is not a bijection.
\bit
\itemit{i} If $|P^a|\geq\aleph_0$, then $\rank(P^a)=|P^a|$.
\itemit{ii} If $|P^a|<\aleph_0$ (so that $[\al=1$ and $|Y|<\aleph_0]$ or $[\al\geq2$, $|X|<\aleph_0$ and $\lam_i<\aleph_0$ for all $i\in I]$, as noted above), then
\[
\rank(P^a) = 
\begin{cases}
|Y| &\text{if $\al=1$}\\
1+\max(\al^\be,\Lam_I) &\text{if $\al\geq2$.}
\end{cases}
\]
\end{itemize}
\end{thm}

\pf The proof is essentially identical to that of Theorem \ref{thm:rank_P_PT}.  However, in part (ii), if $\al\geq2$, then we obtain 
\[
\rank(P^a) = \rank(\T_\al:\S_\al) + \max\big( \al^\be, \Lam_I, \rank(\S_\al) \big).
\]
It is well known that $\rank(\T_\al:\S_\al)=1$ for $\al\geq2$.  If $a$ is not injective, then $\Lam_I\geq2$.  If $a$ is not surjective, then $\be\geq1$ and so $\al^\be\geq2$.  Thus, since $\rank(\S_\al)\leq2$, it follows that $\max\big( \al^\be, \Lam_I, \rank(\S_\al) \big)=\max( \al^\be, \Lam_I)$.~\epf

As with Remark \ref{rem:rank_Reg_PTXA}, by taking $a$ to be injective or surjective, we obtain as corollaries the ranks of the regular subsemigroups of $\T(X,A)$ and~$\T(Y,\si)$, respectively.  If $X$ is a finite set, and if $\emptyset\subsetneq A\subsetneq X$, then
\begin{align*}
\rank(\Reg(\T(X,A))) &= \begin{cases}
1 &\hspace{-.7cm}\text{if $|A|=1$}\\
1+|A|^{|X|-|A|} \phantom{1+\Lam_I} &\hspace{-.7cm}\text{if $|A|\geq2$.}
\end{cases}
\intertext{The above was proved in \cite[Theorem 3.6]{SS2013}.  If $Y$ is a finite non-empty set, and if $\si$ is an equivalence relation on $Y$, then}
\rank(\Reg(\T(Y,\si))) &= \begin{cases}
|Y| &\hspace{-.7cm}\text{if $|Y/\si|=1$}\\
1+\Lam_I \phantom{1+|A|^{|X|-|A|}} &\hspace{-.7cm}\text{if $|Y/\si|\geq2$.}
\end{cases}
\end{align*} 
The latter appears to be a new result.

\subsectiontitle{Idempotents and idempotent-generation}\label{sect:EaTXYa}

As in Proposition \ref{prop:EaPTXYa}, the idempotents of $\TXYa$ may be characterised and enumerated as follows:
\bit
\itemnit{i} $E_a(\TXYa) = \set{f\in\TXY}{(af)|_{\im(f)}=\id_{\im(f)}}$.
\itemnit{ii} If $|P^a|\geq\aleph_0$, then $|E_a(\TXYa)|=|P^a|$.
\itemnit{iii} If $|P^a|<\aleph_0$, then $\ds |E_a(\TXYa)| = \sum_{\mu=1}^\al\mu^{|X|-\mu}\sum_{J\sub I\atop|J|=\mu}\Lam_J$.
\eit
As with Theorems \ref{thm:EXYa} and \ref{thm:finite_EXYa_PT}, we may describe the idempotent-generated subsemigroup of the sandwich semigroup $\TXYa$; for this subsection only, we denote this subsemigroup by $\EXYa=\E_a(\TXYa)$.  To give this description, we use the $\T_A$ analogue of Proposition \ref{prop:PTBSB}, which comes from \cite{Gomes1987,Howie1978,Howie1966}, and states the following:
\bit
\itemnit{i} If $|\B|<\aleph_0$, then $\E(\T_\B) = \{\id_\B\} \cup (\T_\B\sm\S_\B)$, and $\rank(\E(\T_\B))=\idrank(\E(\T_\B))=\tbinom\al2+1$.
\itemnit{ii} If $|\B|\geq\aleph_0$, then $\E(\T_\B)=\{\id_\B\}\cup\set{f\in\T_\B\sm\S_\B}{\sh(f)<\aleph_0}\cup\set{f\in\T_\B}{\sh(f)=\col(f)=\defect(f)\geq\aleph_0}$, and $\rank(\E(\T_\B))=\idrank(\E(\T_\B))=|\T_\B|=2^{|\B|}$.
\eit

\begin{thm}\label{thm:EXYa_T}
We have $\EXYa=\E_a(\TXYa)=\E(\T_\B)\varphi^{-1}$.  If $\al<\aleph_0$, then $\EXYa=E_a(D_\al^a) \cup (P^a\sm D_\al^a)$.  Further,
\[
\rank(\EXYa)=\idrank(\EXYa)= \begin{cases}
|\EXYa|=|P^a| &\text{if $|P^a|\geq\aleph_0$}\\
\binom\al2+\max(\al^\be,\Lam_I) &\text{if $|P^a|<\aleph_0$.}
\end{cases}
\]
\end{thm}

Again, we may deduce formulae for the number of idempotents in the semigroups $\T(X,A)$ and $\T(Y,\si)$, and for the (idempotent) ranks of the idempotent-generated subsemigroups of $\T(X,A)$ and $\T(Y,\si)$, but we leave the details for the reader.

\subsectiontitle{The rank of a sandwich semigroup $\TXYa$}\label{sect:rank_TXYa}

As in Section \ref{sect:rank_PTXYa}, we may give formulae for the rank of a sandwich semigroup $\TXYa$.  Again, we eliminate some easy special cases:
\bit
\item If $|X|=1$, then $\TXYa$ is a right zero semigroup, and so $\rank(\TXYa)=|\TXYa|=|Y|$.
\item If $|Y|=1$, then $\rank(\TXYa)=|\TXYa|=1$.
\item If $|Y|\geq2$ and $|X|\geq\aleph_0$, then $\TXYa$ is uncountable, and so $\rank(\TXYa)=|\TXYa|=|Y|^{|X|}$.
\item If $|Y|>\aleph_0$, then $\TXYa$ is uncountable, so again $\rank(\TXYa)=|\TXYa|=|Y|^{|X|}$.
\eit
Thus, for the duration of Section \ref{sect:rank_TXYa}, we assume that $2\leq|X|<\aleph_0$ and $2\leq|Y|\leq\aleph_0$.  In the case that $\al=\rank(a)<\xi=\min(|X|,|Y|)$, we have the following (note that $\al<\xi$ forces $a$ to be non-injective and non-surjective); its proof is parallel to that of Theorem \ref{thm:rankPTXYa1}, moving through a similar series of preliminary results.

\begin{thm}%[cf.~Theorem \ref{thm:rankPTXYa1}]
\label{thm:rankTXYa1}
Suppose $\al<\xi$, and that $2\leq|X|<\aleph_0$ and $2\leq|Y|\leq\aleph_0$.  Then
\[
\epfreseq
\rank(\TXYa) = \sum_{\mu=\al+1}^\xi \mu! \tbinom{|Y|}{\mu}S(|X|,\mu).
\]
\end{thm}

For the $\al=\xi$ case, the proofs of Theorems \ref{thm:rankPTXYa2} and \ref{thm:rankPTXYa3} may easily be modified to yield the following.

\begin{thm}
%[cf.~Theorems \ref{thm:rankPTXYa2} and \ref{thm:rankPTXYa3}]
\label{thm:rankTXYa2}
%\leavevmode\newline
%~\vspace{-4truemm}
%\bit
\begin{itemize}
\itemit{i} If $\al=|Y|<|X|<\aleph_0$, then $\rank(\TXYa) = S(|X|,\al)$. 
\itemit{ii} If $\al=|X|<|Y|\leq\aleph_0$, then $\rank(\TXYa) = \binom{|Y|}\al$. \epfres
\end{itemize}
\end{thm}

Theorem \ref{thm:rankTXYa2} yields results concerning the ranks of finite $\T(X,A)$ and $\T(Y,\si)$:
\bit
\itemnit{i} If $X$ is a finite set, and $\emptyset\not=A\subsetneq X$, then $\rank(\T(X,A))=S(|X|,|A|)$.  This was originally proved in \cite[Theorem 2.3]{FS2014}.
\itemnit{ii} If $Y$ is finite, and $\si$ is a non-diagonal equivalence on $Y$, then $\rank(\T(Y,\si)) = \tbinom{|Y|}{|Y/\si|}$.  This appears to be a new result.
\eit

\subsectiontitle{Egg-box diagrams}\label{sect:eggbox_T}

Figures \ref{fig:T_1} and \ref{fig:T_2} give egg-box diagrams (in the sense of Section \ref{sect:eggbox}) for various sandwich semigroups $\TXYa$ and their regular subsemigroups.

\begin{figure}[ht]
\begin{center}
\includegraphics[height=6.1cm]{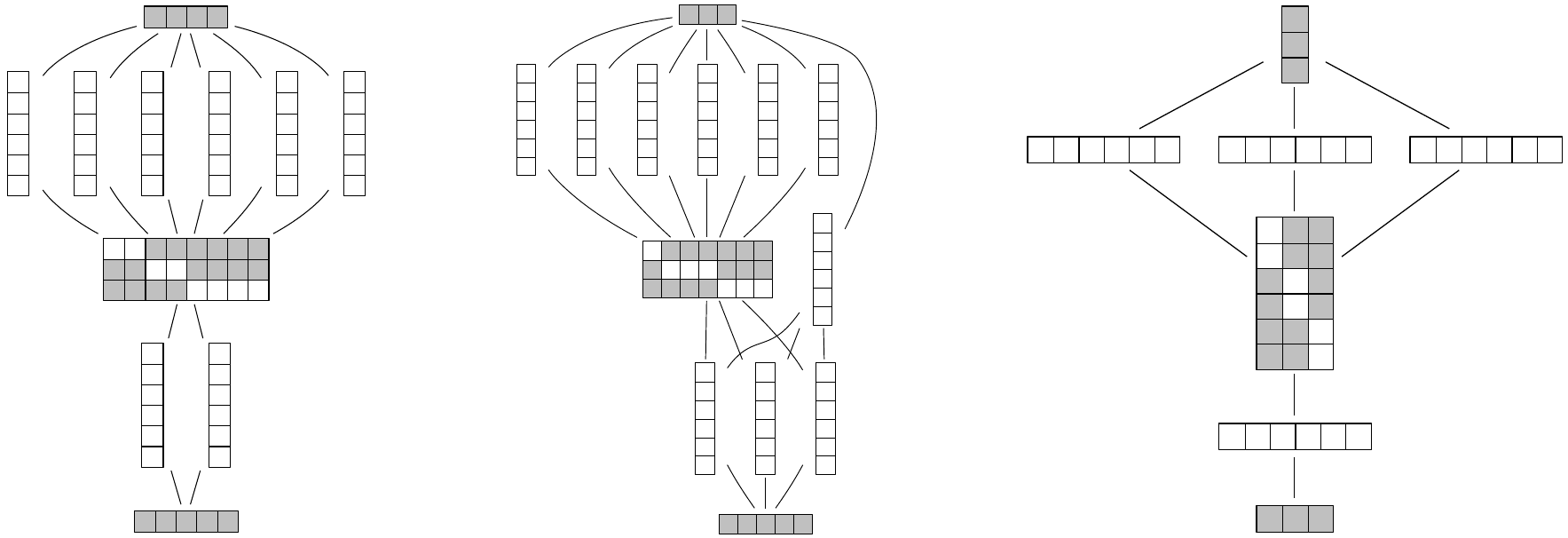} 
%\vspace{-5mm}
\caption[blah]{Left to right: sandwich semigroups $\T_{35}^a$, $\T_{35}^b$, $\T_{43}^c$, where $a=\trans{1&2&3&4&5}{1&2&2&3&3}\in\T_{53}$, $b=\trans{1&2&3&4&5}{1&2&3&3&3}\in\T_{53}$ and $c=\trans{1&2&3}{1&2&3}\in\T_{34}$.  Note that $a,b$ are surjective (but not injective), and $c$ is injective (but not surjective).}
\label{fig:T_1}
\end{center}
\end{figure}

\vspace{-6mm}

\begin{figure}[ht]
\begin{center}
\includegraphics[height=6.1cm]{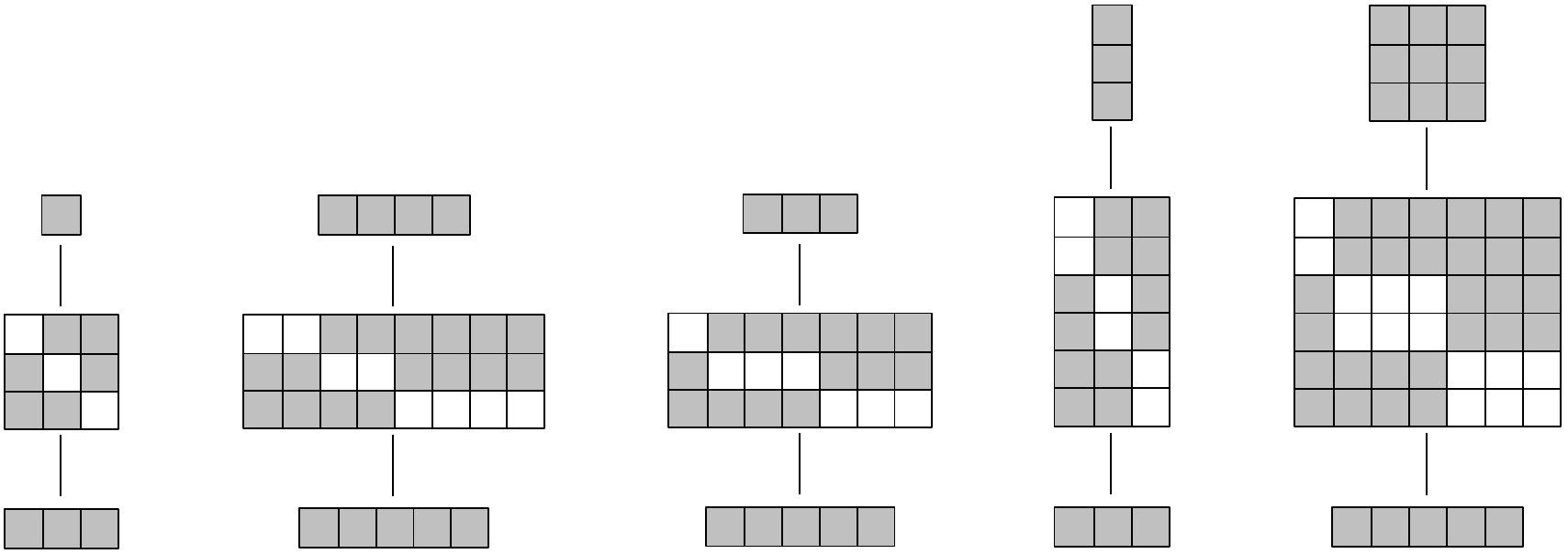} 
%\vspace{-5mm}
\caption[blah]{Left to right: the full transformation semigroup $\T_3$, and the regular sandwich semigroups $\Reg(\T_{35}^a)$, $\Reg(\T_{35}^b)$, $\Reg(\T_{43}^c)$ and $\Reg(\T_{45}^d)$, where $a,b,c$ are as in Figure \ref{fig:T_1}, and $d=\trans{1&2&3&4&5}{1&2&3&3&3}\in\T_{54}$. The four regular sandwich semigroups are all inflations of $\T_3$ (in the sense of \cite[Remark 2.15]{Sandwiches1}).}
\label{fig:T_2}
\end{center}
\end{figure}

\vspace{-4mm}

\sectiontitle{The category $\I$}\label{sect:I}

Recall that the category $\I = \bigset{(A,f,B)}{A,B\in\Set,\ f\in\bI_{AB}}$ is the (regular monoidal) partial subsemigroup of the category $\PT$ consisting of all \emph{injective} partial transformations.  Again, we may adapt the approach of Section \ref{sect:PT} to study the category $\I$.  However, since $\I$ is an \emph{inverse category} (to be explained and proved shortly), almost all of the results of \cite[Sections 2 and 3]{Sandwiches1}, concerning regular subsemigroups $\Reg(\IXYa)$ of the sandwich semigroups $\IXYa$, reduce to known results concerning symmetric inverse monoids $\I_A$.

\begin{lemma}
Let $X,Y\in\Set$, and let $f\in\bI_{XY}$ and $g\in\bI_{YX}$.  Then $f=fgf$ and $g=gfg$ both hold if and only if $g=f^{-1}$ is the inverse mapping of $f$.  
\end{lemma}

\pf From $f=fgf$ and Lemma \ref{lem:PT_prelim}(i), we obtain $\dom(f)=\dom(fgf)\sub\dom(fg)\sub\dom(f)$, so $\dom(fg)=\dom(f)$.  It follows from the same lemma that $\im(f)\sub\dom(g)$.  Similarly, considering images of $f=fgf$, and applying Lemma \ref{lem:PT_prelim}(ii), we obtain $\dom(f)\sub\im(g)$; note that the latter is equivalent to the assertion that $\im(g)$ saturates $\ker(f)$.  Reversing the roles of $f,g$, and using $g=gfg$, we obtain 
the reverse inclusions, and conclude that
%
%$\im(f)\supseteq\dom(g)$ and $\dom(f)\supseteq\im(g)$.  Thus, 
%
$\im(f)=\dom(g)$ and $\dom(f)=\im(g)$.  So we may write $f=\binom{f_i}{g_i}_{i\in I}$ and $g=\binom{g_i}{f_{ih}}$ for some permutation $h\in\S_I$.  But $\binom{f_i}{g_i}=f=fgf=\binom{f_i}{g_{ih}}$ gives $h=\id_I$, and so $g=\binom{g_i}{f_i}=f^{-1}$. \epf

It follows that every element of $\I$ is \emph{uniquely regular}, and so \emph{uniquely sandwich-regular}, using the terminology of \cite[Section 4]{Sandwiches1}.  In other words, $\I$ is an \emph{inverse category}, as defined in \cite{Kastl1979} and \cite[Section 2.3.2]{CL2002}.

As in Section \ref{sect:Green_T}, we may deduce descriptions of Green's relations and preorders on $\I$ from Proposition~\ref{prop:GreenPT}.  We leave it to the reader to supply the details (but note that statements concerning kernels now become redundant).
Again, the~${\J}={\D}$-classes of $\I_{AB}$ are the sets
\[
D_\mu = D_\mu(\I_{AB}) = \set{(A,f,B)}{f\in\bI_{AB},\ \rank(f)=\mu} \qquad\text{for each cardinal \ $0\leq\mu\leq\min(|A|,|B|)$,}
\]
%(Note that the minimum rank of a full transformation is $1$.)
and again these form a chain.  If $|A|=\al$ and $|B|=\be$, then
\[
|D_\mu/{\R}| = \tbinom\al\mu \COMMA
|D_\mu/{\L}| = \tbinom\be\mu \COMMA
|D_\mu/{\H}| = \tbinom\al\mu\tbinom\be\mu \COMMA
|D_\mu| = \mu! \tbinom\al\mu\tbinom\be\mu ,
\]
and each $\H$-class in $D_\mu$ has size $\mu!$.  This yields the standard formula $|\I_{AB}|=\sum_{\mu=0}^{\min(\al,\be)}\mu! \tbinom\al\mu\tbinom\be\mu$; to the authors' knowledge, no closed formula for $|\I_{AB}|$ exists. Again, $\R$- and/or $\L$-stable elements of $\I$ are easily described.

In order to discuss sandwich semigroups in $\I$, we fix some sets $X,Y\in\Set$, and as usual identify $\I_{XY}$ with~$\bI_{XY}$, and so on.  Fix $a=\binom{b_i}{a_i}_{i\in I}\in\I_{YX}$ and $b=a^{-1}=\binom{a_i}{b_i}\in\I_{XY}$.  Note that $f\mt f^{-1}$ determines an anti-isomorphism $\I_{XY}^a\to\I_{YX}^b$.  Write
\[
A=\im(a)=\set{a_i}{i\in I} \AND B=\dom(a)=\set{b_i}{i\in I}.
\]
Then we have
%$P_1^a=\set{f\in\I_{XY}}{\im(f)\sub B}$, $P_2^a=\set{f\in\I_{XY}}{\dom(f)\sub A}$ and $P_3^a=\set{f\in\I_{XY}}{\rank(afa)=\rank(f)}$, so that
\[
P_1^a=\set{f\in\I_{XY}}{\im(f)\sub B} \COMMa
P_2^a=\set{f\in\I_{XY}}{\dom(f)\sub A} \COMMa
%P^a=\set{f\in\I_{XY}}{\im(f)\sub A} \COMMA
P_3^a=\set{f\in\I_{XY}}{\rank(afa)=\rank(f)}.
\]
%and so $\Reg(\IXYa)=P^a=\set{f\in\I_{XY}}{\dom(f)\sub A,\ \im(f)\sub B}$.
These yield descriptions of Green's relations on $\IXYa$ (which we will not explicitly state), and show that the regular subsemigroup of $\IXYa$ is
\[
\Reg(\IXYa)=P^a=\set{f\in\I_{XY}}{\dom(f)\sub A,\ \im(f)\sub B}.
\]
The commutative diagrams \cite[diagrams (2.1) and (2.7)]{Sandwiches1} again lead to relationships between the sandwich semigroup~$\IXYa$ and certain (non-sandwich) transformation semigroups: namely,
\[
\IXY a = \I_X(ba) = \set{f\in\I_X}{\im(f)\sub A} 
\AND
a\IXY = (ab)\I_Y = \set{f\in\I_Y}{\dom(f)\sub B}.
\]
The former semigroup has been studied extensively in the literature (see for example \cite{FS2014}), and is usually denoted $\I(X,A)$; the latter semigroup is anti-isomorphic to $\I(Y,B)$, and we will denote it by $\I(Y,B)^*$.  Again, the regular monoid $(a\IXY a,\starb)\sub\I_{YX}^b$ is isomorphic to $ba\IXY b=(ba)\I_X(ba)$, the local monoid of $\I_X$ with respect to the idempotent $ba=\binom{a_i}{a_i}\in\I_X$; the latter submonoid is isomorphic to, and will be identified with, $\I_A$.  Thus, we obtain the diagrams:
\begin{equation}\label{eq:CD_I_2}
\includegraphics{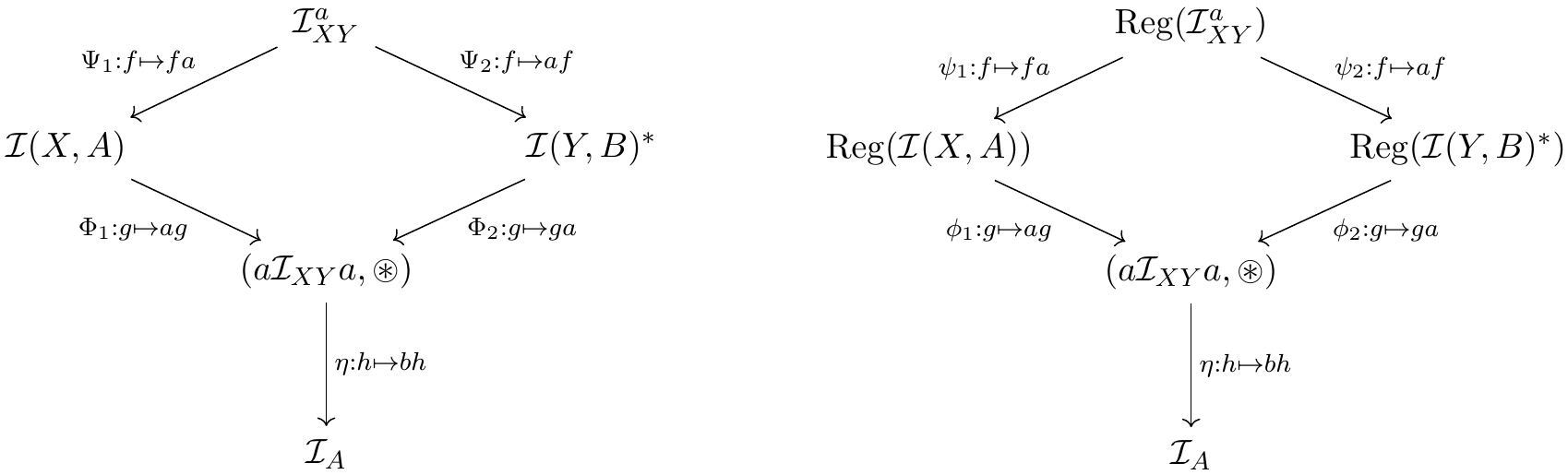}
%\begin{tikzcd} %[column sep=small]
%~ & \IXYa \arrow[swap]{dl}{\Psi_1:f\mt fa} \arrow{dr}{\Psi_2:f\mt af} & \\
%\I(X,\B) \arrow[swap]{dr}{\Phi_1:g\mt ag} & & \I(Y,B)^* \arrow{dl}{\Phi_2:g\mt ga}\\
%& (a\IXY a,\starb)\arrow{dd}{\eta:h\mt bh} & \\
%\\
%& \I_\B & 
%\end{tikzcd}
%\qquad \qquad
%\begin{tikzcd} [column sep=small]
%~ & \Reg(\IXYa) \arrow[swap]{dl}{\psi_1:f\mt fa} \arrow{dr}{\psi_2:f\mt af} & \\
%\Reg(\I(X,\B)) \arrow[swap]{dr}{\phi_1:g\mt ag} & & \Reg(\I(Y,B)^*) \arrow{dl}{\phi_2:g\mt ga}\\
%& (a\IXY a,\starb)\arrow{dd}{\eta:h\mt bh}  & \\
%\\
%& \I_\B & 
%\end{tikzcd}
\end{equation}
Since $a$ (indeed, every element of $\I$) is uniquely sandwich-regular, it follows from \cite[Proposition 4.1]{Sandwiches1} that all maps in the right-hand diagram of \eqref{eq:CD_I_2} are isomorphisms.  Consequently, 
\[
\Reg(\IXYa) \cong \Reg(\I(X,A)) \cong \Reg(\I(Y,B)^*) \cong \I_A.
\]
(In fact, it was shown in \cite[Theorem 3.1]{FS2014} that $\Reg(\I(X,A))=\I_A$, giving also $\Reg(\I(Y,B)^*)=\I_B\cong\I_A$.)  Consequently, any result concerning $\Reg(\IXYa)$ is equivalent to a result concerning $\I_A$; cf.~\cite[Remark 4.2]{Sandwiches1}.  Since the symmetric inverse monoids $\I_A$ are well understood, we will not state any such results, but instead refer the reader to the monographs \cite{Lipscombe1996,Lawson1998,Petrich1984}.  (Note that since $\I_A$ is an inverse semigroup, its idempotent-generated subsemigroup is precisely its semilattice of idempotents, and this is isomorphic to the power set of~$A$ under intersection.)

Nevertheless, the problem of calculating the rank of a sandwich semigroup $\IXYa$ itself is not covered by the general theory developed in \cite{Sandwiches1}, although the method of Section \ref{sect:rank_PTXYa} (of the current paper) may be easily adapted to yield the next result; since $\IXYa$ is anti-isomorphic to $\I_{YX}^b$, as noted above, we may assume that~$|Y|\leq|X|$ in the statement.  For part (ii), note that $\IXYa$ has a maximum $\J^a$-class if $\al=\min(|X|,|Y|)$, and this is a group isomorphic to $\S_\al$ (cf.~Theorem \ref{thm:rankPTXYa3} and Lemma \ref{lem:rank_PTXYa_lower}(iv)).  %Since $\IXYa$ is anti-isomorphic to $\I_{YX}^{a^{-1}}$, as noted above, we may assume that $|Y|\leq|X|$ in the next result.

\begin{thm}
\label{thm:rankIXYa1}
\leavevmode
\vspace{-5truemm}
\begin{itemize}
\itemit{i} If $\al<|Y|\leq|X|\leq\aleph_0$, then $\ds\rank(\IXYa) = \sum_{\mu=\al+1}^{|Y|} \mu! \tbinom{|X|}{\mu}\tbinom{|Y|}{\mu}$.
\itemit{ii} If $\al=|Y|<|X|\leq\aleph_0$, then $\rank(\IXYa) = \binom{|X|}{|Y|} + \begin{cases}
0 &\text{if $\al\leq2$}\\
1 &\text{if $\al\geq3$.}
\end{cases}$
\epfres
\end{itemize}
\end{thm}

For values of $\al,|X|,|Y|$ outside of the stated ranges, the value of $\rank(\IXYa)$ is easy to calculate.  
Egg-box diagrams for some sandwich semigroups $\IXYa$ are given in Figures \ref{fig:I_1} and \ref{fig:I_2}.

\begin{figure}[ht]
\begin{center}
\includegraphics[width=\textwidth]{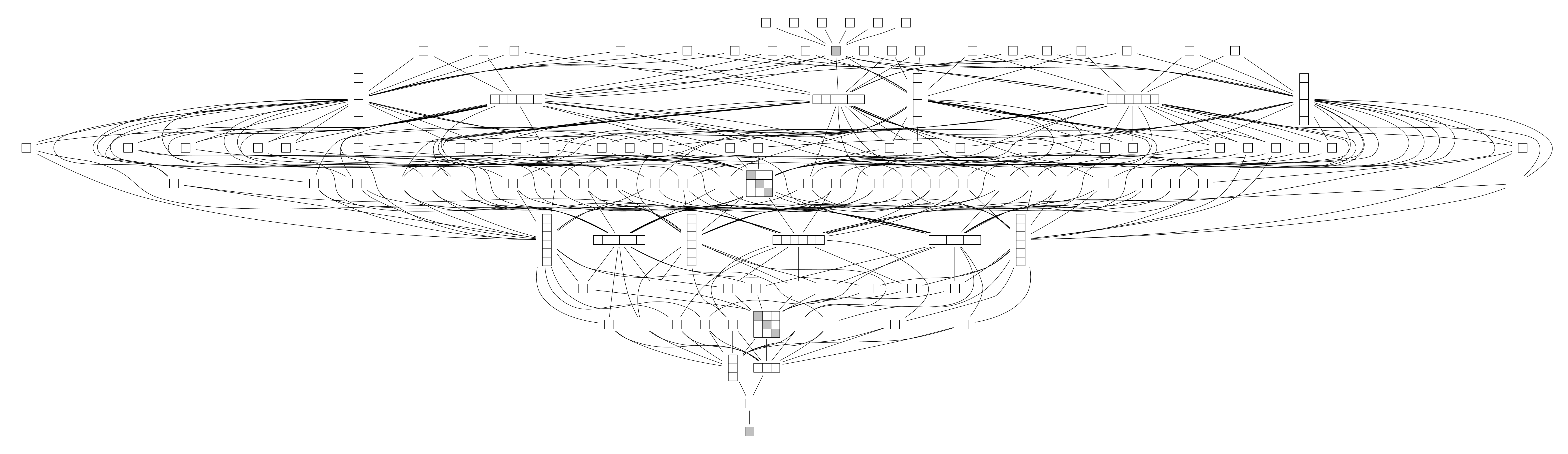} \\~\\
\vspace{-5mm}
\caption[blah]{Egg-box diagram of the sandwich semigroup $\I_{44}^a$, where $a=\trans{1&2&3&4}{1&2&3&-}\in\I_{4}$. Note that $a$ is neither full nor surjective.}
\label{fig:I_1}
\end{center}
\end{figure}

\vspace{-5mm}

\begin{figure}[ht]
\begin{center}
\includegraphics[height=7cm]{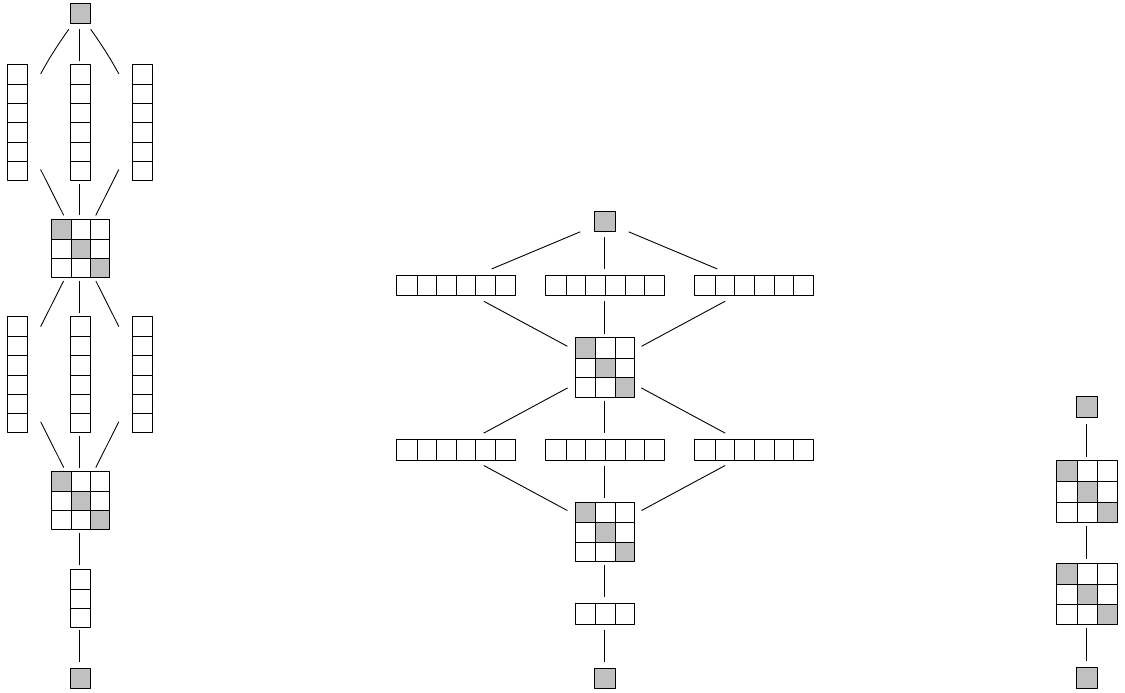}
%\vspace{-5mm}
\caption[blah]{Left to right: egg-box diagrams of the sandwich semigroups $\I_{43}^b$ and $\I_{34}^c$, and the symmetric inverse monoid $\I_3$, where $b=\trans{1&2&3}{1&2&3}\in\I_{34}$ and $c=\trans{1&2&3&4}{1&2&3&-}\in\I_{43}$. Note that $b$ is full but not surjective, while $c$ is surjective but not full.  Note also that $\Reg(\I_{44}^a)$, $\Reg(\I_{43}^b)$ and $\Reg(\I_{34}^c)$ are all isomorphic to $\I_3$, where $a$ is as in Figure \ref{fig:I_1}.}
\label{fig:I_2}
\end{center}
\end{figure}

\section*{Acknowledgements}

This work was initiated during visits of the third author to Chiang Mai in 2015, and to Novi Sad in 2016; he thanks these institutions for their generous hospitality.  
The first and second authors are supported by Grant Nos.~174019 and 174018, respectively, of the Ministry of Education, Science, and Technological Development of the Republic of Serbia.
The first and third authors would like to thank \emph{Macchiato Fax} for providing us with many a \emph{Klub sendvi\v{c}}.
Finally, we thank the referee for his/her careful reading of a lengthy pair of articles, and for a number of helpful suggestions.

\footnotesize
\def\bibspacing{-1.1pt}
\bibliography{biblio}

\begin{thebibliography}{10}

\bibitem{CL2002}
J.~R.~B. Cockett and Stephen Lack.
\newblock Restriction categories. {I}. {C}ategories of partial maps.
\newblock {\em Theoret. Comput. Sci.}, 270(1-2):223--259, 2002.

\bibitem{DEvariants}
Igor Dolinka and James East.
\newblock Variants of finite full transformation semigroups.
\newblock {\em Internat. J. Algebra Comput.}, 25(8):1187--1222, 2015.

\bibitem{DElinear}
Igor Dolinka and James East.
\newblock Semigroups of rectangular matrices under a sandwich operation.
\newblock {\em Semigroup Forum}, to appear, {\tt arXiv:1503.03139}.

\bibitem{Sandwiches1}
Igor Dolinka, Ivana~\DJ ur\dj ev, James East, Preeyanuch Honyam, Kritsada
  Sangkhanan, Jintana Sanwong, and Worachead Sommanee.
\newblock Sandwich semigroups in locally small categories {I}: {F}oundations.
\newblock {\em Preprint}, 2017, {\tt arXiv:1710.01890}.

\bibitem{EP1970}
A.~E. Evseev and N.~E. Podran.
\newblock Transformation semigroups of a finite set that are generated by
  idempotents with given projective characteristics.
\newblock {\em Izv. Vys\v s. U\v cebn. Zaved. Matematika}, 1970(12
  (103)):30--36, 1970.

\bibitem{FS2014}
V{\'{\i}}tor~H. Fernandes and Jintana Sanwong.
\newblock On the ranks of semigroups of transformations on a finite set with
  restricted range.
\newblock {\em Algebra Colloq.}, 21(3):497--510, 2014.

\bibitem{GMbook}
Olexandr Ganyushkin and Volodymyr Mazorchuk.
\newblock {\em Classical finite transformation semigroups, an introduction},
  volume~9 of {\em Algebra and Applications}.
\newblock Springer-Verlag London, Ltd., London, 2009.

\bibitem{Garba1990}
G.~U. Garba.
\newblock Idempotents in partial transformation semigroups.
\newblock {\em Proc. Roy. Soc. Edinburgh Sect. A}, 116(3-4):359--366, 1990.

\bibitem{Gomes1987}
Gracinda Gomes and John~M. Howie.
\newblock On the ranks of certain finite semigroups of transformations.
\newblock {\em Math. Proc. Cambridge Philos. Soc.}, 101(3):395--403, 1987.

\bibitem{Howie1966}
J.~M. Howie.
\newblock The subsemigroup generated by the idempotents of a full
  transformation semigroup.
\newblock {\em J. London Math. Soc.}, 41:707--716, 1966.

\bibitem{Howie1978}
J.~M. Howie.
\newblock Idempotent generators in finite full transformation semigroups.
\newblock {\em Proc. Roy. Soc. Edinburgh Sect. A}, 81(3-4):317--323, 1978.

\bibitem{Howie}
John~M. Howie.
\newblock {\em Fundamentals of semigroup theory}, volume~12 of {\em London
  Mathematical Society Monographs. New Series}.
\newblock The Clarendon Press, Oxford University Press, New York, 1995.
\newblock Oxford Science Publications.

\bibitem{Jech2003}
Thomas Jech.
\newblock {\em Set theory}.
\newblock Springer Monographs in Mathematics. Springer-Verlag, Berlin, 2003.
\newblock The third millennium edition, revised and expanded.

\bibitem{Kastl1979}
J.~Kastl.
\newblock Inverse categories.
\newblock In {\em Algebraische {M}odelle, {K}ategorien und {G}ruppoide},
  volume~7 of {\em Stud. Algebra Anwendungen}, pages 51--60. Akademie-Verlag,
  Berlin, 1979.

\bibitem{Lawson1998}
Mark~V. Lawson.
\newblock {\em Inverse semigroups}.
\newblock World Scientific Publishing Co., Inc., River Edge, NJ, 1998.
\newblock The theory of partial symmetries.

\bibitem{Lipscombe1996}
Stephen Lipscomb.
\newblock {\em Symmetric inverse semigroups}, volume~46 of {\em Mathematical
  Surveys and Monographs}.
\newblock American Mathematical Society, Providence, RI, 1996.

\bibitem{MacLane1998}
Saunders Mac~Lane.
\newblock {\em Categories for the working mathematician}, volume~5 of {\em
  Graduate Texts in Mathematics}.
\newblock Springer-Verlag, New York, second edition, 1998.

\bibitem{MS1975}
K.~D. Magill, Jr. and S.~Subbiah.
\newblock Green's relations for regular elements of sandwich semigroups. {I}.
  {G}eneral results.
\newblock {\em Proc. London Math. Soc. (3)}, 31(2):194--210, 1975.

\bibitem{Magill1967}
Kenneth~D. Magill, Jr.
\newblock Semigroup structures for families of functions. {I}. {S}ome
  homomorphism theorems.
\newblock {\em J. Austral. Math. Soc.}, 7:81--94, 1967.

\bibitem{MGS2013}
Suzana Mendes-Gon{\c{c}}alves and R.~P. Sullivan.
\newblock Regular elements and {G}reen's relations in generalized
  transformation semigroups.
\newblock {\em Asian-Eur. J. Math.}, 6(1):1350006, 11, 2013.

\bibitem{MGS2010}
Suzana Mendes-Gon{\c{c}}alves and Robert~P. Sullivan.
\newblock Semigroups of transformations restricted by an equivalence.
\newblock {\em Cent. Eur. J. Math.}, 8(6):1120--1131, 2010.

\bibitem{GAP}
J.~D. Mitchell et~al.
\newblock {\em Semigroups - GAP package, Version 3.0.5}, Aug 2017.

\bibitem{Petrich1984}
Mario Petrich.
\newblock {\em Inverse semigroups}.
\newblock Pure and Applied Mathematics (New York). John Wiley \& Sons, Inc.,
  New York, 1984.
\newblock A Wiley-Interscience Publication.

\bibitem{SS2008}
Jintana Sanwong and Worachead Sommanee.
\newblock Regularity and {G}reen's relations on a semigroup of transformations
  with restricted range.
\newblock {\em Int. J. Math. Math. Sci.}, 2008, Art. ID 794013, 11 pp.

\bibitem{SS2013}
Worachead Sommanee and Jintana Sanwong.
\newblock Rank and idempotent rank of finite full transformation semigroups
  with restricted range.
\newblock {\em Semigroup Forum}, 87(1):230--242, 2013.

\bibitem{Symons1975}
J.~S.~V. Symons.
\newblock Some results concerning a transformation semigroup.
\newblock {\em J. Austral. Math. Soc.}, 19(4):413--425, 1975.

\bibitem{Tirasupa1979}
Yupaporn Tirasupa.
\newblock Factorizable transformation semigroups.
\newblock {\em Semigroup Forum}, 18(1):15--19, 1979.

\end{thebibliography}
\bibliographystyle{plain}
\end{document}